\documentclass[review, sort&compress, numbers]{elsarticle}
\usepackage[fleqn]{amsmath}
\usepackage{graphicx}
\usepackage{lineno,hyperref}
\usepackage{latexsym}
\usepackage{amssymb}
\usepackage{amsmath}
\usepackage[center]{subfigure}
\usepackage{float}
\usepackage{color}
\usepackage{geometry}
\usepackage{diagbox}
\usepackage[section]{placeins}

\usepackage{bm}
\usepackage{algorithm}
\usepackage{algorithmicx}
\usepackage{algpseudocode}
\usepackage{amsmath}

\usepackage{float}
\usepackage{lipsum}
\makeatletter
\usepackage{array}
\usepackage{multirow}
\usepackage{longtable}
\usepackage{rotating}
\newtheorem{definition}{Definition}

\modulolinenumbers[5]

\newtheorem{theorem}{Theorem}

\graphicspath{{./figs/},{./logo/}}

\journal{}

\begin{document}
\begin{frontmatter}
\title{Two New Low Rank Tensor Completion Methods Based on Sum Nuclear Norm}

% Group authors per affiliation:
%\author{Hongbing Zhang,
%	Xinyi Liu, Hongtao Fan, Yajing Li, Yinlin Ye
%	\thanks{This work was supported by the National Natural Science Foundation of
%		China (Nos. 11701456, 11801452, 11571004), Fundamental Research Project of Natural Science in Shaanxi Province General Project (Youth) (Nos. 2019JQ-415, 2019JQ-196), the Initial Foundation for Scientific Research of Northwest A\&F University (Nos.2452017219, 2452018017), and Innovation and Entrepreneurship Training Program for College Students of Shaanxi Province (S201910712132).}
%	\thanks{H. Zhang, X. Liu, H. Fan, Y. Li and Y. Ye are with the College of Science, Northwest A\&F University, Yangling, Shaanxi 712100, China(e-mail: zhanghb@nwafu.edu.cn; Lxy6x1@163.com; fanht17@nwafu.edu.cn; hliyajing@163.com; 13314910376@163.com).}
\author[mymainaddress]{Hongbing~Zhang}
\author[mymainaddress]{Xinyi~Liu}
\author[mymainaddress]{Hongtao~Fan\corref{mycorrespondingauthor}}
\cortext[mycorrespondingauthor]{Corresponding author}
\ead{fanht17@nwafu.edu.cn}
\author[mymainaddress]{Yajing~Li}
\author[mymainaddress]{Yinlin~Ye}
\author[mysecondaryaddress]{Xinyun~Zhu}
\address[mymainaddress]{College of Science, Northwest A\&F University, Yangling, Shaanxi 712100, PR China}
\address[mysecondaryaddress]{Department of Mathematics, University of Texas of the Permian Basin, Odessa, TX 79762, USA}

\begin{abstract}
The low rank tensor completion (LRTC) problem has attracted great attention in computer vision and signal processing. How to acquire high quality image recovery effect is still an urgent task to be solved at present. This paper proposes a new tensor $L_{2,1}$ norm minimization model (TLNM) that integrates sum nuclear norm (SNN) method, differing from the classical tensor nuclear norm (TNN)-based tensor completion method, with $L_{2,1}$ norm and Qatar Riyal decomposition for solving the LRTC problem. To improve the utilization rate of the local prior information of the image, a total variation (TV) regularization term is introduced, resulting in a new class of tensor $L_{2,1}$ norm minimization with total variation model (TLNMTV). Both proposed models are convex and therefore have global optimal solutions. Moreover, we adopt the Alternating Direction Multiplier Method (ADMM) to obtain the closed-form solution of each variable, thus ensuring the feasibility of the algorithm. Numerical experiments show that the two proposed algorithms are convergent and outperform compared methods. In particular, our method significantly outperforms the contrastive methods when the sampling rate of hyperspectral images is 2.5\%.
\end{abstract}


\end{frontmatter}

\section{Introduction}
The low rank tensor completion (LRTC) problem, especially when dealing with extremely high-dimensional data, such as appearing in color image and video processing \cite{1472018163}, \cite{3762018t397}, magnetic resonance image \cite{341202131}, \cite{1242020783}, hyperspectral image \cite{XUE2019109}, \cite{7467446}, \cite{8657368}, pattern recognition \cite{6213108}, \cite{6862897} face modeling and analysis \cite{6995956} and other fields, has received extensive attention. The fact that the imaging of visual data is based on poor acquisition conditions or severe data corruption during transmission, results in incomplete or severely corrupted data being acquired. This will facilitate LRTC \cite{8941238} to become crucially important for collecting data, particularly how to utilize the internal structural information between collected observations and missing data to achieve LRTC problem.
%due to various factors, such as poor imaging acquisition conditions for visual data, or severe data corruption during transmission, the acquired data is incomplete or severely damaged. This will promote the research of LRTC on the aquired data, especially how to utilize the internal structure information between the acquired observation data and the missing data to realize LRTC becomes particularly important.

%However, the obtained different ranks have their own shortcomings. For examples, the calculation of the rank of tensor CP is NP-hard; the determination of Tucker rank requires the operation of the mode-$n$ unfolding and folding of one tensor , which will cause the curse of dimensionality; the multi-rank and tube rank \cite{124120101253}, \cite{6412011658} based on t-SVD decomposition induced by tensor product only apply to third-order tensors due to the current limitations of tensor product.

Mathematically, the above LRTC problems can be can be characterized as follows:
\begin{eqnarray}
\min_{\mathcal{X}}\quad rank(\mathcal{X}), \quad s.t.\quad \mathcal{X}_{\Omega}=\mathcal{T}_{\Omega}
\end{eqnarray}
where $\mathcal{T}$ and $\mathcal{X}$ represent the incomplete and complete tensors,
respectively, $\Omega$ denotes an index set of observed data
of $\mathcal{T}$, $rank(\mathcal{X})$ implies tensor rank of any tensor $\mathcal{X}$, conditions to be satisfied $\mathcal{X}_{\Omega}=\mathcal{T}_{\Omega}$ mean that the elements at the same positions in two sets $\mathcal{X}$ and $\mathcal{T}$ remain equal according to index set $\Omega$. 

It is well-known, on the one hand, that the definition of tensor rank is not unique, and its various definitions of tensor rank have been proposed according to different ways of tensor decomposition, such as CANDECOMP/PARAFAC(CP) rank \cite{ACAR201141}, \cite{TICHAVSKY2017362} generated by CP decomposition \cite{43413413431}, Tucker-rank established by Tucker decomposition \cite{1996279311}, \cite{7488247}, and tube rank and multi-rank based on tensor products \cite{21321312}, \cite{6909886}, etc. However, specific features, such as the calculation of the rank of a tensor CP is NP-hard \cite{NPhard}, unfolding and folding of mode-$n$ on a tensor in Tucker rank \cite{3153564}, \cite{7460200} lead to dimensional disaster, and the limitation of the tensor product definition results in the multi-rank and tube rank \cite{124120101253}, \cite{6412011658} based on t-SVD decomposition are only applicable to third-order tensors, indicate that the definitions of different tensor rank show their respective inadequacies. The other half is that tensor rank is a non-convex function.

In this regard, Liu et al. \cite{6138863} turned to the convex envelope nuclear norm of the rank function to approximate it, and proposed a tensor nuclear norm minimization method to solve such problems, specifically, for image recovery problems

\begin{eqnarray}
\min_{\mathcal{X}} \|\mathcal{X}\|_{\ast},\quad s.t.\quad \mathcal{X}_{\Omega}=\mathcal{T}_{\Omega} \label{tensor1}
\end{eqnarray}
where $\|\mathcal{X}\|_{\ast}=\sum_{n=1}^{N}\alpha_{n}\|\mathbf{X}_{(n)}\|_{\ast}$ with $\alpha_{n}\geq0$ and $\sum_{n=1}^{N}\alpha_{n}=1$. 

To further enhance, despite the model (\ref{tensor1}) has ideal effect in global information recovery, the performance of the model for local information extraction, especially in processing high-dimensional image recovery, by incorporating local piece-wise smoothness prior, and proposed a large class of Total Variation (TV) minimization methods \cite{20062434113} and applied to image processing and pattern recognition \cite{7488247},  \cite{7089247}, \cite{8454775}, \cite{MADATHIL2018376}. In fact, for LRTC problems, TV terms are usually included in a low-rank (LR) framework to jointly represent local piecewise continuity and global LR structure along different dimensions \cite{10235753224}, \cite{6253256}, \cite{10.5555/3298483.3298558}, \cite{7298983}, which can greatly improve the local information effect of the model. Its specific modeling is as follows:

%Problem (\ref{tensor1}) has a good effect on global information recovery, but the local information extraction performance is not good enough to meet the requirements of practical image recovery applications, especially for high-dimensional image recovery. In view of this, a local piece-wise smoothness prior is employed to achieve the regularization of the image recovery, and a large class of total variation (TV) minimization methods is proposed and applied to image processing and pattern recognitionn. Generally, for the LRTC problem, the TV term is usually included in the low rank (LR) framework to jointly represent the local piece-wise continuity and global LR structure along different dimensions, thus greatly improving the local information effect of the model. Such specific model is shown as follows:
\begin{eqnarray}
\min_{\mathcal{X}} \|\mathcal{X}\|_{\ast}+\lambda\|\mathcal{X}\|_{TV},\quad s.t.\quad \mathcal{X}_{\Omega}=\mathcal{T}_{\Omega} \label{tensor22}
\end{eqnarray}
where $\|\mathcal{X}\|_{TV}$ is TV regularization term, and $\lambda>0$ is a parameter used to balance the relationship between the nuclear norm and the regularization term. 

Recently, Liu et al. \cite{8421060} proposed an approximate SVD computation technique based on QR decomposition (CSVD-QR), and integrated L2,1norm to construct a fast and accurate matrix completion method. Inspired by the high performance of this method, Zheng et al. \cite{2021108240} designed a numerical method that approximates t-SVD combined with Qatari Riyal Decomposition (CTSVD-QR), which further generalizes the matrix-form QR method to the tensor case. Through in-depth comparative analysis and thorough experimental exploration, it is found that although the speed of the algorithm has been significantly improved, compared with the classical tensor nuclear norm (TNN) method, the performance of the algorithm in terms of image restoration effect or quality is still less ideal. In view of this, to improve and develop this method, we consider integrating tensor SNN, not TNN, with QR decomposition and $L_{2,1}$ norm to design and build a more efficient and feasible method, namely, tensor $L_{2,1}$ norm minimization model (TLNM) . In order to further improve the effect of image restoration, that is, to better protect the local detail information while better representing the global information of the image, a TV regularization term is introduced and then a class of tensor $L_{2,1}$ norm minimizes the total variation model is proposed. To sum up, our main contributions are as follows:

Firstly, we adopt a new method different from TNN, that is, the integration of SNN and CSVD-QR method, to improve and enhance the restoration effect of the combination of LRTC and CSVD-QR method. For the sake of fairness, we still choose the $L_{2,1}$ norm instead of the classical nuclear norm in the objective function as in the recently published papers \cite{8421060} and \cite{2021108240}. In the meantime, by introducing the TV regularization term into the model constructed above, a new method named tensor $L_ {2, 1} $ norm minimization with TV (TLNMTV) is established to better preserve some details of the image, thereby improving the utilization of the local prior information of the image.

%Firstly, to improve the recovery effect of the combination of LRTC and CSVD-QR method, we adopt classical method different from TNN, namely, the combination of the SNN and CSVD-QR methods. For fairness, we still replace the nuclear norm with the $L_{2,1}$ norm to ensure the consistency with the previous paper \cite{8421060}, \cite{2021108240}. Total variation makes some details of the image better retained by improving the utilization rate of the local prior information of the image. In this regard, we propose the TLNMTV method by introducing the TV regularization term.

Secondly, for the above two established LRTC models, namely TLNM and TLMNTV, two efficient alternating direction multiplier (ADMM) algorithms \cite{683201156}, \cite{6122011620} are designed to solve them respectively. The closed-form solution of each variable can be obtained by the ADMM algorithm, so that the algorithm can be effectively implemented, thus maintaining the accuracy. Furthermore, the experimental results show that the relative error (RE) of these two algorithms always decreases, hence ensuring that they are indeed convergent.

%we propose two LRTC models, i.e., TLNM and TLNMTV, and design two effective alternating direction multiplier method (ADMM) algorithms \cite{683201156}, \cite{6122011620} to solve them respectively. The closed-form solution of each variable is obtained by the ADMM algorithm, leading to the algorithm can be implemented effectively and the accuracy is therefore guaranteed. Experiments results illustrate that the proposed algorithm can guarantee the decrease of relative error (RE), thus ensuring the convergence of the algorithm.

Thirdly, extensive real-data experiments show that the TLNM method obviously outperforms the classical TNN method, which indicates that our proposed integrating QR decomposition and SNN method outperforms the TNN-based TLNM-TQR \cite{2021108240} method. More importantly, the TLMNTV method further enhances the local sparsity, with higher quantitative numerical result and better visual restoration effect than other comparison methods in this paper. It is also particularly worth mentioning that when the sampling rate of hyperspectral images is 2.5\%, our methods can perform image restoration clearly and very well and are shown to be the best.

%a large number of real data experiments show that the TLNM method is obviously superior to the classical TNN method, which indicates that the QR decomposition combined with SNN method proposed by us is superior to TLNM-TQR \cite{2021108240} method based on the TNN method. The TLNMTV method, which further enhances local sparsity, has higher quantitative numerical effect and better visual recovery effect compared with the comparison method in this paper. In addition, it is worth mentioning that when the sampling rate of hyperspectral image is 2.5\%, the image recovery effect of our method is still the best.

%The summary of this article is as follows: in Section II, some preliminary knowledge of tensors are given. The main results, including the proposed model and algorithm, are shown in Section III. The results of extensive experiments and discussion are presented in Section IV. Conclusions are drawn in section V.

\section{PRELIMINARIES}
We are now in a position to give the basic notation of the tensors involved and lay out the basic definitions and theorems for constructing the two proposed new methods.
\subsection{Tensor Notations and their Definitions}
Generally, the uppercase and lowercase of any letter represent respectively vectors and matrices, such as $x$ represents a vector and $X$ represents a matrix. We use a calligraphic upper case letter $\mathcal{X}\in \mathbb{R}^{\mathit{I}_{1}\times\mathit{I}_{2}\times\cdots\times\mathit{I}_{N}}$ to represent an $N$th-order tensor, and $\mathit{x}_{i_{1},i_{2},\cdots,i_{N}}$ to represent its $(i_{1},i_{2},\cdots,i_{N})$-th element. Unless otherwise stated, we still use $\|\mathcal{X}\|_{F}=(\sum_{i_{1},i_{2},\cdots,i_{N}}\mathit{x}_{i_{1},i_{2},\cdots,i_{N}}^{2})^{1/2}$ as the definition of the Frobenius norm and define the inner product of two $N$th-order tensors $\mathcal{Y}$ and $\mathcal{Z}$ as $\langle\mathcal{Y},\mathcal{Z}\rangle=\sum_{i_{1},i_{2},\cdots,i_{N}}\mathit{y}_{i_{1},i_{2},\cdots,i_{N}}\mathit{z}_{i_{1},i_{2},\cdots,i_{N}}$, where $\mathit{y}_{i_{1},i_{2},\cdots,i_{N}}$ and $\mathit{z}_{i_{1},i_{2},\cdots,i_{N}}$ stand for the $(i_{1},i_{2},\cdots,i_{N})$-th element of $\mathcal{Y}$ and $\mathcal{Z}$, respectively.
\begin{definition}[\cite{12345152009}]
	The mode-$n$ unfolding of a tensor $\mathcal{X}\in \mathbb{R}^{\mathit{I}_{1}\times\mathit{I}_{2}\times\cdots\times\mathit{I}_{N}}$ is denoted as a matrix $\mathbf{X}_{(n)}\in\mathbb{R}^{\mathit{I}_{n}\times\mathit{I}_{1}\cdots\mathit{I}_{n-1}\mathit{I}_{n+1}\cdots\mathit{I}_{N}} $. Tensor element $(i_{1}, i_{2},...,i_{N} )$ maps to matrix element $(i_{n}, j)$, where
	\begin{eqnarray}
	j=1+\sum_{k=1,k\neq n}^{N}(i_{k}-1)\mathit{J}_{k}\quad with\quad \mathit{J}_{k}=\prod_{m=1,m\neq n}^{k-1}\mathit{I}_{m}.
	\end{eqnarray}
	
The mode-$n$ unfolding operator and its corresponding inverse operator are abbreviated as $unfold_{n}$ and $fold_{n}$, and satisfy transformation relation $\mathcal{X}=fold_{n}(\mathbf{X}_{(n)})=fold_{n}(unfold_{n}(\mathcal{X}))$.
\end{definition}
\begin{definition}[\cite{12345152009}]
The mode-$n$ product operation between tensor $\mathcal{X}\in \mathbb{R}^{\mathit{I}_{1}\times\mathit{I}_{2}\times\cdots\times\mathit{I}_{N}}$ and matrix $U\in\mathbb{R}^{\mathit{J}_{n}\times \mathit{I}_{n}}$ is expressed as $\mathcal{Y}=\mathcal{X}\times_{n}U$, where $\mathcal{Y}\in\mathbb{R}^{\mathit{I}_{1}\times\mathit{I}_{2}\times\cdots\mathit{I}_{n-1}\mathit{J}_{n}\mathit{I}_{n+1}\cdots\mathit{I}_{N}}$. Elementwisely, we have
	\begin{eqnarray}
	\mathcal{Y}=\mathcal{X}\times_{n}U\quad\Leftrightarrow\quad \mathbf{Y}_{(n)}=U\cdotp unfold_{n}(\mathbf{X}_{(n)}).
	\end{eqnarray}
\end{definition}
\begin{definition}[\cite{1028361472}]
	The $L_{2,1}$ norm of a matirx $M\in\mathbb{R}^{\mathit{I}\times\mathit{J}}$ can be defined as
	\begin{eqnarray}
	\|M\|_{2,1}=\sum_{j=1}^{\mathit{J}}\sqrt{\sum_{i=1}^{\mathit{I}}m_{ij}^{2}};
	\end{eqnarray}
	where $m_{i,j}$ represents the element at row i-th and column j-th of the matrix $M$.
	
\end{definition}

\begin{algorithm}[t]
	\caption{CSVD-QR\cite{8421060}}
	\hspace*{0.02in} {\bf Input:} $X$, a real matrix $C\in\mathbb{R}^{m\times n}$; $It_{max}>0$\\
	\hspace*{0.02in} {\bf Initialization:} %???????
	$r>0$, $q>0$, $k=1$; $\varepsilon$ is a positive tolerance, $L_{0}=eye(m,r)$, $D_{0}=eye(r,r)$, $R_{0}=eye(r,n)$.
	\begin{algorithmic}
		%\State some description % \State ??????		
		\While{$\|L_{k}D_{k}R_{k}-X\|\geq\varepsilon$ or $k<It_{max}$} % While??, ???EndWhile??
		%\For{k=1,2,$\cdots$,K} % For ??, ???EndFor??
		\State $[Q,T]=qr(XR_{k}^{T})$;
		\State $L_{k+1}=Q(:,1:r)$.
		\State $[Q,T]=qr(X^{T}L_{j+1})$;
		\State $R_{k+1}=Q(:,1:r)^{T}$.	
		\State $D_{k+1}=T(1:r,1:r)^{T}$.
		\State $k=k+1$.
		\EndWhile
		\State \Return $L=L_{k}$, $D=D_{k}$, $R=R_{k}$
	\end{algorithmic}
	\hspace*{0.02in} {\bf Output:} %???????
	$L$, $D$, $R$ ($X=LDR$)\label{CSVDQR1}
\end{algorithm}

The $L_{2,1}$ norm has been successfully applied to low-rank representation \cite{1028361472} to optimize the noise data matrix $K\in\mathbb{R}^{m\times n}$, which can be updated by solving the minimization problem as follows:
\begin{eqnarray}
\min_{K} \tau\|K\|_{2,1}+\frac{1}{2}\|K-C\|_{F}^{2} \label{l21problem},
\end{eqnarray}
where $C\in\mathbb{R}^{m\times n}$ is a pre-determined real matrix and $\tau>0$.

Here we will describe an approximate SVD (CSVD-QR) method based on QR decomposition mainly used in this paper, see \cite{8421060}, and the corresponding CSVD-QR algorithm can be specifically shown in Algorithm 1 below.
\begin{theorem}[$L_{2,1}$ norm minimization solver (LNMS) \cite{8421060}]
	The optimal $K(;,j)$, i.e., the $j$th column of $K$, of the problem in (\ref{l21problem}) follows
	\begin{eqnarray}
	K(;,j)=\dfrac{max\{{\|C(:,j)\|_{2}-\tau,0}\}}{\|C(:,j)\|_{2}}C(:,j), \label{lnms}
	\end{eqnarray}
	where $\|C(:,j)\|_{2}=\sqrt{\sum_{i=1}^{m}C_{ij}^{2}}$.\label{LNMSS}
\end{theorem}
\section{Models and Algorithms}
\subsection{Tensor $L_{2,1}$ Norm Minimization Model With Its Algorithm}
For LRTC problem, to further improve the quality and visual effect of image restoration, we integrate the $L_{2,1}$ norm and CSVD-QR method into the tensor SNN minimization problem instead of the classical tensor TNN method, resulting in a tensor $L_{2,1}$ norm minimization model as below:

\begin{eqnarray}
&&\min_{\mathcal{X}}\sum_{n=1}^{N}\alpha_{n}\|D_{n}\|_{2,1}\nonumber
\\&&s.t.\quad\{\mathcal{M}_{n}=\mathcal{X}\}_{n=1}^{N}, \quad\{\mathbf{M}_{n(n)}=L_{n}D_{n}R_{n}\}_{n=1}^{N},\nonumber
\\&&\qquad\quad\mathcal{X}_{\Omega}=\mathcal{T}_{\Omega},\quad L_{n}^{T}L_{n}=I_{n},\quad R_{n}R_{n}^{T}=I_{n}, \label{tensor2}
\end{eqnarray}
where $\mathcal{M}_{n},n=1,2,\dots,N$ indicate auxiliary parameters, $\sum_{n=1}^{N}\alpha_{n}=1$, $L_{n}$, $D_{n}$, $R_{n}$ are induced by $\mathbf{M}_{n(n)}$ when completing CSVD-QR decomposition, $\Omega$ represents the exact locations of known observations. It is worth mentioning here that we still adopt the $L_{2,1}$ norm accommodating with CSVD-QR method and the tensor SNN minimization method to replace the nuclear norm in the objective function, which exactly shows the novelty and superiority of proposed method compared with the existing state-of-the-art methods and ensures the fairness in performing trial comparisons between two recently published methods \cite{8421060} and \cite{2021108240}.
\subsection{Optimization TLNM Algorithm Based on ADMM}
With the augmented Lagrangian formula, the above model (\ref{tensor2}) can be transformed into the following optimization form:
\begin{eqnarray}
&&Lag(\mathcal{X},\{\mathcal{M}_{n},\mathcal{Q}_{n},L_{n},R_{n},D_{n},\Phi_{n}\}_{n=1}^{N})\nonumber
\\&&=\sum_{n=1}^{N}\alpha_{n}\|D_{n}\|_{2,1}+\frac{\mu}{2}\|\mathbf{M}_{n(n)}-L_{n}D_{n}R_{n}+\frac{\Phi_{n}}{\mu}\|_{F}^{2}\nonumber
\\&&+\frac{\mu}{2}\|\mathcal{M}_{n}-\mathcal{X}+\frac{\mathcal{Q}_{n}}{\mu}\|_{F}^{2}, \label{TLNM111}
\end{eqnarray}
where matrices $\mathcal{Q}_{n}$ and $\Phi_{n}$ denote the Lagrange multipliers and $\mu>0$. Now, we start to optimize the problem (\ref{TLNM111}) under the ADMM framework. Here we use these matrix symbols $L_{n}^{+}, R_{n}^{+}, D_{n}^{+}, \mathcal{M}_{n}^{+}, \mathcal{X}^{+}$ to represent the final results of iterative updates in ADMM respectively.
\subsubsection{Optimize \{$L_{1},...,L_{N}$\}}
Let the initial values of all parameter matrices $L_{n},D_{n},R_{n}$ in the proposed algorithm be $L_{n}^{0}=eye(\mathit{I}_{n},r_{n})$, $D_{n}^{0}=eye(r_{n},r_{n})$, $R_{n}^{0}=eye(r,t_{n})$ and $t_{n}=\mathit{I}_{1}\cdots\mathit{I}_{n-1}\mathit{I}_{n+1}\cdots\mathit{I}_{N}$, respectively, where $eye$ is the built-in command for the identity matrix in MATLAB. Keep other variables unchanged in (\ref{TLNM111}), the optimization objective function with respect to $L_{n}$ becomes
\begin{eqnarray}
\min_{L_{n},R_{n}} \|\mathbf{M}_{n(n)}-L_{n}D_{n}R_{n}+\frac{\Phi_{n}}{\mu}\|_{F}^{2}.\label{CSVD111}
\end{eqnarray}
Here, we regard the solution of problem (\ref{CSVD111}) as the iterative solution process of the CSVD-QR method based on matrix $\mathbf{M}_{n(n)}+\frac{\Phi_{n}}{\mu} $, which can be accurately expressed as follows:
\begin{eqnarray}
\mathbf{M}_{n(n)}+\frac{\Phi_{n}}{\mu}=L_{n}D_{Tn}R_{n},\label{LR}
\end{eqnarray}
where $D_{Tn}\in\mathbb{R}^{r_{n}\times r_{n}}$. If we let $\mathbf{M}_{n(n)}^{k}+\frac{\Phi_{n}^{k}}{\mu}=G_{n}$, then $L_{n}^{k+1}$ can be obtained recursively as follows:
\begin{eqnarray}
[Q,T]=qr(G_{n}R_{n}^{T})\label{UL1}\\
L_{n}^{+}=Q(q_{1},...,q_{r_{n}})\label{UL2}
\end{eqnarray}
where matrices $Q\in\mathbb{R}^{\mathit{I}_{n}\times \mathit{I}_{n}}$ and $T\in\mathbb{R}^{\mathit{I}_{n}\times r_{n}}$ are treated as intermediate variables.
\subsubsection{Optimize \{$R_{1},...,R_{N}$\}}
In a similar manner, $R_{n}^{k+1}$ can be derived iteratively as follows:
\begin{eqnarray}
[Q,T]=qr(G_{n}^{T}L_{n})\label{UR1}\\
R_{n}^{+}=Q(q_{1},...,q_{r_{n}})^{T}\label{UR2}
\end{eqnarray}
where $Q\in\mathbb{R}^{t_{n}\times t_{n}}$ and $T\in\mathbb{R}^{t_{n}\times r_{n}}$ are considered as intermediate variables.
\subsubsection{Optimize \{$D_{1},...,D_{N}$\}}
With $\mathbf{M}_{n(n)}$, $\Phi_{n}$, $L_{n}$, and $R_{n}$ held fixed, the variable $D_{n}$ can be optimally solved in the following way:
\begin{eqnarray}
D_{n}=\arg\min_{D_{n}}\frac{\alpha_{n}}{\mu}\|D_{n}\|_{2,1}+\frac{1}{2}\|D_{n}-L_{n}^{T}(\mathbf{M}_{n(n)}+\frac{\Phi_{n}}{\mu})R_{n}^{T}\|_{F}^{2}.\label{D1}
\end{eqnarray}
Combine (\ref{LR}) with (\ref{D1}), we get:
\begin{eqnarray}
D_{Rn}=L_{n}^{T}(\mathbf{M}_{n(n)}+\frac{\Phi_{n}}{\mu})R_{n}^{T}
\end{eqnarray}
As a result, (\ref{D1}) can be equivalently expressed as:
\begin{eqnarray}
D_{n}=\arg\min_{D_{n}}\frac{\alpha_{n}}{\mu}\|D_{n}\|_{2,1} +\frac{1}{2}\|D_{n}-D_{Rn}\|_{F}^{2}.\label{D2}
\end{eqnarray}
With the result of Theorem \ref{LNMSS}, we can approximate and update $D_{n}$ through the following formula:
\begin{eqnarray}
D_{n}^{+}=D_{Rn}G_{n}\label{UD1},
\end{eqnarray}
where 
\begin{eqnarray}
G_{n}=diag(g_{n1},...,g_{nr_{n}}),
\end{eqnarray}
and the jth entry $g_{nj}$ can be calculated by using the following formula:
\begin{eqnarray}
g_{nj}=\dfrac{max\{\|D_{Rn}(:,j)\|_{F}-\frac{\alpha_{n}}{\mu},0\}}{\|D_{Rn}(:,j)\|_{F}}.
\end{eqnarray}
\subsubsection{Optimize \{$\mathcal{M}_{1},...,\mathcal{M}_{N}$\}}
Under the assumption that other variables remain unchanged, the optimization objective function with respect to $\mathcal{M}_{n}$ can degenerate into
\begin{eqnarray}
&&\mathcal{M}_{n}=\min_{\mathcal{M}_{n}}\|\mathbf{M}_{n(n)}-L_{n}D_{n}R_{n}+\frac{\Phi_{n}}{\mu}\|_{F}^{2}\nonumber
\\&&+\|\mathcal{M}_{n}-\mathcal{X}+\frac{\mathcal{Q}_{n}}{\mu}\|_{F}^{2}.\label{tensorM1}
\end{eqnarray}
Using the definitions and properties of tensor (matrix) Frobenius norm, $\mathcal{M}_{n}$ of (\ref{tensorM1}) can be equivalently reformulated as
\begin{eqnarray}
&&\mathbf{M}_{n(n)}=\min_{\mathcal{M}_{n}}\|\mathbf{M}_{n(n)}-L_{n}D_{n}R_{n}+\frac{\Phi_{n}}{\mu}\|_{F}^{2}\nonumber
\\&&+\|\mathbf{M}_{n(n)}-\mathbf{X}_{(n)}+\frac{\mathbf{Q}_{n(n)}}{\mu}\|_{F}^{2}.\label{matrixM}
\end{eqnarray}
Accordingly, $\mathcal{M}_{n}^{+}$ is deduced as follows:
\begin{eqnarray}
\mathcal{M}_{n}^{+}=fold_{n}(\frac{1}{2}(\mathbf{X}_{(n)}-\frac{\mathbf{Q}_{n(n)}}{\mu}+L_{n}D_{n}R_{n}-\frac{\Phi_{n}}{\mu})). \label{UM1}
\end{eqnarray}

\begin{algorithm}[t]
	\caption{TLNM} %?????
	\hspace*{0.02in} {\bf Input:} %?????,  \hspace*{0.02in}??????, ???? \\ ????
	any incomplete tensor $\mathcal{T}$, the index set representing the position of known tensor elements $\Omega$, convergence criteria $\epsilon$, prescribed iteration number $K$. \\
	\hspace*{0.02in} {\bf Initialization:} %???????
	$\mathcal{X}^{0}=\mathcal{T}_{\Omega}$, $\{\mathcal{M}_{n}^{0}=\mathcal{X}^{0}\}_{n=1}^{N}$, $\mu^{0}>0$, $\rho>1$, $\{\alpha_{n},r_{n},L_{n}^{0},D_{n}^{0},R_{n}^{0}\}_{n=1}^{N}$.
	\begin{algorithmic}[1]
		%\State some description % \State ??????		
		\While{not converged and $k<K$} % While??, ???EndWhile??
		%\For{k=1,2,$\cdots$,K} % For ??, ???EndFor??
		\For{n=1:N}
		\State Compute $L_{n}^{k}$ via (\ref{UL1}-\ref{UL2});		
		\State Compute $R_{n}^{k}$ via (\ref{UR1}-\ref{UR2});		
		\State Compute $D_{n}^{k}$ via (\ref{UD1});
		\State Compute $\mathcal{M}_{n}^{k}$ via (\ref{UM1});
		%		\State Updating the multipliers $\Phi_{n}^{k}$ via (\ref{Uphi1});
		\EndFor
		\State Compute $\mathcal{X}^{k}$ via (\ref{UX1});
		\State Compute the multipliers $\mathcal{Q}_{n}^{k}$ and $\Phi_{n}^{k}$ via (\ref{UP1});
		\State $\mu^{k}=\rho\mu^{k-1}$, $k=k+1$;
		\State Check whether the convergence condition $\|\mathcal{X}^{k+1}-\mathcal{X}^{k}\|_{\infty}\leq\epsilon$ is satisfied
		%	\EndFor
		\EndWhile
		\State \Return $\mathcal{X}^{k+1}$
	\end{algorithmic}
	\hspace*{0.02in} {\bf Output:} %???????
	Completed tensor $\mathcal{X}=\mathcal{X}^{k+1}$\label{TLNMA}
\end{algorithm}
\subsubsection{Optimize $\mathcal{X}$}
The optimal solution for the variable $\mathcal{X}$ can be obtained by optimizing the minimization function of the following problem:
\begin{eqnarray}
\mathcal{X}^{+}=\min_{\mathcal{X}}\sum_{n=1}^{N}\frac{\mu}{2}\|\mathcal{M}_{n}-\mathcal{X}+\frac{\mathcal{Q}_{n}}{\mu}\|_{F}^{2}\nonumber\quad s.t.\quad \mathcal{X}_{\Omega}=\mathcal{T}_{\Omega}.
\end{eqnarray}
its specific explicit closed-form solution can be written as
\begin{eqnarray}
\left\{\begin{array}{l}
\mathcal{X}^{+}_{\Omega}=\mathcal{T}_{\Omega},\\\mathcal{X}^{+}_{\Omega^{\perp}}=\frac{1}{N\mu}(\sum_{n=1}^{N}\mu\mathcal{M}_{n}+\mathcal{Q}_{n})_{\Omega^{\perp}},
\end{array}
\right.\label{UX1}
\end{eqnarray}
where the vertical symbol $^{\perp}$ in the upper right corner of $\Omega$ indicates the complement set $\Omega^{\perp}$ of $\Omega$.
\subsubsection{ Update the multipliers $\mathcal{Q}_{n}$, $\Phi_{n}$}
\begin{eqnarray}
\left\{\begin{array}{l}
\mathcal{Q}_{n}^{+}=\mathcal{Q}_{n}+\mu(\mathcal{M}_{n}-\mathcal{X}),
\\\Phi_{n}^{+}=\Phi_{n}+\mu(\mathbf{M}_{n(n)}-L_{n}D_{n}R_{n}),
\\\mu^{+}=\rho\mu,
\end{array}
\right.\label{UP1}
\end{eqnarray}
where $\rho\geq1$. 

So far, all the above variables have been updated and solved, and the pseudocode corresponding to all algorithms for TLNM has been reported in detail in Algorithm \ref{TLNMA}.

\subsection{Tensor $L_{2,1}$-Norm Minimization with Total variation (TLNMTV)}
Total variation method is inherently remarkable in the protection of image details so as to improve the utilization rate of the local prior information of the image. In this regard, we investigate and propose the tensor $L_ {2, 1} $ norm minimization with TV (TLNMTV) method to obtain the following model by adding the regularization term of the form in \cite{10.5555/3298483.3298558}.
\begin{eqnarray}
\min_{\mathcal{X}}  \sum_{n=1}^{N}\alpha_{n}\|\mathbf{X}_{(n)}\|_{\ast} +\lambda\sum_{n=1}^{N}\beta_{n}\rvert C_{n}\mathbf{X}_{(n)}\rvert \quad s.t.\quad\mathcal{X}_{\Omega}=\mathcal{T}_{\Omega}.
\end{eqnarray}
Here, the latter term is the TV regularization term, where the dimension of $C_{n}$ is $(\mathit{I}_{n}-1,\mathit{I}_{n})$ with $(C_{n})_{i,i}=1$ and $(C_{n})_{i,i+1}=-1$. The operator $\rvert\cdot\rvert$ of any matrix $X$ is defined as $\rvert X\rvert=\sum_{i=1}\sum_{j=1}\rvert X_{i,j}\rvert$. The parameter $\lambda$ is represented as a tunable hyperparameter. These parameters $\beta_{1},\cdots,\beta_{N}$ take either the value 0 or 1, which means whether the smoothing and piecewise prior information on the $n$-th mode of the completed tensor is known or not. However, as we know, different dimensions of different high-dimensional data represent different data information, and some information is not smooth. For this, it is not necessary to compute smoothing and piecewise priors in this dimension. It should be pointed out that, quite different from the traditional regular meaning, the TV term used here is only able to determine whether smoothing and piecewise priors need to be computed for different modes. Therefore, it is more reasonable to use the TV regularization term in combination with the actual situation. In particular, when the TV regularization term is removed, the problem degenerates itself into a low-rank tensor problem (\ref{tensor1}). 

In this connection, by introducing auxiliary matrix and tensor variable sets $\{Q_{n}\}_{n=1}^{N}$,  $\{A_{n}\}_{n=1}^{N}$ $\{\mathcal{Z}_{n}\}_{n=1}^{N}$ for these variables $\mathbf{X}_{(n)}$, $C_{n}A_{n}$ and $\mathcal{X}$, respectively, i.e., $A_{n}=\mathbf{X}_{(n)}$, $Q_{n}=C_{n}A_{n}$ and $\mathcal{Z}_{n}=\mathcal{X}$, the proposed tensor $L_{2,1}$-norm minimization with TV model (TLNMTV) can be induced as follows:
\begin{eqnarray}
&&\min_{\mathcal{X}}\sum_{n=1}^{N}\alpha_{n}\|D_{n}\|_{2,1}+\lambda\sum_{n=1}^{N}\beta_{n}\rvert Q_{n}\rvert\nonumber
\\&&\quad s.t.\quad\mathcal{X}_{\Omega}=\mathcal{T}_{\Omega},\{Q_{n}=C_{n}A_{n},\mathcal{Z}_{n}=\mathcal{X},A_{n}=\mathbf{X}_{(n)}\}_{n=1}^{N},
\nonumber
\\&&\qquad\{\mathbf{Z}_{n(n)}=L_{n}D_{n}R_{n}\}_{n=1}^{N}, L_{n}^{T}L_{n}=I_{n}, R_{n}R_{n}^{T}=I_{n}. \label{TLNMTV1}
\end{eqnarray}
Due to the fact that this model (\ref{TLNMTV1}) is convex we can solve it using the ADMM framework.

\subsection{Optimization TLNMTV Algorithm Based on ADMM}
With the augmented Lagrange thinking, the optimization problem (\ref{TLNMTV1}) can be transformed into:
\begin{eqnarray}
&&Lag(\mathcal{X},\{\mathcal{Z}_{n},\mathcal{G}_{n},Q_{n},\Lambda_{n},A_{n},\Gamma_{n},L_{n},R_{n},D_{n},\Phi_{n}\}_{n=1}^{N})\nonumber
\\&&=\sum_{n=1}^{N}(\alpha_{n}\|D_{n}\|_{2,1}+\frac{\mu}{2}\|\mathcal{Z}_{n}-\mathcal{X}+\frac{\mathcal{G}_{n}}{\mu}\|_{F}^{2})\nonumber
\\&&+\sum_{n=1}^{N}\beta_{n}(\lambda\rvert Q_{n}\rvert+\frac{\mu}{2}\|Q_{n}-C_{n}A_{n}+\frac{\Lambda_{n}}{\mu}\|_{F}^{2})\nonumber
\\&&+\sum_{n=1}^{N}\beta_{n}(\frac{\mu}{2}\|A_{n}-\mathbf{X}_{(n)}+\frac{\Gamma_{n}}{\mu}\|_{F}^{2})\nonumber
\\&&+\sum_{n=1}^{N}\frac{\mu}{2}\|\mathbf{Z}_{n(n)}-L_{n}D_{n}R_{n}+\frac{\Phi_{n}}{\mu}\|_{F}^{2},
\end{eqnarray}
where matrices $\{\Lambda_{n}\}_{n=1}^{N}$, $\{\Gamma_{n}\}_{n=1}^{N}$, $\{\Phi_{n}\}_{n=1}^{N}$ and tensors $\{\mathcal{G}_{n}\}_{n=1}^{N}$ represent, respectively, Lagrange multipliers, $\mu>0$. Next, we explicitly give their respective formulas for these optimization variables $\{Q_{n}\}_{n=1}^{N}$, $\{L_{n}\}_{n=1}^{N}$, $\{R_{n}\}_{n=1}^{N}$, $\{D_{n}\}_{n=1}^{N}$, $\{\mathcal{Z}_{n}\}_{n=1}^{N}$, $\{A_{n}\}_{n=1}^{N}$ and $\mathcal{X}$.
\subsubsection{Optimize \{$Q_{1},\cdots,Q_{N}$\}}
Under the premise of ensuring that other variables do not change, we can transform the optimization problem of \{$Q_{1},\cdots,Q_{N}$\} into a form of:
\begin{eqnarray}
\min_{\{Q_{1},\cdots,Q_{N}\}}\sum_{n=1}^{N}\beta_{n}(\lambda\rvert Q_{n}\rvert+\frac{\mu}{2}\|Q_{n}-C_{n}A_{n}+\frac{\Lambda_{n}^{k}}{\mu}\|_{F}^{2}).
\end{eqnarray}
Based on the fact that the matrix groups \{$Q_{1},\cdots,Q_{N}$\} in the optimization problem are independent of each other, the solution formula for each of them can be easily obtained:
\begin{eqnarray}
Q_{n}=\beta_{n}\cdot shrinkage_{\frac{\lambda}{\mu}}(C_{n}A_{n}-\frac{1}{\mu}\Lambda_{n}),
\label{UQ1}\end{eqnarray}
where the operator symbol $shrinkage_{\alpha}(\cdot)$ represents the elementwise shrinkage thresholding operator of a matrix, i.e.,
\begin{eqnarray}
shrinkage_{\alpha}(X)_{i,j}=(X)_{i,j}-\min(\alpha,\rvert (X)_{i,j}\rvert)\cdot\frac{(X)_{i,j}}{\rvert (X)_{i,j}\rvert}\nonumber,
\end{eqnarray}
and $\frac{(X)_{i,j}}{\rvert (X)_{i,j}\rvert}$ is equal to zero when $(X)_{i,j}=0$.
\subsubsection{Optimize \{$L_{1},...,L_{N}$\},\{$R_{1},...,R_{N}$\},\{$D_{1},...,D_{N}$\}}
Since the update of these variables $L_{1},...,L_{N}$, $R_{1},...,R_{N}$ and $D_{1},...,D_{N}$ is consistent with the TLNM algorithm, it is omitted here and will not be repeated.
\subsubsection{Optimize \{$\mathcal{Z}_{1},...,\mathcal{Z}_{N}$\}}
By keeping the rest of the variables unchanged, the optimization problem with respect to $\mathcal{Z}_{n}$ can be precisely described as below
\begin{eqnarray}
\mathcal{Z}_{n}=\min_{\mathcal{Z}_{n}}\|\mathbf{Z}_{n(n)}-L_{n}D_{n}R_{n}+\frac{\Phi_{n}}{\mu}\|_{F}^{2}+\|\mathcal{Z}_{n}-\mathcal{X}+\frac{\mathcal{Q}_{n}}{\mu}\|_{F}^{2}.\label{tensorZ1}
\end{eqnarray}
With the help of the definition of tensor Frobenius norm and matrix Frobenius norm, the above problem (\ref{tensorZ1}) can be equivalently transformed into
\begin{eqnarray}
\mathbf{Z}_{n(n)}=\min_{\mathcal{Z}_{n}}\|2\mathbf{Z}_{n(n)}-L_{n}D_{n}R_{n}+\frac{\Phi_{n}}{\mu}-\mathbf{X}_{(n)}+\frac{\mathbf{Q}_{n(n)}}{\mu}\|_{F}^{2}.\label{matrixZ}
\end{eqnarray}
Therefore $\mathcal{Z}_{n}^{+}$ is obtained as follows:
\begin{eqnarray}
\mathcal{Z}_{n}^{+}=fold_{n}(\frac{1}{2\mu}(\mu\mathbf{X}_{(n)}-\mathbf{Q}_{n(n)}+\mu L_{n}D_{n}R_{n}-\Phi_{n})).\label{UZ1}
\end{eqnarray}

\subsubsection{Optimize \{$A_{1},\cdots,A_{N}$\}}

Make sure that the rest of the above variables do not change, we can reformulate the above optimization problem of \{$A_{1},\cdots,A_{N}$\} into the following form:
\begin{eqnarray}
\min\sum_{n=1}^{N}(\|Q_{n}-C_{n}A_{n}+\frac{\Lambda_{n}}{\mu}\|_{F}^{2}+\|A_{n}-\mathbf{X}_{(n)}+\frac{\Gamma_{n}}{\mu}\|_{F}^{2})\label{A11}.
\end{eqnarray}
This type of problem (\ref{A11}) can be solved by dividing it into $N$ sub-problems, each of which can be expressed in this form:
\begin{eqnarray}
A_{n}=\arg\min_{A_{n}}\|Q_{n}-C_{n}A_{n}+\frac{\Lambda_{n}}{\mu}\|_{F}^{2}+\|A_{n}-\mathbf{X}_{(n)}+\frac{\Gamma_{n}}{\mu}\|_{F}^{2}.
\end{eqnarray}
In this regard, by solving the above minimization problem, the following update formula is derived:
\begin{eqnarray}
A_{n}^{+}=(\mu C_{n}^{T}C_{n}+\mu I)^{-1}(C_{n}^{T}\Lambda_{n}+\mu C_{n}^{T}Q_{n}+\mu\mathbf{X}_{(n)}-\Gamma_{n}),
\label{UA1}\end{eqnarray}
where $I$ is the identity matrix of the appropriate order.
\subsubsection{Optimize $\mathcal{X}$}
We obtain the variable $\mathcal{X}$ by solving the following constraint minimization function:
\begin{eqnarray}
&&\mathcal{X}=\min_{\mathcal{X}}\sum_{n=1}^{N}(\|\mathcal{Z}_{n}-\mathcal{X}+\frac{\mathcal{G}_{n}}{\mu}\|_{F}^{2}+\beta_{n}\|A_{n}-\mathbf{X}_{(n)}+\frac{\Gamma_{n}}{\mu}\|_{F}^{2})\nonumber \\&&s.t.\quad \mathcal{X}_{\Omega}=\mathcal{T}_{\Omega},
\end{eqnarray}
whose closed-form solution can be expressed as
\begin{eqnarray}
\left\{\begin{array}{l}
\mathcal{X}^{+}_{\Omega}=\mathcal{T}_{\Omega},\\\mathcal{X}^{+}_{\Omega^{\perp}}=\frac{[\sum_{n=1}^{N}\mu\mathcal{Z}_{n}+\mathcal{G}_{n}+fold_{n}(\mu A_{n}+\Gamma_{n})]_{\Omega^{\perp}}}{(N+\sum_{n=1}^{N}\beta_{n})\mu}.
\end{array}
\right.\label{UX2}
\end{eqnarray}

\begin{algorithm}[t]
	\caption{TLNMTV} %?????
	\hspace*{0.02in} {\bf Input:} %?????,  \hspace*{0.02in}??????, ???? \\ ????
	an arbitrary incomplete tensor $\mathcal{T}$, the set of indices of the known elements $\Omega$, convergence stop criteria $\epsilon$, maximum number of iteration steps $K$. \\
	\hspace*{0.02in} {\bf Initialization:} %???????
	$\mathcal{X}^{0}=\mathcal{T}_{\Omega}$, $\{\mathcal{Z}_{n}^{0}=\mathcal{X}^{0}\}_{n=1}^{N}$, $\mu^{0}>0$, $\rho>1$, $\{\alpha_{n},r_{n},L_{n}^{0},D_{n}^{0},R_{n}^{0}\}_{n=1}^{N}$.
	\begin{algorithmic}[1]
		%\State some description % \State ??????		
		\While{not converged and $k<K$} % While??, ???EndWhile??
		%\For{k=1,2,$\cdots$,K} % For ??, ???EndFor??
		\For{n=1:$N$}
		\State Calculate $Q_{n}^{k}$ via (\ref{UQ1});
		\State Calculate $L_{n}^{k}$ via (\ref{UL1}-\ref{UL2});		
		\State Calculate $R_{n}^{k}$ via (\ref{UR1}-\ref{UR2});		
		\State Calculate $D_{n}^{k}$ via (\ref{UD1});
		\State Calculate $\mathcal{Z}_{n}^{k}$ via (\ref{UZ1});
		%		\State Updating the multipliers $\phi_{n}^{k}$ via (\ref{Uphi2});
		\State Calculate $A_{n}^{k}$ via (\ref{UA1});
		\EndFor
		\State Calculate $\mathcal{X}^{k}$ via (\ref{UX2});
		\State Calculate the multipliers $\mathcal{G}_{n}$, $\Lambda_{n}$, $\Phi_{n}$, and $\Gamma_{n}$ via (\ref{UP2});
		\State  $\mu^{k+1}=\rho\mu^{k}$, $k=k+1$;
		\State Check whether the convergence condition $\|\mathcal{X}^{k+1}-\mathcal{X}^{k}\|_{\infty}\leq\epsilon$ is satisfied
		%	\EndFor
		\EndWhile
		\State \Return $\mathcal{X}^{k+1}$
	\end{algorithmic}
	\hspace*{0.02in} {\bf Output:} %???????
	Completed tensor $\mathcal{X}=\mathcal{X}^{k+1}$\label{TLNMTVA}
\end{algorithm}
\subsubsection{Update the multipliers $\mathcal{G}_{n}$, $\Lambda_{n}$, $\Phi_{n}$ and $\Gamma_{n}$}
\begin{eqnarray}
\left\{\begin{array}{l}
\mathcal{G}_{n}^{+}=\mathcal{G}_{n}+\mu(\mathcal{Z}_{n}-\mathcal{X}),
\\\Lambda_{n}^{+}=\Lambda_{n}+\mu(Q_{n}-C_{n}A_{n}),
\\\Gamma_{n}^{+}=\Gamma_{n}+\mu(A_{n}-\mathbf{X}_{(n)}),
\\\Phi_{n}^{+}=\Phi_{n}+\mu(\mathbf{Z}_{n(n)}-L_{n}D_{n}R_{n}),
\end{array}
\right.\label{UP2}
\end{eqnarray}
and $\mu^{+}=\rho\mu$, where $\rho\geq1$. The pseudo-code corresponding to our proposed algorithm for TLNMTV is displayed in the table, i.e., Algorithm \ref{TLNMTVA}.

\subsection{Computational Complexity Analysis of Corresponding Algorithms}
For any input tensor $\mathcal{T}\in \mathbb{R}^{\mathit{I}_{1}\times\mathit{I}_{2}\times\cdots\times\mathit{I}_{N}}$ of size $N$-order, the entire computational complexity of Algorithm \ref{TLNMA} includes the following five major aspects.

\par Updating the matrix $L_{n}$ needs to give the decomposition of QR [see (\ref{UL1})-(\ref{UL2})], and the complexity of which is $O(\mathit{I}_{n}r_{n}^{2})$, and the total computational complexity of these variable matrices $L_{n}(n=1,2,\dots,N)$ is $O(\sum_{n=1}^{N}\mathit{I}_{n}r_{n}^{2})$.

\par Likewise, the complexity of calculating $R_{n}$ is $O(\sum_{n=1}^{N}r_{n}^{2}\prod_{i=1,i\neq n}^{N}\mathit{I}_{i})$.
\par
When updating matrix $D_{n}$, it is necessary to perform the product operation between two $r_{n}\times r_{n}$ matrices for solving the LNMS problem, thus the computational complexity of all these variable matrices $D_{n}(n=1,2,\dots,N)$ is $O(\sum_{n=1}^{N}r_{n}^{2})$.
\par
Updating each matrix $\mathcal{M}_{n}$ requires the product operation of these three matrices $L_{n}$, $D_{n}$, and $R_{n}$ in advance, for which the computational complexity of all the corresponding variable matrices $\mathcal{M}_{n}(n=1,2,\dots,N)$ is $O(\sum_{n=1}^{N}(\mathit{I}_{n}r_{n}^{2}+r\prod_{i=1}^{N}\mathit{I}_{i}))$.
\par
In addition, the complexity of updating matrix $\mathcal{X}$ is known as $O(\prod_{n}\mathit{I}_{n})$. Thus, by summing the above complexities, all the computational complexity of Algorithm \ref{TLNMA} is $O(K(\sum_{n=1}^{N}(\mathit{I}_{n}r_{n}^{2}+r_{n}^{2}\prod_{i=1,i\neq n}^{N}\mathit{I}_{i}+r_{n}^{2}+\mathit{I}_{n}r_{n}^{2}+r\prod_{i=1}^{N}\mathit{I}_{i})+\prod_{n}\mathit{I}_{n}))$.

Similarly, the entire computational complexity of Algorithm \ref{TLNMTVA} is mainly composed of the following seven parts.
\par
It is easy to know that the computational complexity required to update the matrix $Q_{n}$ is $O(\sum_{n=1}^{N}(\mathit{I}_{n}-1)\prod_{n}\mathit{I}_{n})$.
\par
The computational complexity required to update these matrices $L_{n}$, $R_{n}$, $D_{n}$ is the same as in Algorithm \ref{TLNMA}.
\par
Updating $\mathcal{Z}_{n}$ requires calculating the product of the three matrices $L_{n}D_{n}R_{n}$ in advance, the complexity of these variable matrices $\mathcal{Z}_{n}(n=1,2,\dots,N)$ is $O(\sum_{n=1}^{N}(\mathit{I}_{n}r_{n}^{2}+r\prod_{i=1}^{N}\mathit{I}_{i}))$.
\par
Updating $A_{n}$ requires calculating the inverse of a matrix of order $\mathit{I}_{n}\times \mathit{I}_{n}$ and the product matrix of any matrix, the complexity of which is $O(\mathit{I}_{n}^{3}+\mathit{I}_{n}^{2}\prod_{i=1,i\neq n}^{N}\mathit{I}_{i})$. The complexity of these variable matrices $A_{n}(n=1,2,\dots,N)$ is $O(\sum_{n=1}^{N}\mathit{I}_{n}^{3}+\mathit{I}_{n}^{2}\prod_{i=1,i\neq n}^{N}\mathit{I}_{i})$
\par
The computational complexity required to update a tensor $\mathcal{X}$ is $O(\prod_{n}\mathit{I}_{n})$. Thus, the total computational complexity of our proposed Algorithm \ref{TLNMTVA} is $O(K(\sum_{n=1}^{N}((\mathit{I}_{n}-1)\prod_{n}\mathit{I}_{n}+\mathit{I}_{n}r_{n}^{2}+r_{n}^{2}\prod_{i=1,i\neq n}^{N}\mathit{I}_{i}+r_{n}^{2}+\mathit{I}_{n}r_{n}^{2}+r\prod_{i=1}^{N}\mathit{I}_{i}+\mathit{I}_{n}^{3}+\mathit{I}_{n}^{2}\prod_{i=1,i\neq n}^{N}\mathit{I}_{i})+\prod_{n}\mathit{I}_{n}))$.
\section{Experiments}

In this section, two typical tensor data, i.e., HSI data and MRI data are employed to illustrate the performance of the proposed method. Four quantitative picture quality indices (PQIs) are used to evaluate the quality of recovery, including peak signal\textendash to\textendash noise ratio (PSNR), structural similarity (SSIM) \cite{1284395} feature similarity (FSIM) \cite{5705575}, and erreur relative globale adimensionnelle de synth\`{e}se (ERGAS) \cite{13531574432}. All tests are implemented on the Windows 10 platform and MATLAB (R2019a) with an Intel Core i7-10875H 2.30 GHz and 32 GB of RAM.

We compare our results with nine recently developed state-of-the-art LRTC
methods, including the tensor trace norm-based LRTC (HaLRTC) \cite{6138863}, the t-SVD-based TC method \cite{7782758}, alternating direction multiplier method(ALM)-based matrix completion (MC-ALM) \cite{366666}, the nonconvex tensor rank constraint-based, i.e., the minimax concave plus penalty-based TC (McpTC) and the smoothly clipped absolute deviation penalty-based TC (ScadTC) method \cite{2612015273}, the parallel matrix factorization-based LRTC method (TMac) \cite{6253256}, the LRTC with TV on tensor unfolding (LRTC-TV) \cite{10.5555/3298483.3298558}, Tensor Completion based Framelet Representation of Tensor Nuclear Norm(FTNN)\cite{9115254}, and the joint trace/TV-based TC method (Trace/TV) \cite{20085525245}. All comparison methods use the corresponding optimal parameters for the experiment.

\begin{figure*}[!h]
	\centering
	\subfigure[]{
		\begin{minipage}{0.056\textwidth}
			\includegraphics[width=1\linewidth]{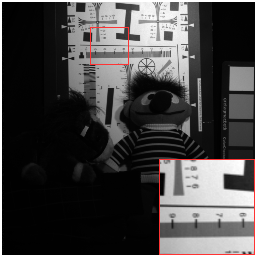}
			\includegraphics[width=1\linewidth]{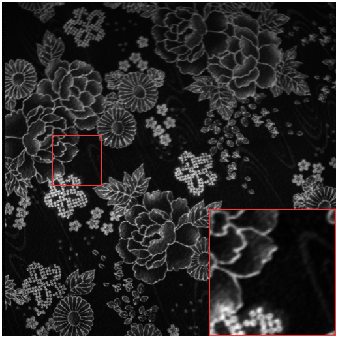}
			\includegraphics[width=1\linewidth]{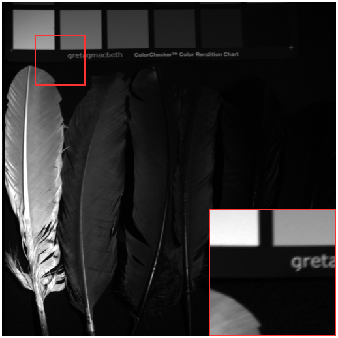}
			\includegraphics[width=1\linewidth]{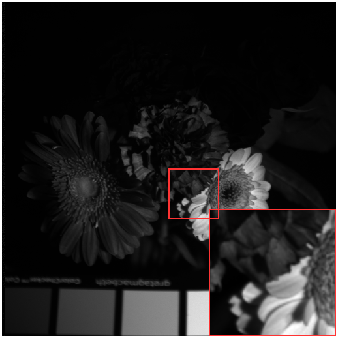}
			\includegraphics[width=1\linewidth]{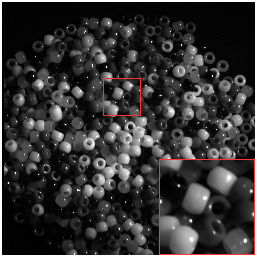}
			\includegraphics[width=1\linewidth]{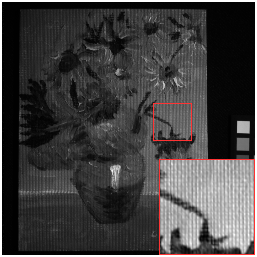}
	\end{minipage}}
	\subfigure[]{
		\begin{minipage}{0.056\textwidth}
			\includegraphics[width=1\linewidth]{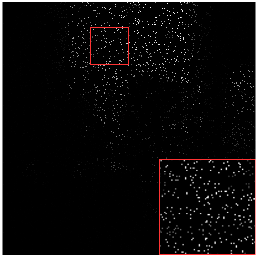}
			\includegraphics[width=1\linewidth]{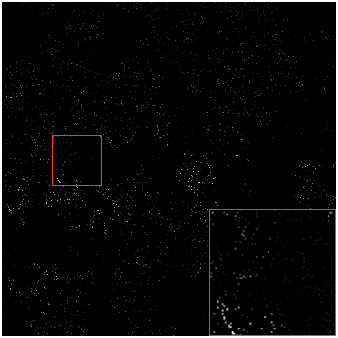}
			\includegraphics[width=1\linewidth]{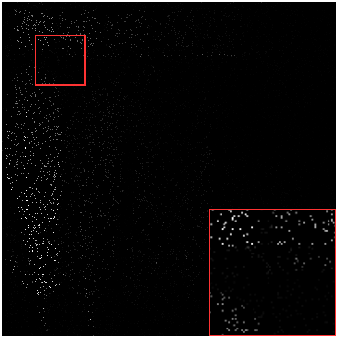}
			\includegraphics[width=1\linewidth]{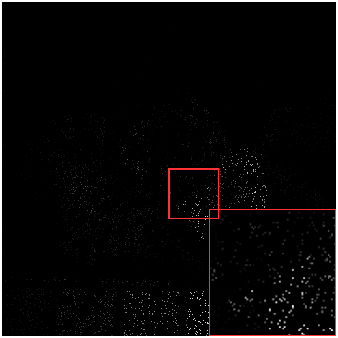}
			\includegraphics[width=1\linewidth]{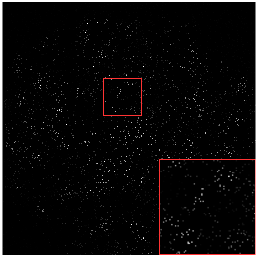}
			\includegraphics[width=1\linewidth]{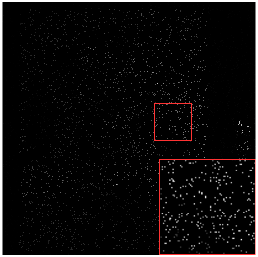}
	\end{minipage}}
	\subfigure[]{
		\begin{minipage}{0.056\textwidth}
			\includegraphics[width=1\linewidth]{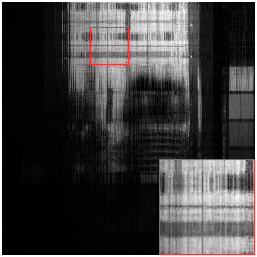}
			\includegraphics[width=1\linewidth]{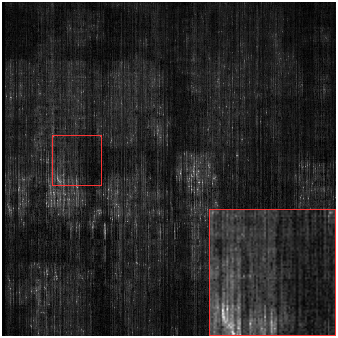}
			\includegraphics[width=1\linewidth]{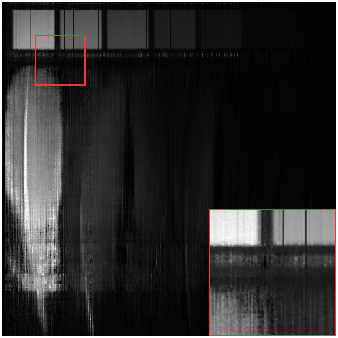}
			\includegraphics[width=1\linewidth]{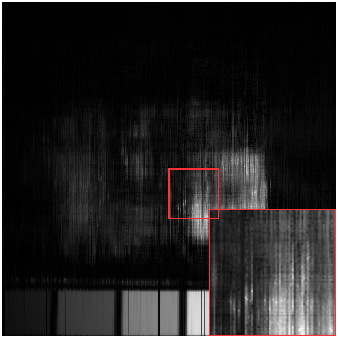}
			\includegraphics[width=1\linewidth]{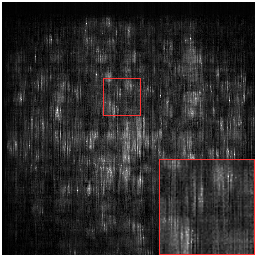}
			\includegraphics[width=1\linewidth]{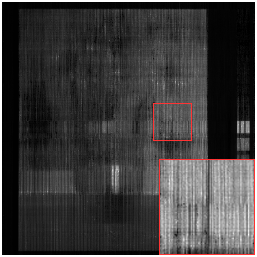}
	\end{minipage}}
	\subfigure[]{
		\begin{minipage}{0.056\textwidth}
			\includegraphics[width=1\linewidth]{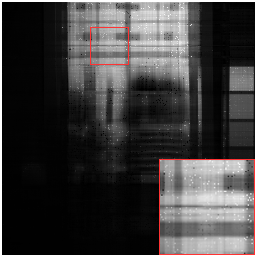}
			\includegraphics[width=1\linewidth]{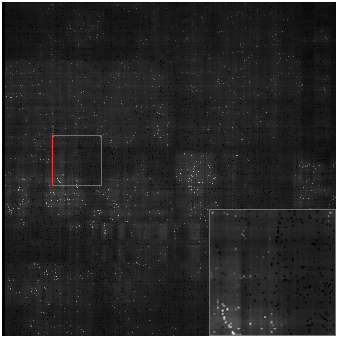}
			\includegraphics[width=1\linewidth]{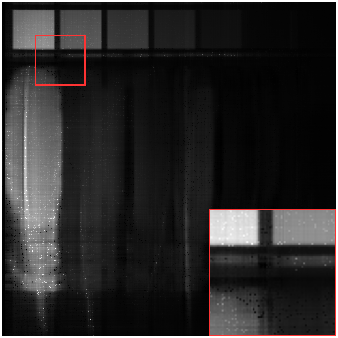}
			\includegraphics[width=1\linewidth]{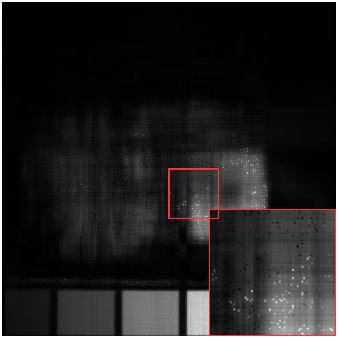}
			\includegraphics[width=1\linewidth]{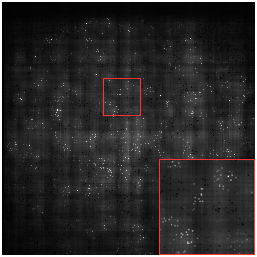}
			\includegraphics[width=1\linewidth]{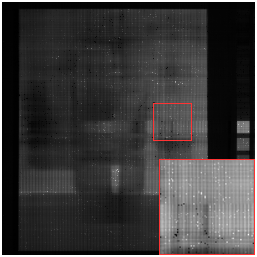}
	\end{minipage}}
	\subfigure[]{
		\begin{minipage}{0.056\textwidth}
			\includegraphics[width=1\linewidth]{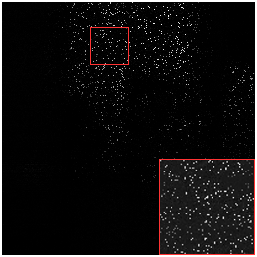}
			\includegraphics[width=1\linewidth]{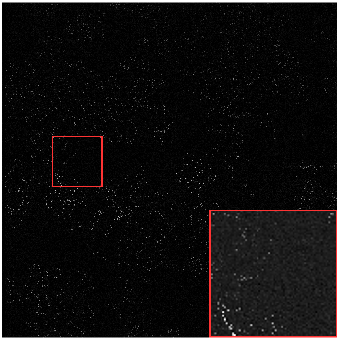}
			\includegraphics[width=1\linewidth]{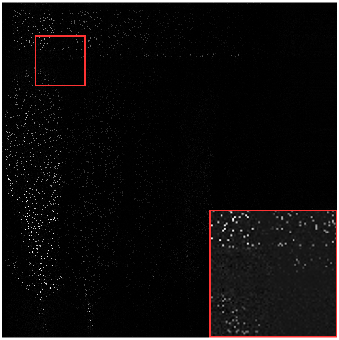}
			\includegraphics[width=1\linewidth]{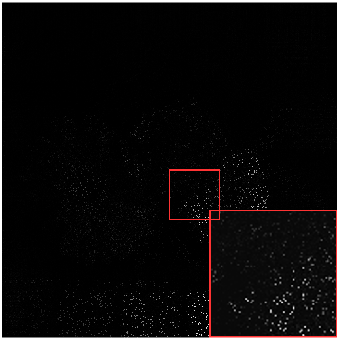}
			\includegraphics[width=1\linewidth]{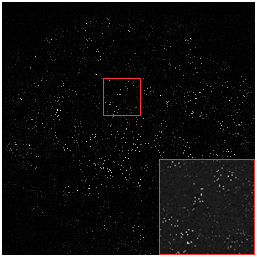}
			\includegraphics[width=1\linewidth]{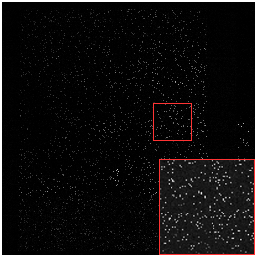}
	\end{minipage}}
	\subfigure[]{
		\begin{minipage}{0.056\textwidth}
			\includegraphics[width=1\linewidth]{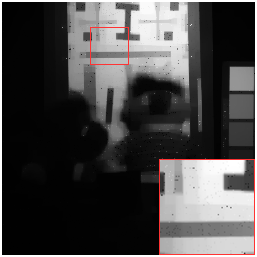}
			\includegraphics[width=1\linewidth]{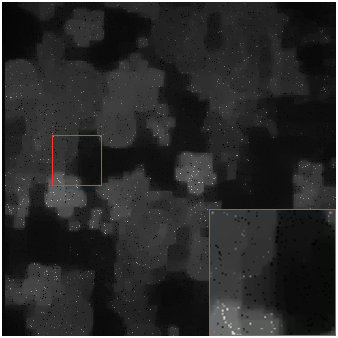}
			\includegraphics[width=1\linewidth]{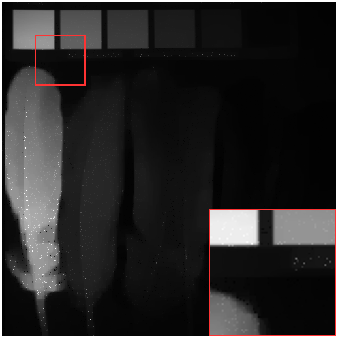}
			\includegraphics[width=1\linewidth]{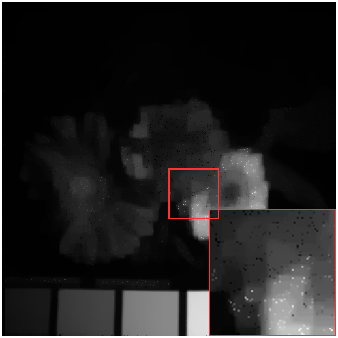}
			\includegraphics[width=1\linewidth]{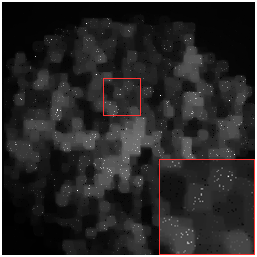}
			\includegraphics[width=1\linewidth]{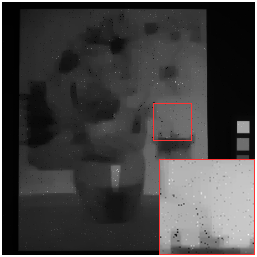}
	\end{minipage}}
	\subfigure[]{
		\begin{minipage}{0.056\textwidth}
			\includegraphics[width=1\linewidth]{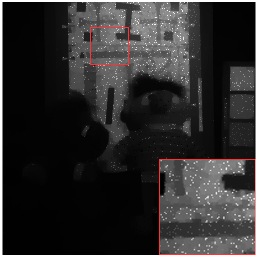}
			\includegraphics[width=1\linewidth]{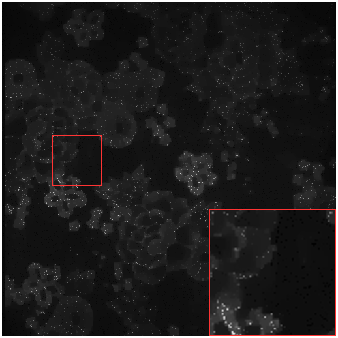}
			\includegraphics[width=1\linewidth]{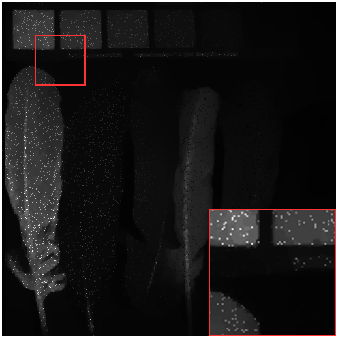}
			\includegraphics[width=1\linewidth]{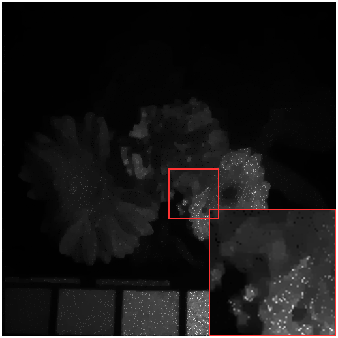}
			\includegraphics[width=1\linewidth]{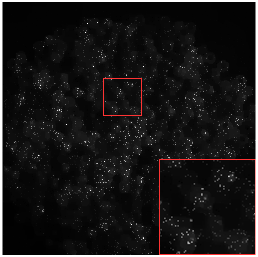}
			\includegraphics[width=1\linewidth]{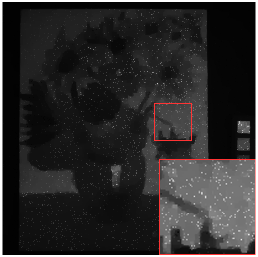}
	\end{minipage}}
	\subfigure[]{
		\begin{minipage}{0.056\textwidth}
			\includegraphics[width=1\linewidth]{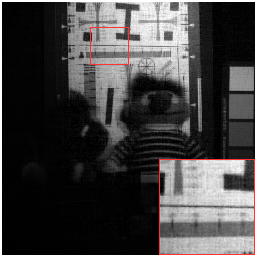}
			\includegraphics[width=1\linewidth]{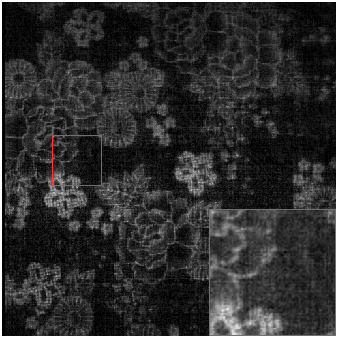}
			\includegraphics[width=1\linewidth]{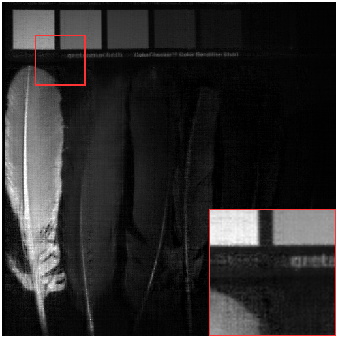}
			\includegraphics[width=1\linewidth]{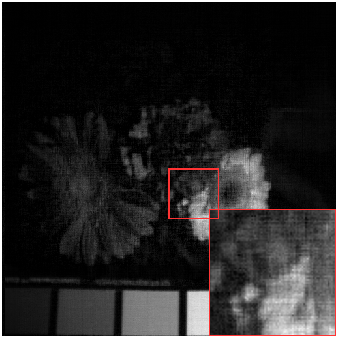}
			\includegraphics[width=1\linewidth]{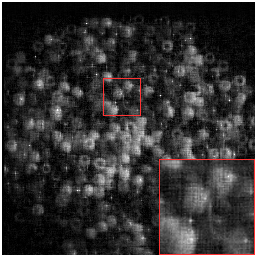}
			\includegraphics[width=1\linewidth]{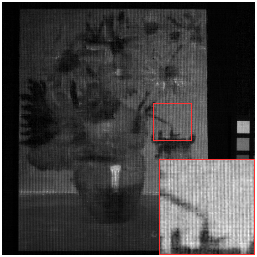}
	\end{minipage}}
	\subfigure[]{
		\begin{minipage}{0.056\textwidth}
			\includegraphics[width=1\linewidth]{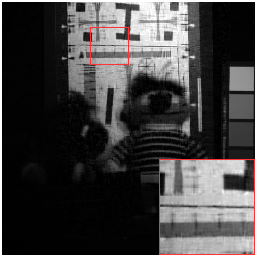}
			\includegraphics[width=1\linewidth]{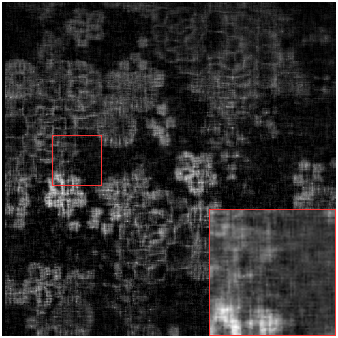}
			\includegraphics[width=1\linewidth]{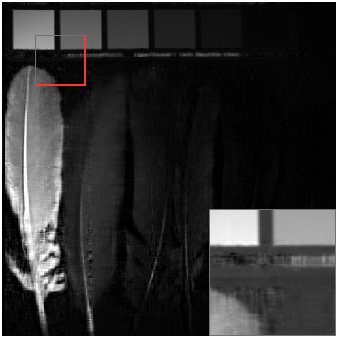}
			\includegraphics[width=1\linewidth]{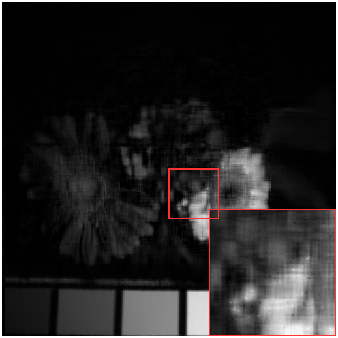}
			\includegraphics[width=1\linewidth]{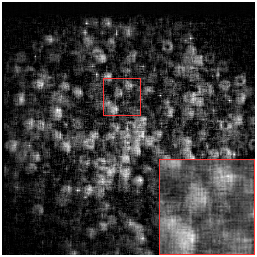}
			\includegraphics[width=1\linewidth]{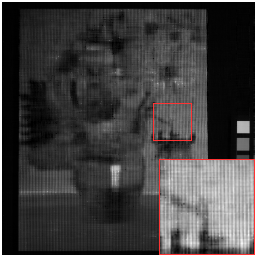}
	\end{minipage}}
	\subfigure[]{
		\begin{minipage}{0.056\textwidth}
			\includegraphics[width=1\linewidth]{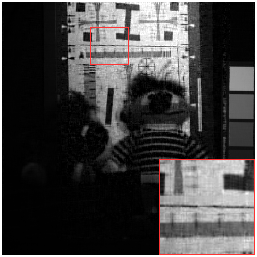}
			\includegraphics[width=1\linewidth]{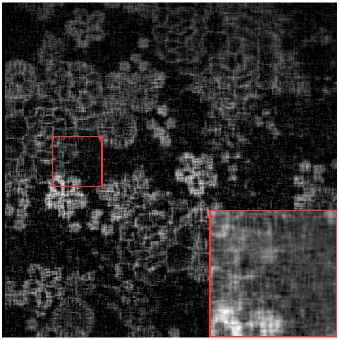}
			\includegraphics[width=1\linewidth]{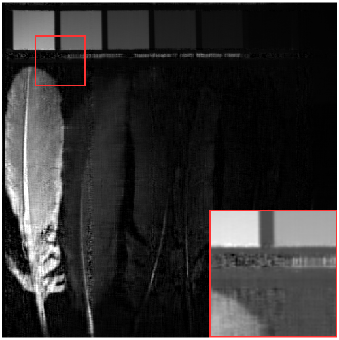}
			\includegraphics[width=1\linewidth]{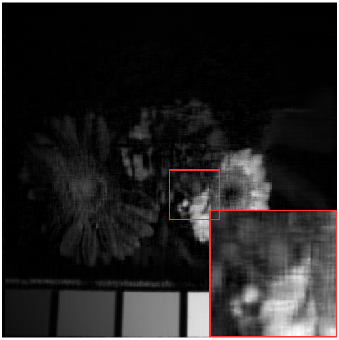}
			\includegraphics[width=1\linewidth]{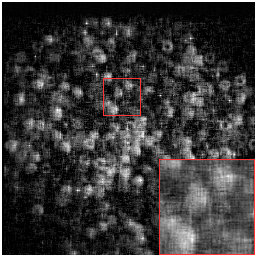}
			\includegraphics[width=1\linewidth]{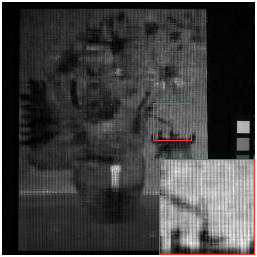}
	\end{minipage}}
	\subfigure[]{
		\begin{minipage}{0.056\textwidth}
			\includegraphics[width=1\linewidth]{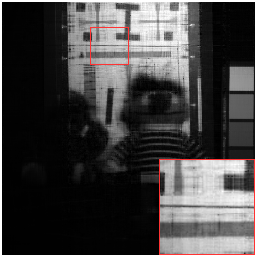}
			\includegraphics[width=1\linewidth]{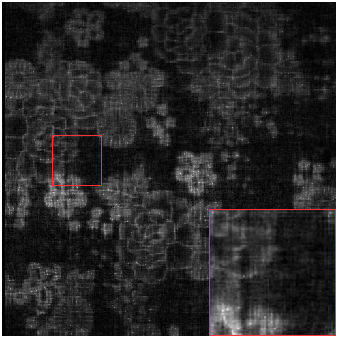}
			\includegraphics[width=1\linewidth]{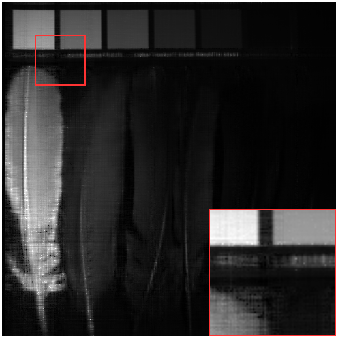}
			\includegraphics[width=1\linewidth]{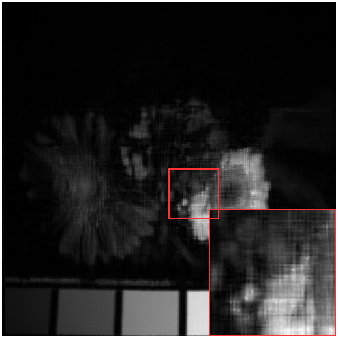}
			\includegraphics[width=1\linewidth]{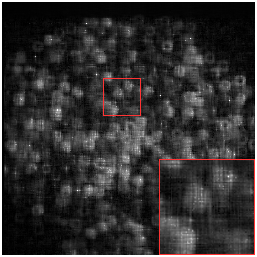}
			\includegraphics[width=1\linewidth]{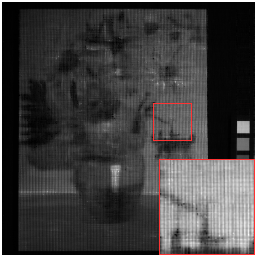}
	\end{minipage}}
	\subfigure[]{
		\begin{minipage}{0.056\textwidth}
			\includegraphics[width=1\linewidth]{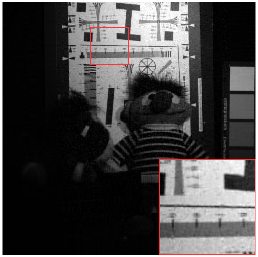}
			\includegraphics[width=1\linewidth]{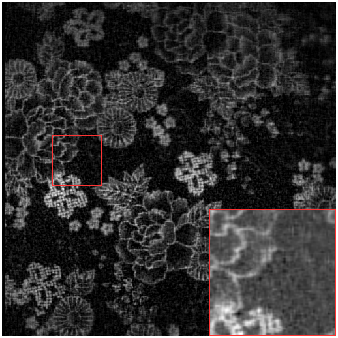}
			\includegraphics[width=1\linewidth]{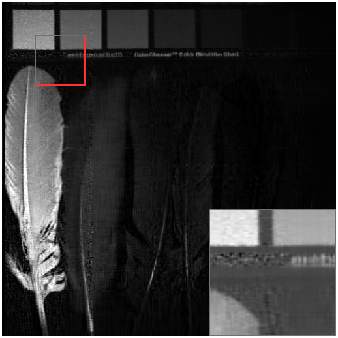}
			\includegraphics[width=1\linewidth]{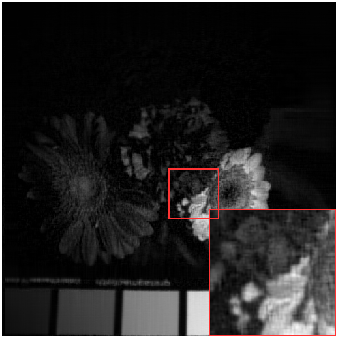}
			\includegraphics[width=1\linewidth]{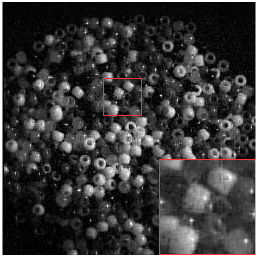}
			\includegraphics[width=1\linewidth]{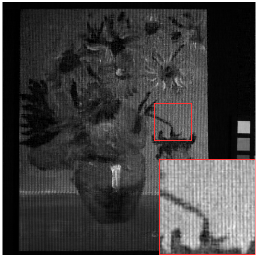}
	\end{minipage}}
	\subfigure[]{
		\begin{minipage}{0.056\textwidth}
			\includegraphics[width=1\linewidth]{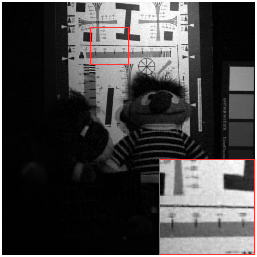}
			\includegraphics[width=1\linewidth]{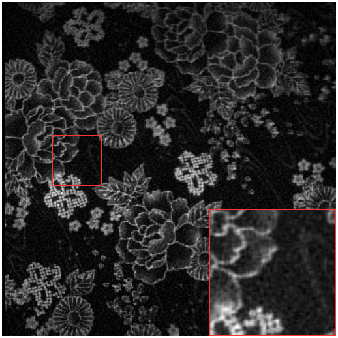}
			\includegraphics[width=1\linewidth]{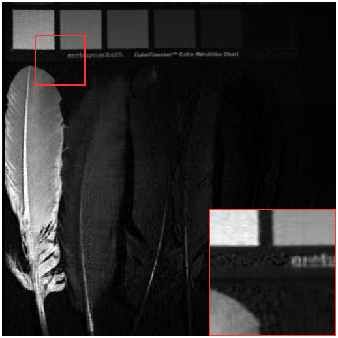}
			\includegraphics[width=1\linewidth]{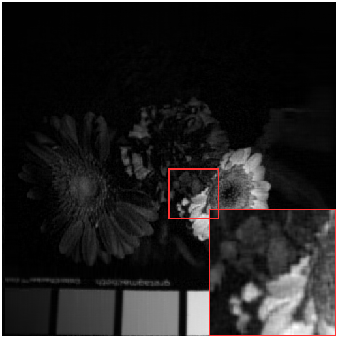}
			\includegraphics[width=1\linewidth]{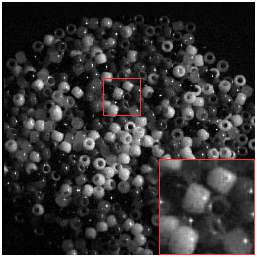}
			\includegraphics[width=1\linewidth]{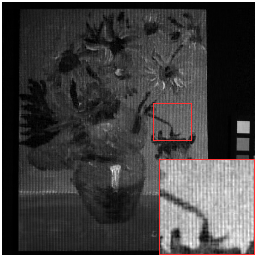}
	\end{minipage}}
	\caption{(a) Original images. (b) Corresponding sampled images with SR 5\%. (c)-(l) and (m) Completed images achieved by nine competing methods and proposed TLNM method and TLNMTV method, respectively. (a) Original image. (b) Observed image. (c) MC-ALM. (d) HaLRTC. (e) TMac. (f) LRTC-TV. (g) Trace/TV. (h) t-SVD. (i) McpTC. (j) ScadTC. (k) FTNN. (l) TLNM. (m) TLNMTV.}
	\label{HSITC5}
\end{figure*}
\begin{figure*}[!h]
	\centering
	\subfigure[]{
		\begin{minipage}{0.056\textwidth}
			\includegraphics[width=1\linewidth]{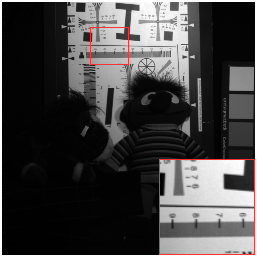}
			\includegraphics[width=1\linewidth]{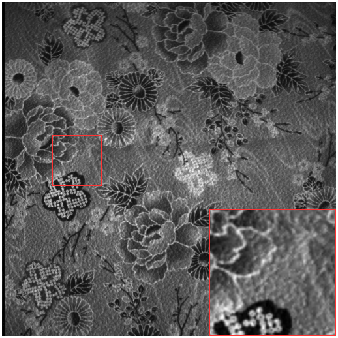}
			\includegraphics[width=1\linewidth]{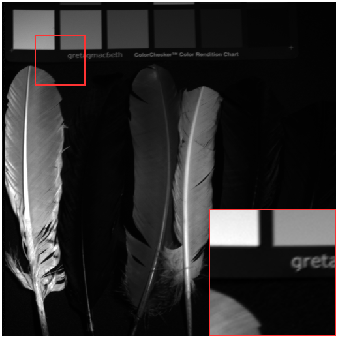}
			\includegraphics[width=1\linewidth]{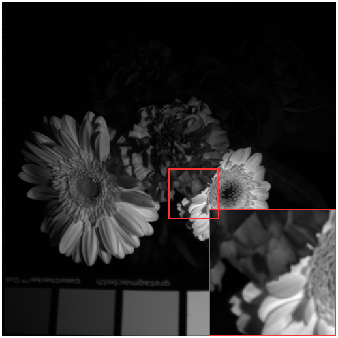}
			\includegraphics[width=1\linewidth]{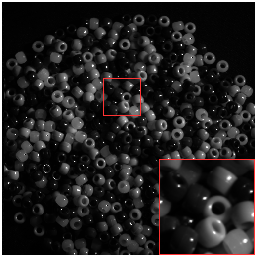}
			\includegraphics[width=1\linewidth]{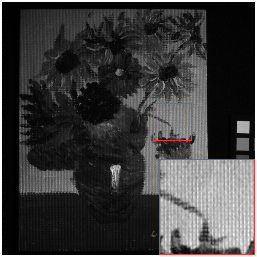}
	\end{minipage}}
	\subfigure[]{
		\begin{minipage}{0.056\textwidth}
			\includegraphics[width=1\linewidth]{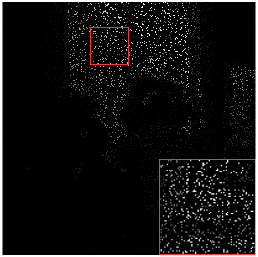}
			\includegraphics[width=1\linewidth]{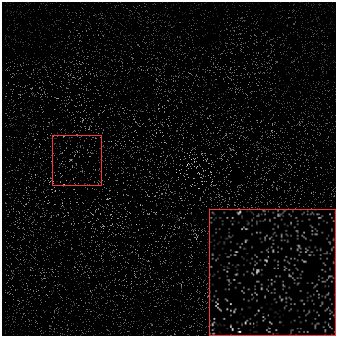}
			\includegraphics[width=1\linewidth]{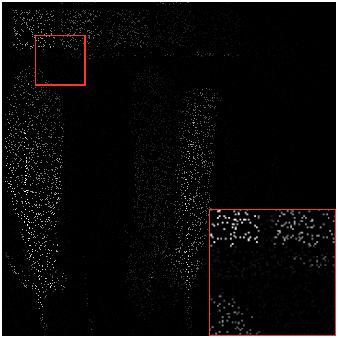}
			\includegraphics[width=1\linewidth]{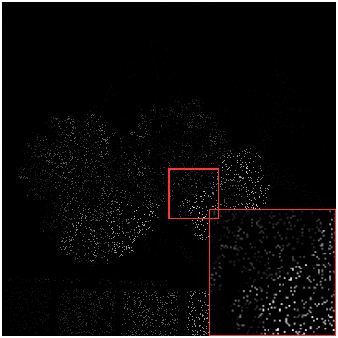}
			\includegraphics[width=1\linewidth]{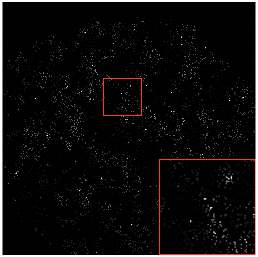}
			\includegraphics[width=1\linewidth]{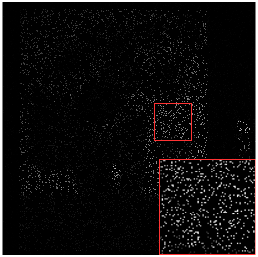}
	\end{minipage}}
	\subfigure[]{
		\begin{minipage}{0.056\textwidth}
			\includegraphics[width=1\linewidth]{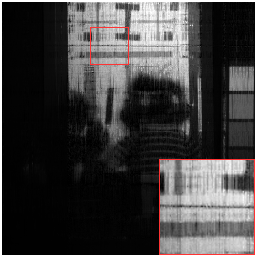}
			\includegraphics[width=1\linewidth]{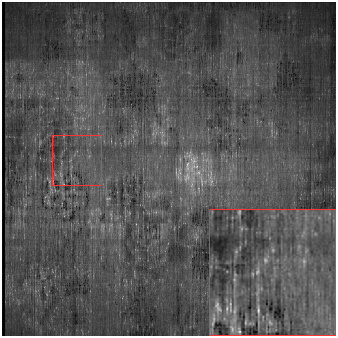}
			\includegraphics[width=1\linewidth]{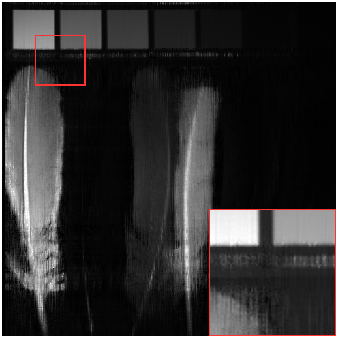}
			\includegraphics[width=1\linewidth]{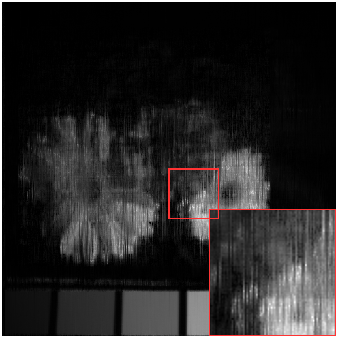}
			\includegraphics[width=1\linewidth]{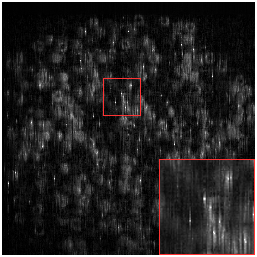}
			\includegraphics[width=1\linewidth]{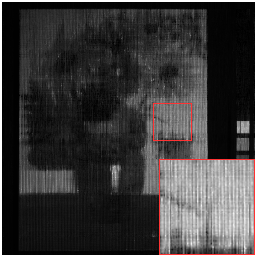}
	\end{minipage}}
	\subfigure[]{
		\begin{minipage}{0.056\textwidth}
			\includegraphics[width=1\linewidth]{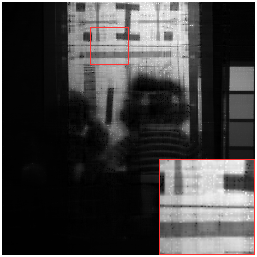}
			\includegraphics[width=1\linewidth]{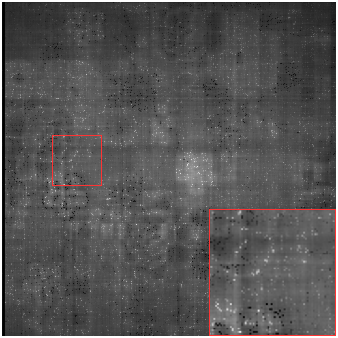}
			\includegraphics[width=1\linewidth]{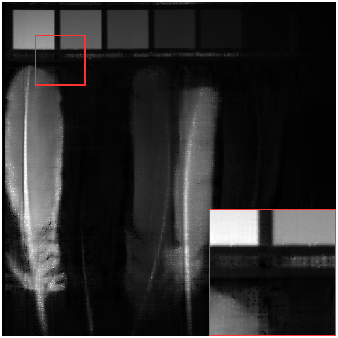}
			\includegraphics[width=1\linewidth]{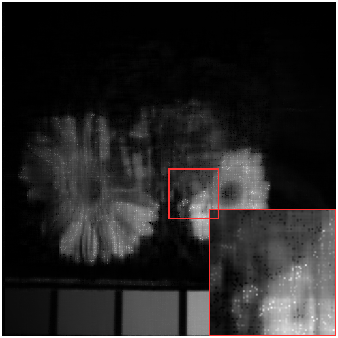}
			\includegraphics[width=1\linewidth]{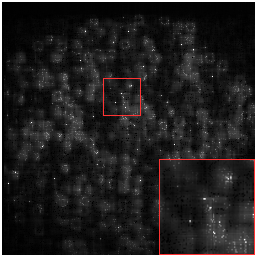}
			\includegraphics[width=1\linewidth]{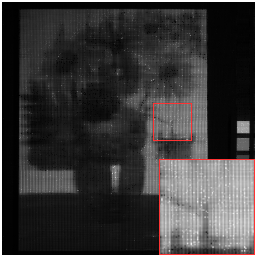}
	\end{minipage}}
	\subfigure[]{
		\begin{minipage}{0.056\textwidth}
			\includegraphics[width=1\linewidth]{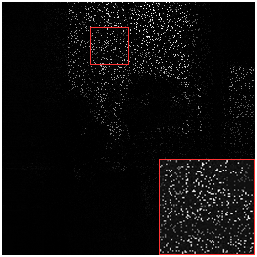}
			\includegraphics[width=1\linewidth]{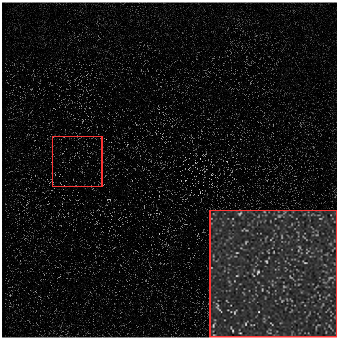}
			\includegraphics[width=1\linewidth]{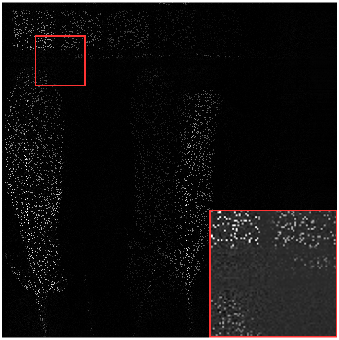}
			\includegraphics[width=1\linewidth]{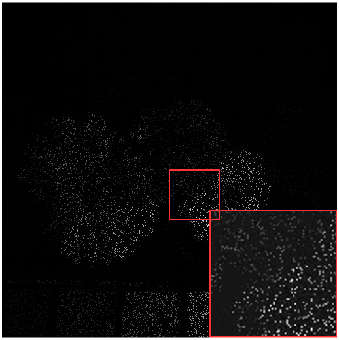}
			\includegraphics[width=1\linewidth]{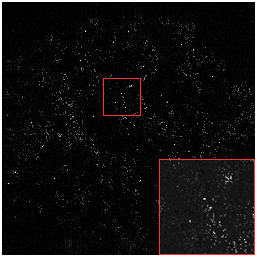}
			\includegraphics[width=1\linewidth]{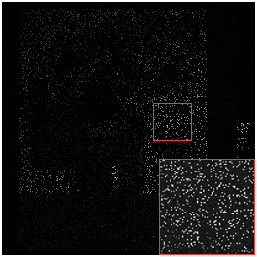}
	\end{minipage}}
	\subfigure[]{
		\begin{minipage}{0.056\textwidth}
			\includegraphics[width=1\linewidth]{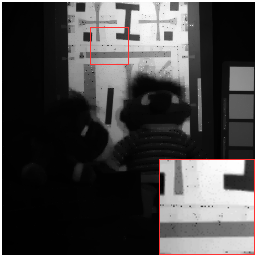}
			\includegraphics[width=1\linewidth]{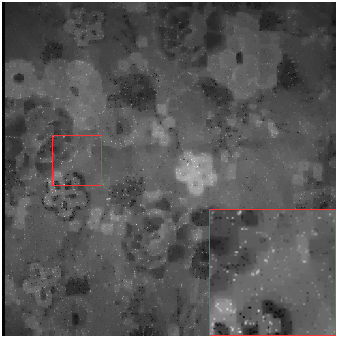}
			\includegraphics[width=1\linewidth]{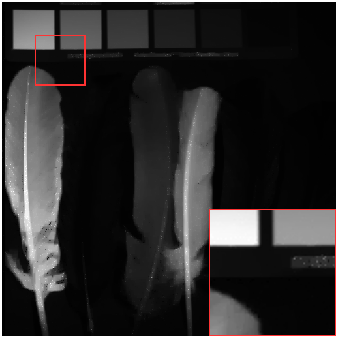}
			\includegraphics[width=1\linewidth]{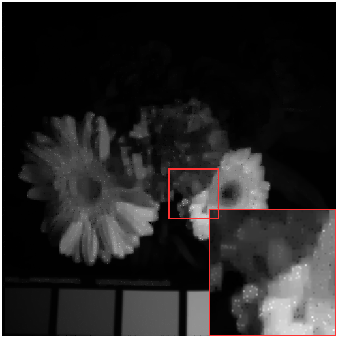}
			\includegraphics[width=1\linewidth]{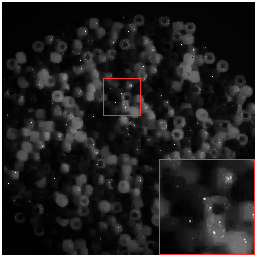}
			\includegraphics[width=1\linewidth]{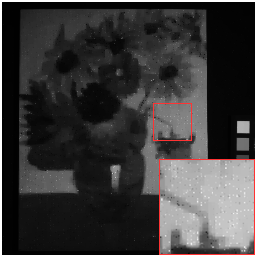}
	\end{minipage}}
	\subfigure[]{
		\begin{minipage}{0.056\textwidth}
			\includegraphics[width=1\linewidth]{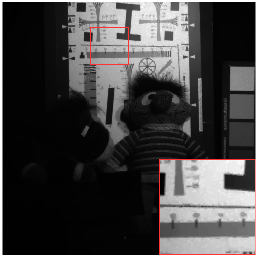}
			\includegraphics[width=1\linewidth]{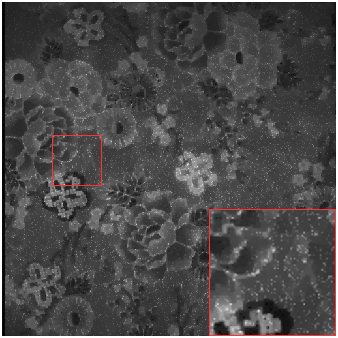}
			\includegraphics[width=1\linewidth]{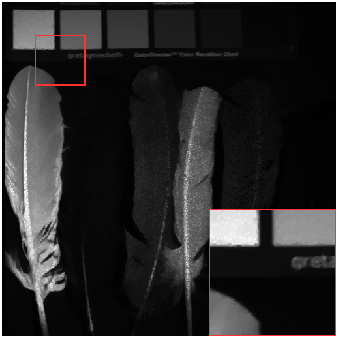}
			\includegraphics[width=1\linewidth]{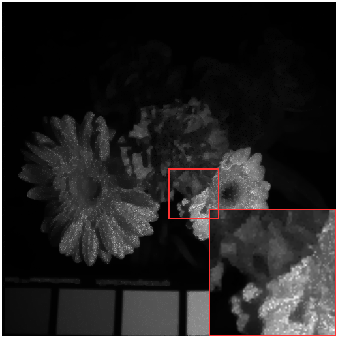}
			\includegraphics[width=1\linewidth]{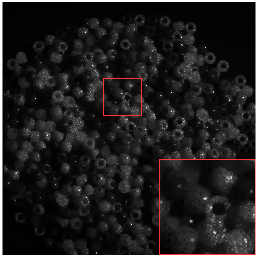}
			\includegraphics[width=1\linewidth]{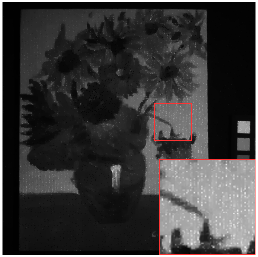}
	\end{minipage}}
	\subfigure[]{
		\begin{minipage}{0.056\textwidth}
			\includegraphics[width=1\linewidth]{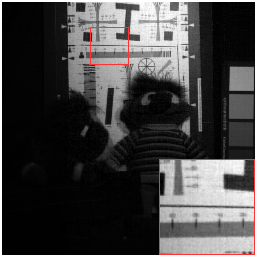}
			\includegraphics[width=1\linewidth]{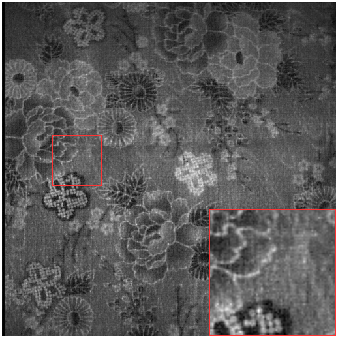}
			\includegraphics[width=1\linewidth]{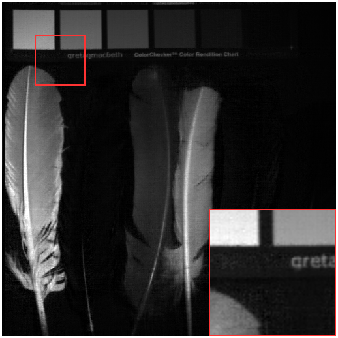}
			\includegraphics[width=1\linewidth]{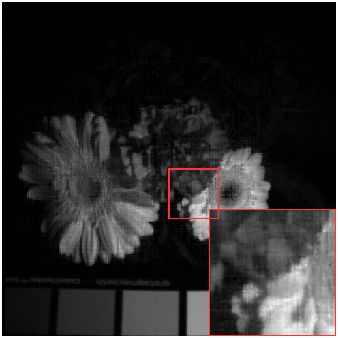}
			\includegraphics[width=1\linewidth]{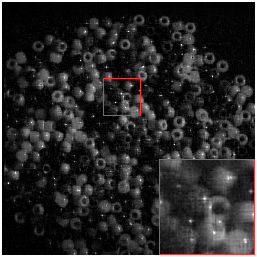}
			\includegraphics[width=1\linewidth]{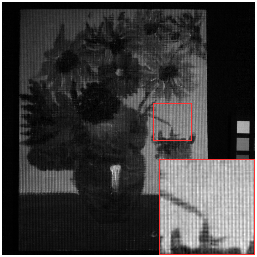}
	\end{minipage}}
	\subfigure[]{
		\begin{minipage}{0.056\textwidth}
			\includegraphics[width=1\linewidth]{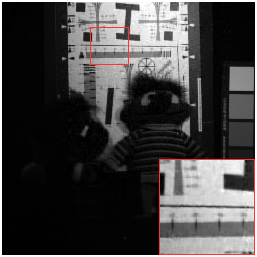}
			\includegraphics[width=1\linewidth]{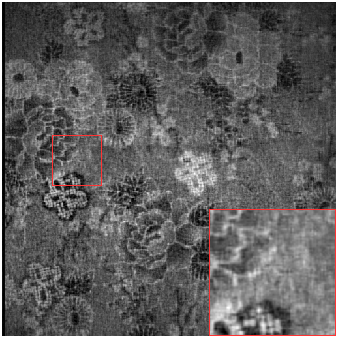}
			\includegraphics[width=1\linewidth]{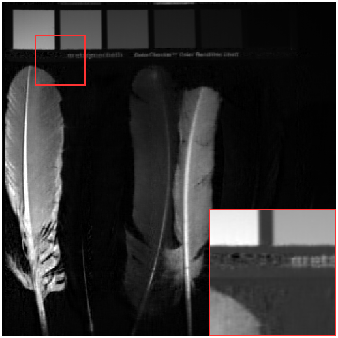}
			\includegraphics[width=1\linewidth]{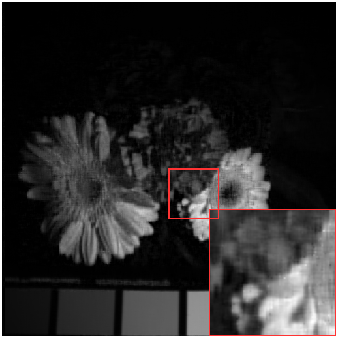}
			\includegraphics[width=1\linewidth]{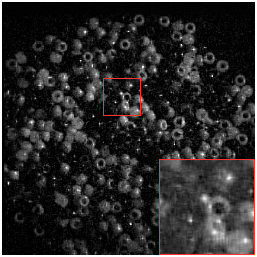}
			\includegraphics[width=1\linewidth]{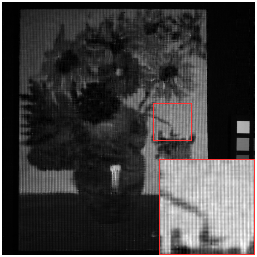}
	\end{minipage}}
	\subfigure[]{
		\begin{minipage}{0.056\textwidth}
			\includegraphics[width=1\linewidth]{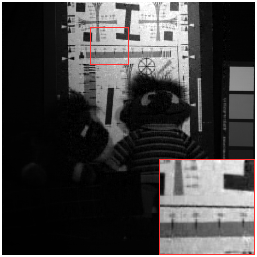}
			\includegraphics[width=1\linewidth]{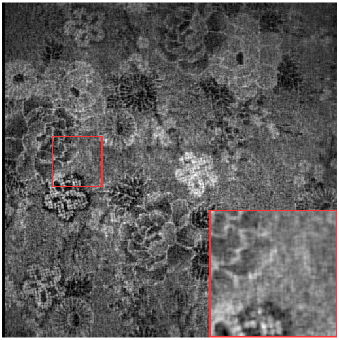}
			\includegraphics[width=1\linewidth]{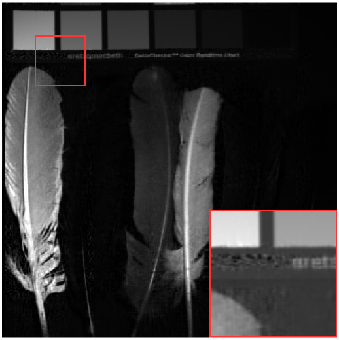}
			\includegraphics[width=1\linewidth]{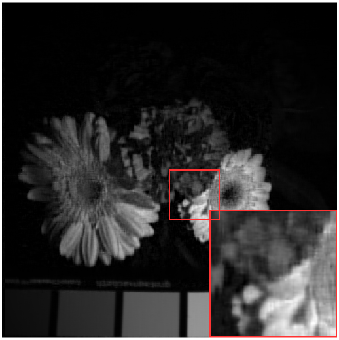}
			\includegraphics[width=1\linewidth]{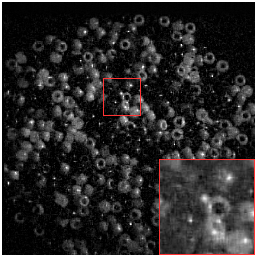}
			\includegraphics[width=1\linewidth]{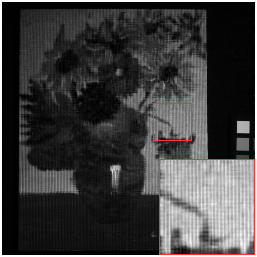}
	\end{minipage}}
	\subfigure[]{
		\begin{minipage}{0.056\textwidth}
			\includegraphics[width=1\linewidth]{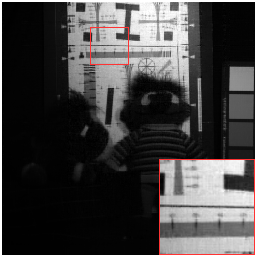}
			\includegraphics[width=1\linewidth]{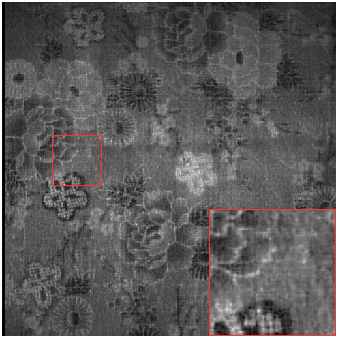}
			\includegraphics[width=1\linewidth]{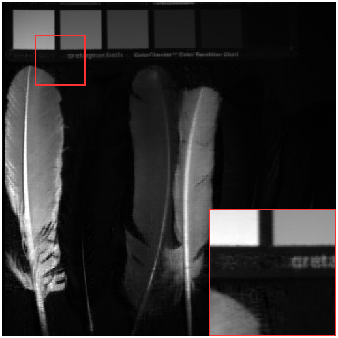}
			\includegraphics[width=1\linewidth]{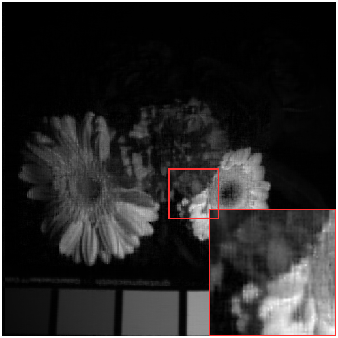}
			\includegraphics[width=1\linewidth]{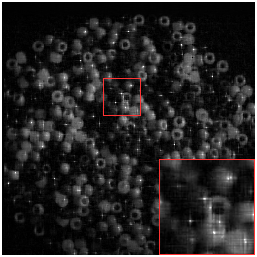}
			\includegraphics[width=1\linewidth]{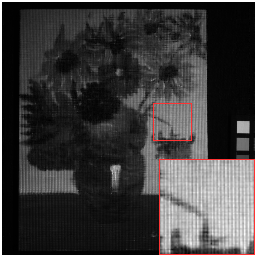}
	\end{minipage}}
	\subfigure[]{
		\begin{minipage}{0.056\textwidth}
			\includegraphics[width=1\linewidth]{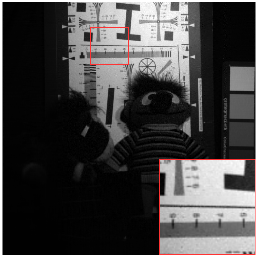}
			\includegraphics[width=1\linewidth]{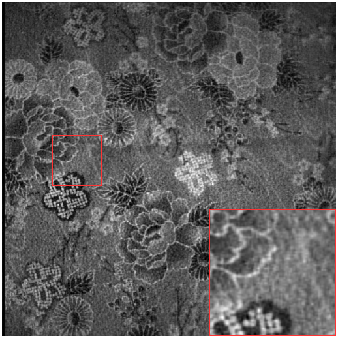}
			\includegraphics[width=1\linewidth]{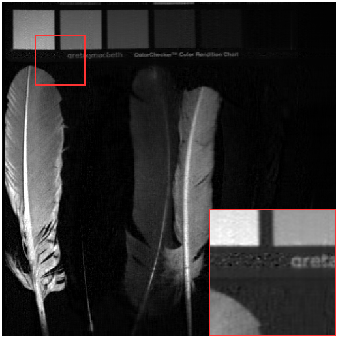}
			\includegraphics[width=1\linewidth]{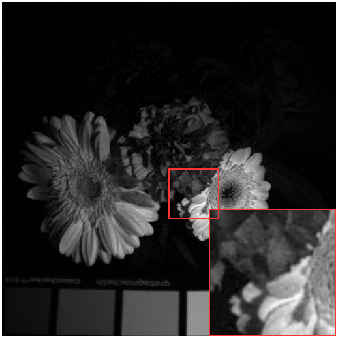}
			\includegraphics[width=1\linewidth]{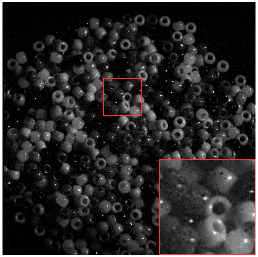}
			\includegraphics[width=1\linewidth]{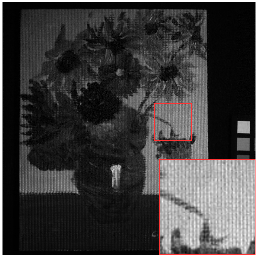}
	\end{minipage}}
	\subfigure[]{
		\begin{minipage}{0.056\textwidth}
			\includegraphics[width=1\linewidth]{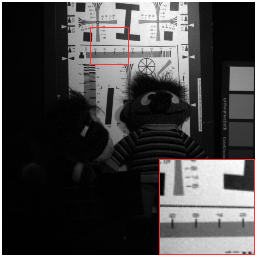}
			\includegraphics[width=1\linewidth]{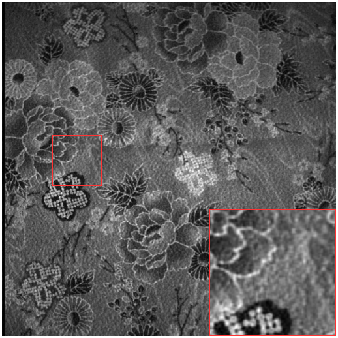}
			\includegraphics[width=1\linewidth]{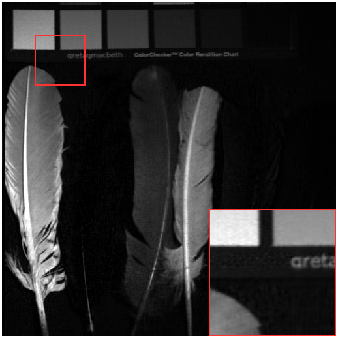}
			\includegraphics[width=1\linewidth]{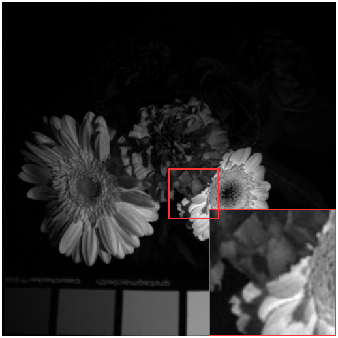}
			\includegraphics[width=1\linewidth]{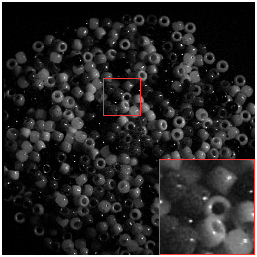}
			\includegraphics[width=1\linewidth]{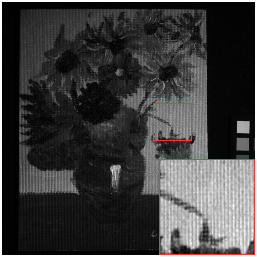}
	\end{minipage}}
	\caption{(a) Original images. (b) Corresponding sampled images with SR 10\%. (c)-(l) and (m) Completed images achieved by nine competing methods and proposed TLNM method and TLNMTV method, respectively. (a) Original image. (b) Observed image. (c) MC-ALM. (d) HaLRTC. (e) TMac. (f) LRTC-TV. (g) Trace/TV. (h) t-SVD. (i) McpTC. (j) ScadTC. (k) FTNN. (l) TLNM. (m) TLNMTV.}
	\label{HSITC10}
\end{figure*}
\begin{figure*}[!h]
	\centering
	\subfigure[]{
		\begin{minipage}{0.056\textwidth}
			\includegraphics[width=1\linewidth]{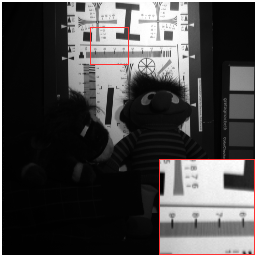}
			\includegraphics[width=1\linewidth]{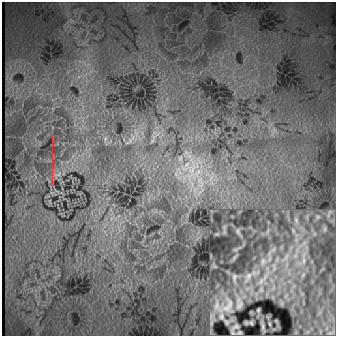}
			\includegraphics[width=1\linewidth]{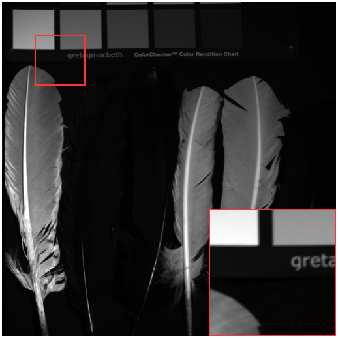}
			\includegraphics[width=1\linewidth]{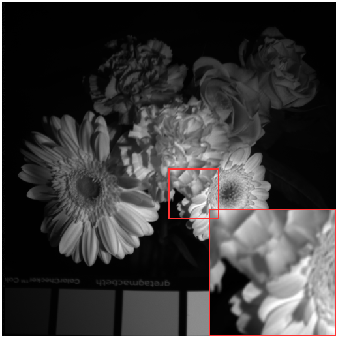}
			\includegraphics[width=1\linewidth]{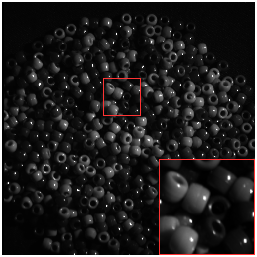}
			\includegraphics[width=1\linewidth]{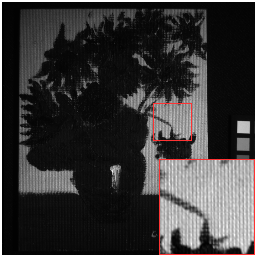}
	\end{minipage}}
	\subfigure[]{
		\begin{minipage}{0.056\textwidth}
			\includegraphics[width=1\linewidth]{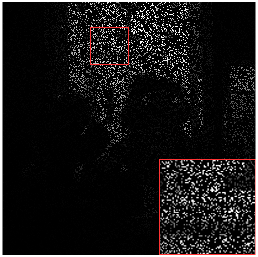}
			\includegraphics[width=1\linewidth]{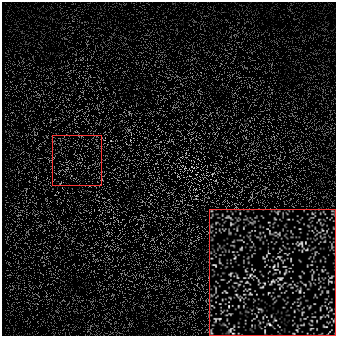}
			\includegraphics[width=1\linewidth]{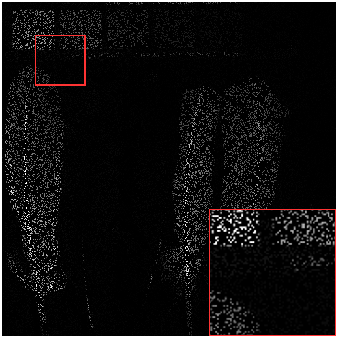}
			\includegraphics[width=1\linewidth]{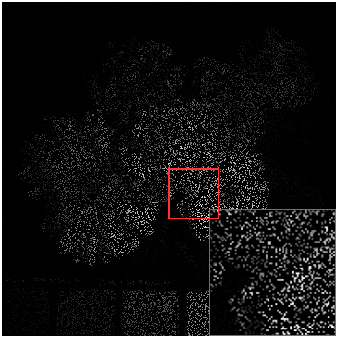}
			\includegraphics[width=1\linewidth]{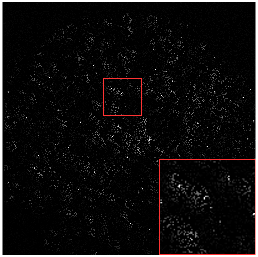}
			\includegraphics[width=1\linewidth]{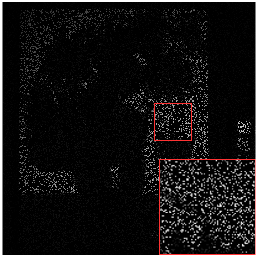}
	\end{minipage}}
	\subfigure[]{
		\begin{minipage}{0.056\textwidth}
			\includegraphics[width=1\linewidth]{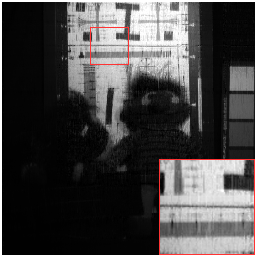}
			\includegraphics[width=1\linewidth]{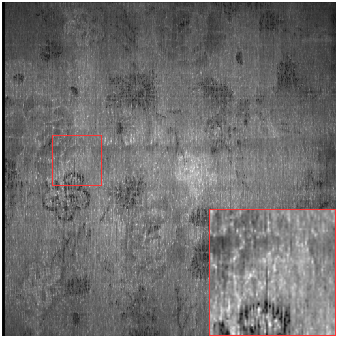}
			\includegraphics[width=1\linewidth]{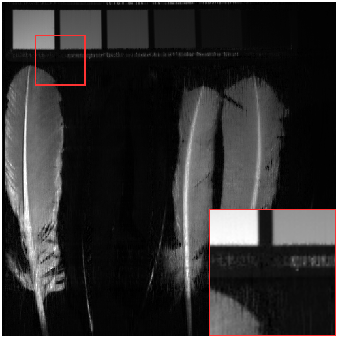}
			\includegraphics[width=1\linewidth]{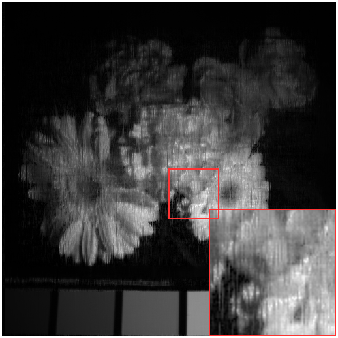}
			\includegraphics[width=1\linewidth]{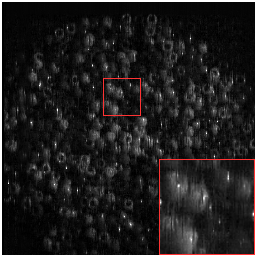}
			\includegraphics[width=1\linewidth]{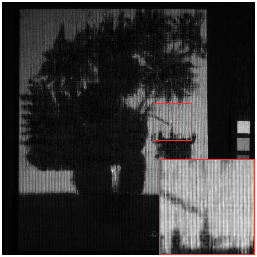}
	\end{minipage}}
	\subfigure[]{
		\begin{minipage}{0.056\textwidth}
			\includegraphics[width=1\linewidth]{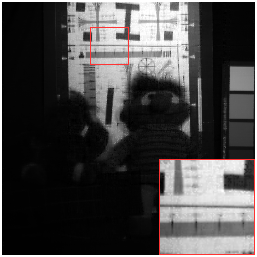}
			\includegraphics[width=1\linewidth]{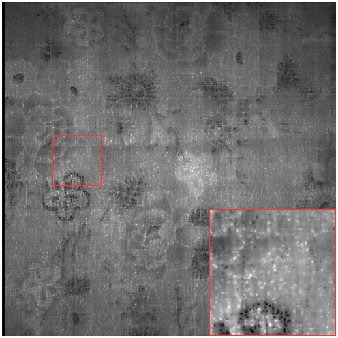}
			\includegraphics[width=1\linewidth]{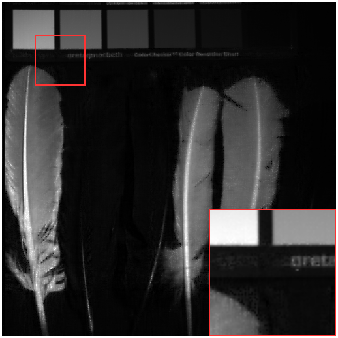}
			\includegraphics[width=1\linewidth]{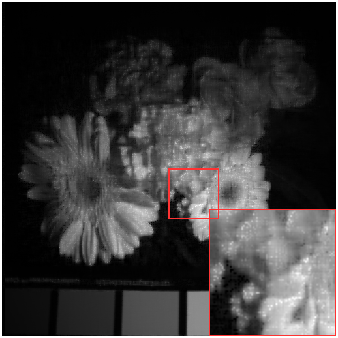}
			\includegraphics[width=1\linewidth]{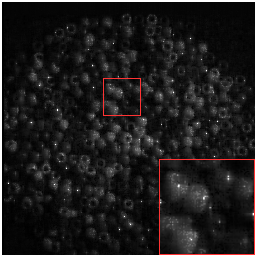}
			\includegraphics[width=1\linewidth]{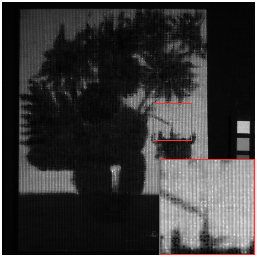}
	\end{minipage}}
	\subfigure[]{
		\begin{minipage}{0.056\textwidth}
			\includegraphics[width=1\linewidth]{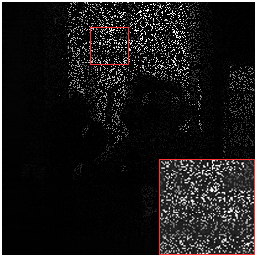}
			\includegraphics[width=1\linewidth]{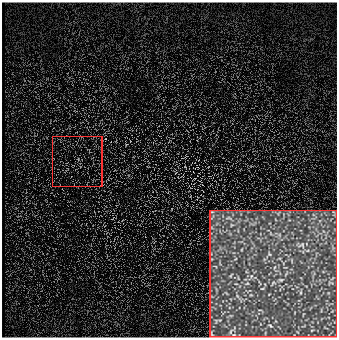}
			\includegraphics[width=1\linewidth]{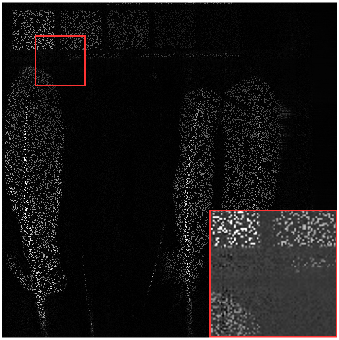}
			\includegraphics[width=1\linewidth]{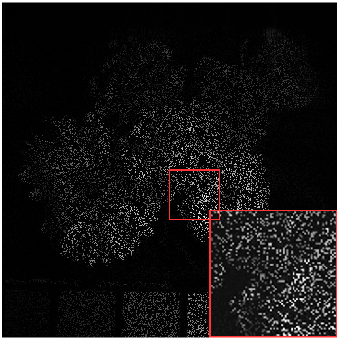}
			\includegraphics[width=1\linewidth]{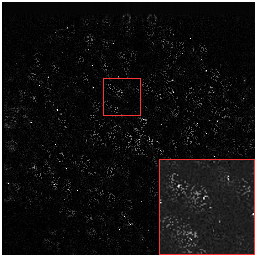}
			\includegraphics[width=1\linewidth]{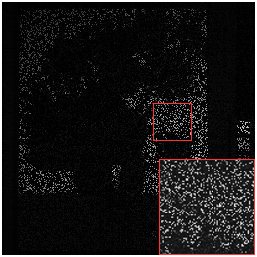}
	\end{minipage}}
	\subfigure[]{
		\begin{minipage}{0.056\textwidth}
			\includegraphics[width=1\linewidth]{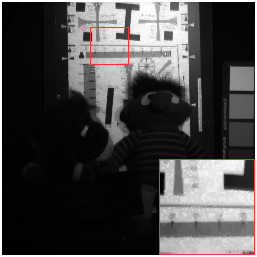}
			\includegraphics[width=1\linewidth]{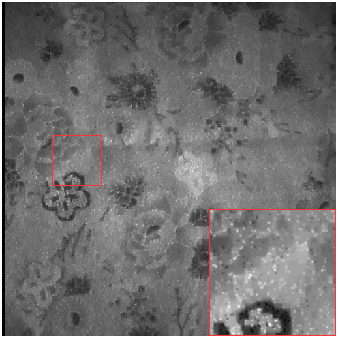}
			\includegraphics[width=1\linewidth]{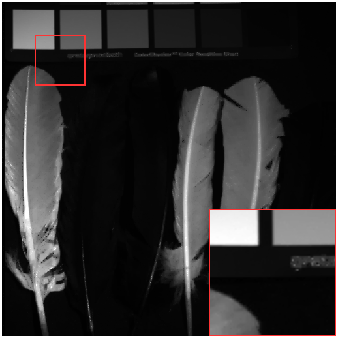}
			\includegraphics[width=1\linewidth]{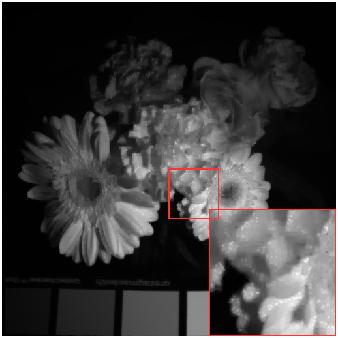}
			\includegraphics[width=1\linewidth]{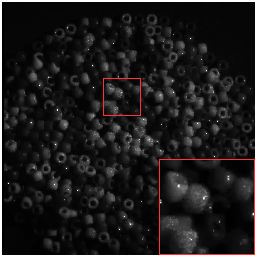}
			\includegraphics[width=1\linewidth]{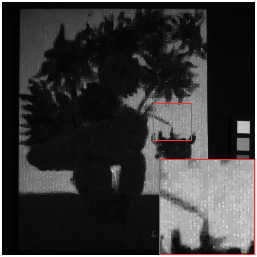}
	\end{minipage}}
	\subfigure[]{
		\begin{minipage}{0.056\textwidth}
			\includegraphics[width=1\linewidth]{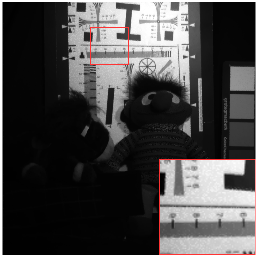}
			\includegraphics[width=1\linewidth]{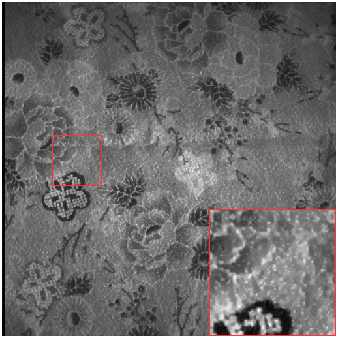}
			\includegraphics[width=1\linewidth]{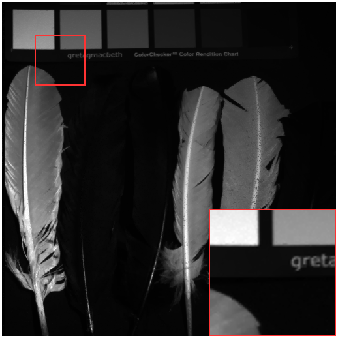}
			\includegraphics[width=1\linewidth]{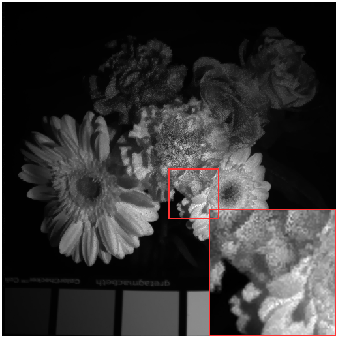}
			\includegraphics[width=1\linewidth]{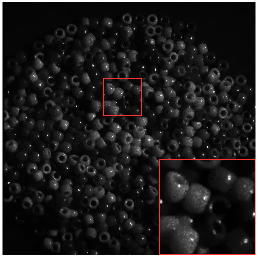}
			\includegraphics[width=1\linewidth]{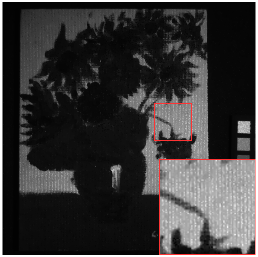}
	\end{minipage}}
	\subfigure[]{
		\begin{minipage}{0.056\textwidth}
			\includegraphics[width=1\linewidth]{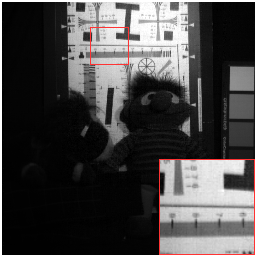}
			\includegraphics[width=1\linewidth]{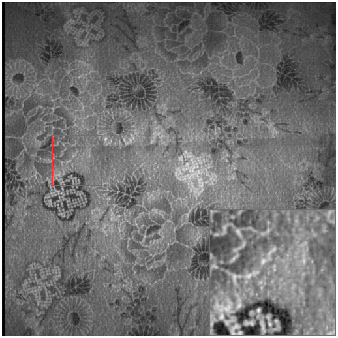}
			\includegraphics[width=1\linewidth]{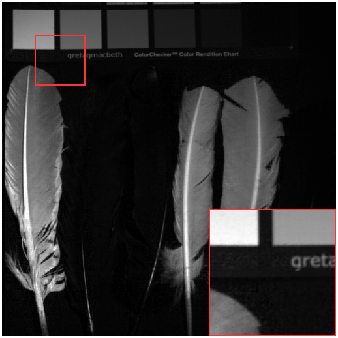}
			\includegraphics[width=1\linewidth]{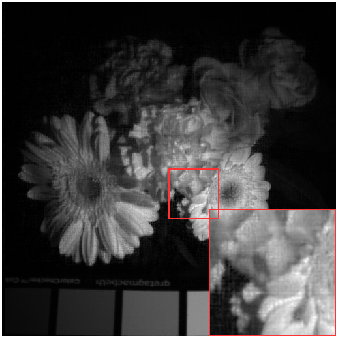}
			\includegraphics[width=1\linewidth]{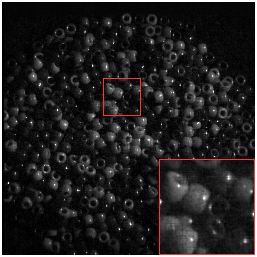}
			\includegraphics[width=1\linewidth]{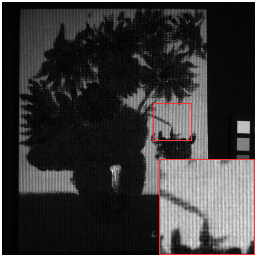}
	\end{minipage}}
	\subfigure[]{
		\begin{minipage}{0.056\textwidth}
			\includegraphics[width=1\linewidth]{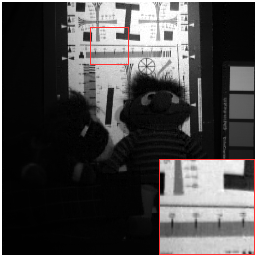}
			\includegraphics[width=1\linewidth]{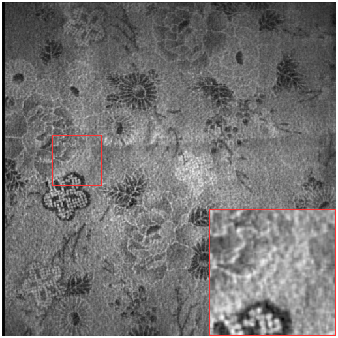}
			\includegraphics[width=1\linewidth]{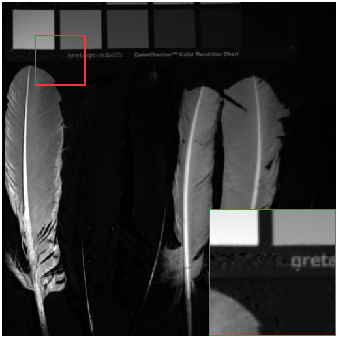}
			\includegraphics[width=1\linewidth]{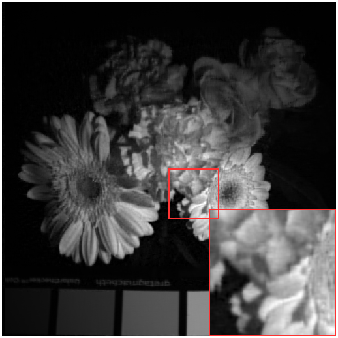}
			\includegraphics[width=1\linewidth]{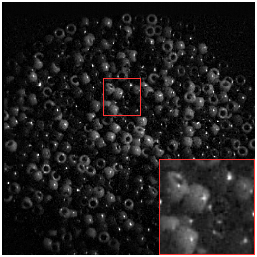}
			\includegraphics[width=1\linewidth]{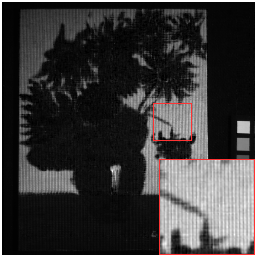}
	\end{minipage}}
	\subfigure[]{
		\begin{minipage}{0.056\textwidth}
			\includegraphics[width=1\linewidth]{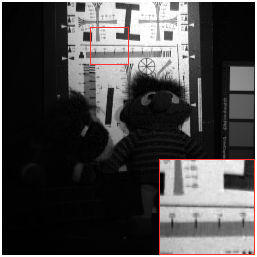}
			\includegraphics[width=1\linewidth]{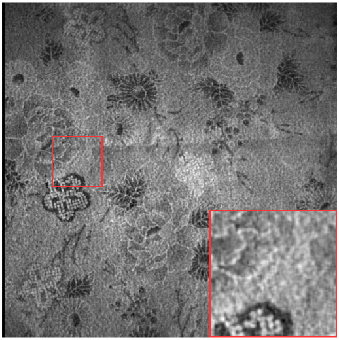}
			\includegraphics[width=1\linewidth]{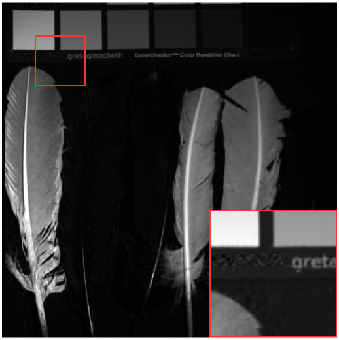}
			\includegraphics[width=1\linewidth]{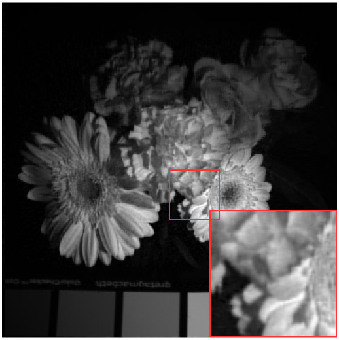}
			\includegraphics[width=1\linewidth]{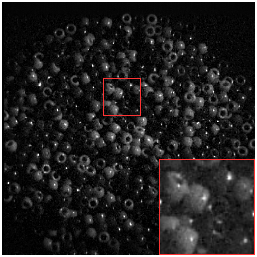}
			\includegraphics[width=1\linewidth]{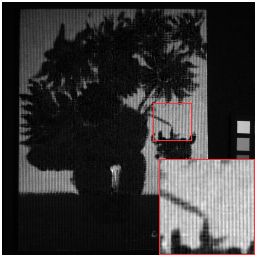}
	\end{minipage}}
	\subfigure[]{
		\begin{minipage}{0.056\textwidth}
			\includegraphics[width=1\linewidth]{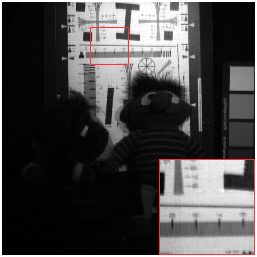}
			\includegraphics[width=1\linewidth]{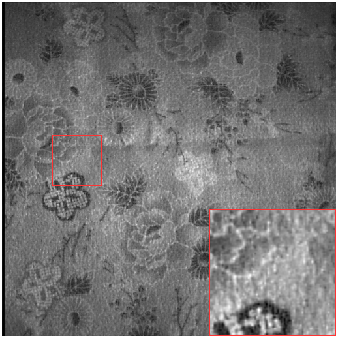}
			\includegraphics[width=1\linewidth]{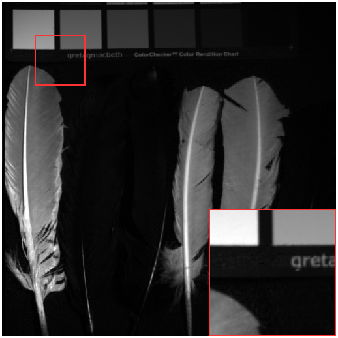}
			\includegraphics[width=1\linewidth]{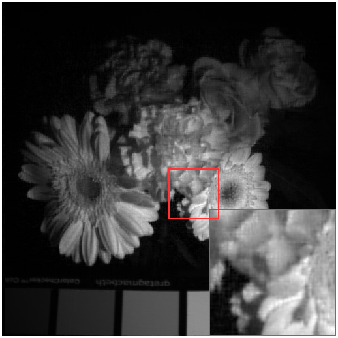}
			\includegraphics[width=1\linewidth]{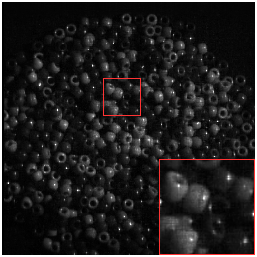}
			\includegraphics[width=1\linewidth]{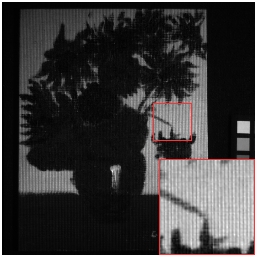}
	\end{minipage}}
	\subfigure[]{
		\begin{minipage}{0.056\textwidth}
			\includegraphics[width=1\linewidth]{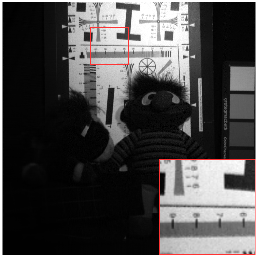}
			\includegraphics[width=1\linewidth]{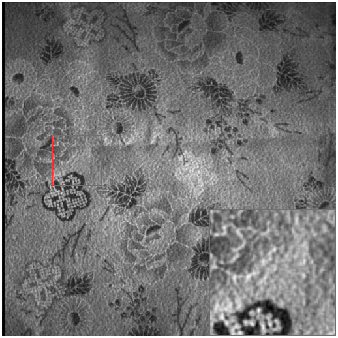}
			\includegraphics[width=1\linewidth]{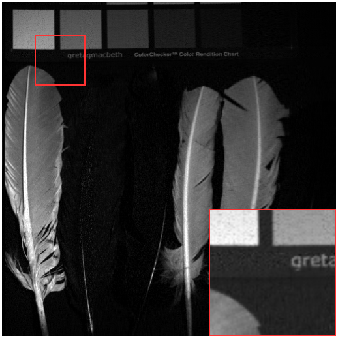}
			\includegraphics[width=1\linewidth]{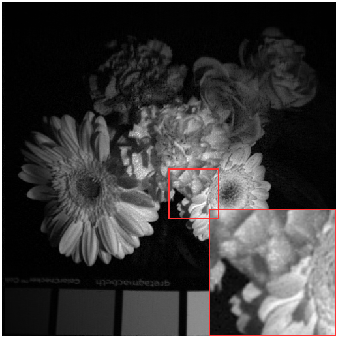}
			\includegraphics[width=1\linewidth]{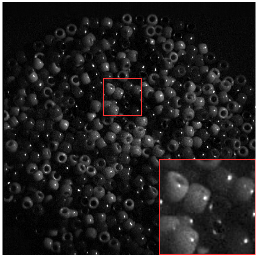}
			\includegraphics[width=1\linewidth]{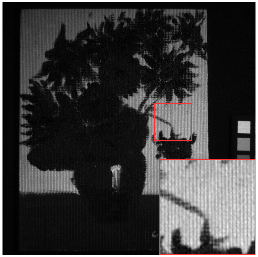}
	\end{minipage}}
	\subfigure[]{
		\begin{minipage}{0.056\textwidth}
			\includegraphics[width=1\linewidth]{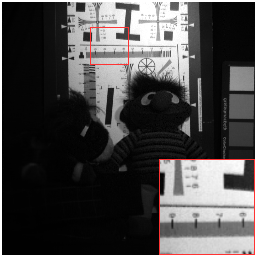}
			\includegraphics[width=1\linewidth]{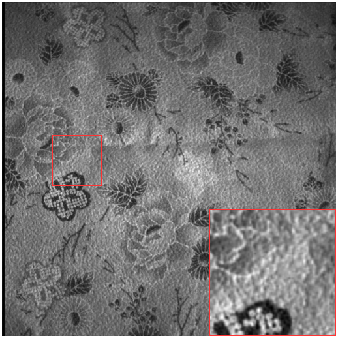}
			\includegraphics[width=1\linewidth]{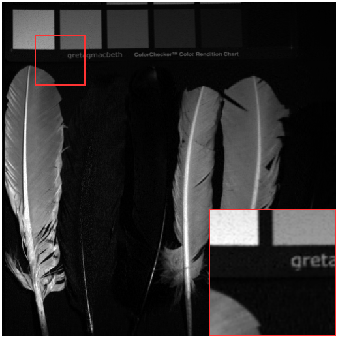}
			\includegraphics[width=1\linewidth]{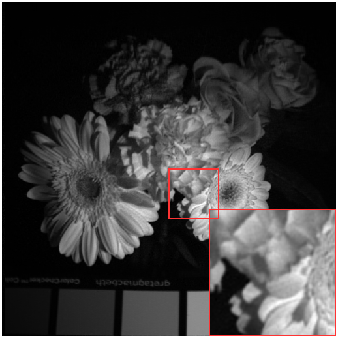}
			\includegraphics[width=1\linewidth]{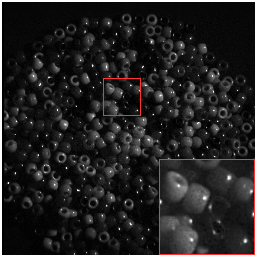}
			\includegraphics[width=1\linewidth]{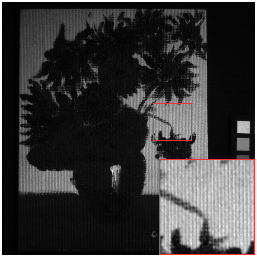}
	\end{minipage}}
	\caption{(a) Original images. (b) Corresponding sampled images with SR 20\%. (c)-(l) and (m) Completed images achieved by nine competing methods and proposed TLNM method and TLNMTV method, respectively. (a) Original image. (b) Observed image. (c) MC-ALM. (d) HaLRTC. (e) TMac. (f) LRTC-TV. (g) Trace/TV. (h) t-SVD. (i) McpTC. (j) ScadTC. (k) FTNN. (l) TLNM. (m) TLNMTV.}
	\label{HSITC20}
\end{figure*}
% Please add the following required packages to your document preamble:
% \usepackage{multirow}
\begin{table*}[!h]
	\caption{Average quantitative evaluation of the results for HSIs with different SRs.}
	\resizebox{\textwidth}{!}{
		\begin{tabular}{|c|c|c|c|c|c|c|c|c|c|c|c|c|c|}
			\hline
			SR                    & index & Observed & MC-ALM  & HaLRTC  & Tmac    & LRTC-TV & Trace-TV & t-SVD   & McpTC   & ScadTC  & FTNN    & TLNM    & TLNMTV  \\ \hline
			\multirow{4}{*}{5\%}  & PSNR  & 14.169   & 23.438  & 22.995  & 14.212  & 24.257  & 18.660   & 27.680  & 28.739  & 28.560  & 27.237  & 31.136  & 32.666  \\ \cline{2-14}
			& SSIM  & 0.203    & 0.593   & 0.650   & 0.226   & 0.709   & 0.621    & 0.727   & 0.725   & 0.701   & 0.763   & 0.790   & 0.857   \\ \cline{2-14}
			& FSIM  & 0.644    & 0.797   & 0.777   & 0.658   & 0.796   & 0.775    & 0.898   & 0.888   & 0.886   & 0.882   & 0.936   & 0.959   \\ \cline{2-14}
			& ERGAS & 910.569  & 332.274 & 350.889 & 906.147 & 308.510 & 562.121  & 210.364 & 195.615 & 199.534 & 221.914 & 139.178 & 116.057 \\ \hline
			\multirow{4}{*}{10\%} & PSNR  & 14.403   & 25.976  & 26.303  & 14.486  & 27.901  & 26.580   & 31.374  & 31.936  & 31.881  & 31.614  & 33.787  & 35.632  \\ \cline{2-14}
			& SSIM  & 0.238    & 0.694   & 0.745   & 0.269   & 0.817   & 0.812    & 0.841   & 0.823   & 0.811   & 0.876   & 0.858   & 0.910   \\ \cline{2-14}
			& FSIM  & 0.688    & 0.861   & 0.856   & 0.696   & 0.894   & 0.919    & 0.949   & 0.937   & 0.935   & 0.950   & 0.964   & 0.979   \\ \cline{2-14}
			& ERGAS & 886.354  & 252.090 & 246.122 & 878.174 & 204.648 & 247.447  & 139.290 & 134.702 & 137.803 & 134.651 & 105.123 & 85.385  \\ \hline
			\multirow{4}{*}{20\%} & PSNR  & 14.915   & 29.336  & 30.575  & 15.262  & 35.053  & 35.053   & 36.054  & 35.728  & 35.705  & 37.290  & 36.500  & 38.578  \\ \cline{2-14}
			& SSIM  & 0.303    & 0.811   & 0.859   & 0.348   & 0.945   & 0.945    & 0.927   & 0.911   & 0.900   & 0.955   & 0.916   & 0.952   \\ \cline{2-14}
			& FSIM  & 0.734    & 0.924   & 0.936   & 0.744   & 0.984   & 0.984    & 0.980   & 0.973   & 0.971   & 0.986   & 0.983   & 0.989   \\ \cline{2-14}
			& ERGAS & 835.570  & 173.900 & 152.819 & 805.944 & 93.710  & 93.710   & 83.575  & 86.751  & 89.127  & 71.704  & 78.446  & 62.725  \\ \hline
	\end{tabular}}\label{THSI}
\end{table*}

\begin{figure*}[!h] %??????????, ???????, ??????????????, ????
	\centering  %??????
	\vspace{0cm} %??????????
	\subfigure[]{
		\begin{minipage}{0.056\textwidth}
			\includegraphics[width=1\linewidth]{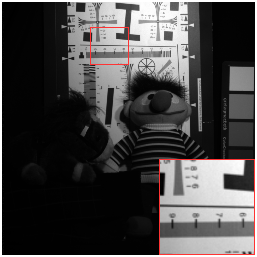}
			\includegraphics[width=1\linewidth]{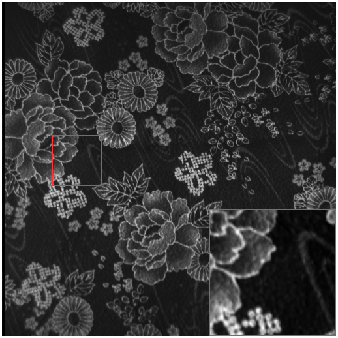}
			\includegraphics[width=1\linewidth]{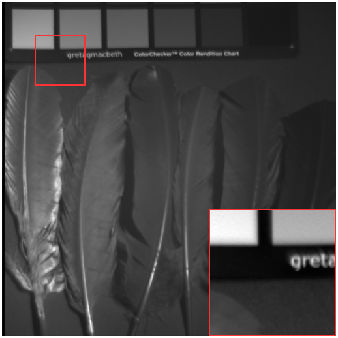}
			\includegraphics[width=1\linewidth]{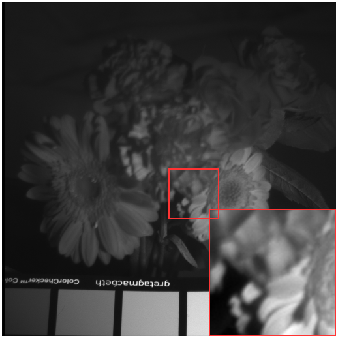}
			\includegraphics[width=1\linewidth]{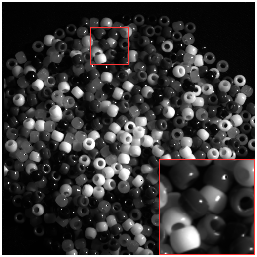}
			\includegraphics[width=1\linewidth]{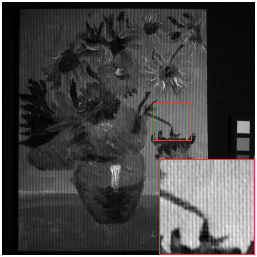}
	\end{minipage}}
	\subfigure[]{
		\begin{minipage}{0.056\textwidth}
			\includegraphics[width=1\linewidth]{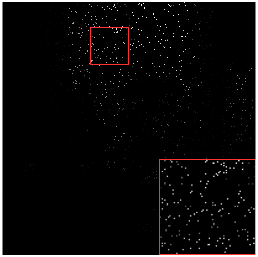}
			\includegraphics[width=1\linewidth]{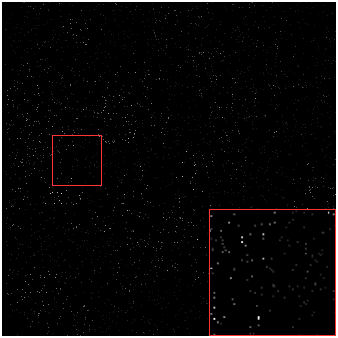}
			\includegraphics[width=1\linewidth]{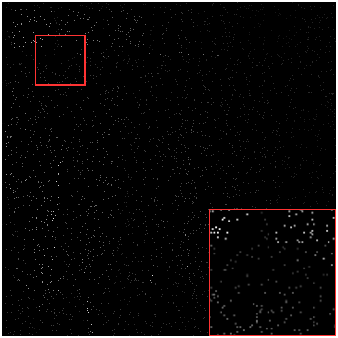}
			\includegraphics[width=1\linewidth]{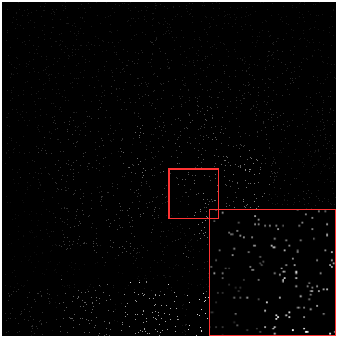}
			\includegraphics[width=1\linewidth]{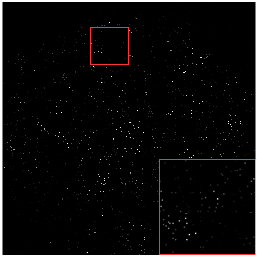}
			\includegraphics[width=1\linewidth]{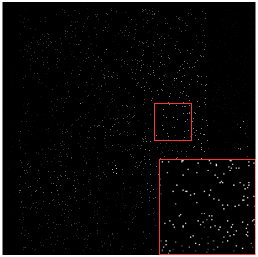}
	\end{minipage}}
	\subfigure[]{
		\begin{minipage}{0.056\textwidth}
			\includegraphics[width=1\linewidth]{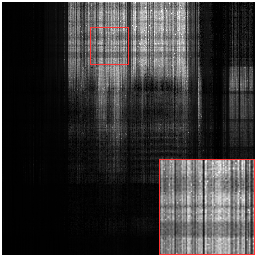}
			\includegraphics[width=1\linewidth]{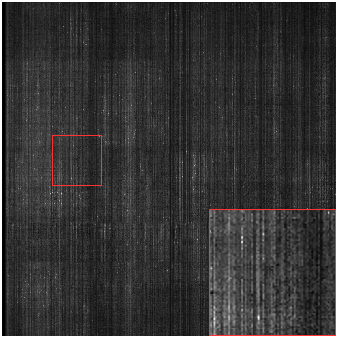}
			\includegraphics[width=1\linewidth]{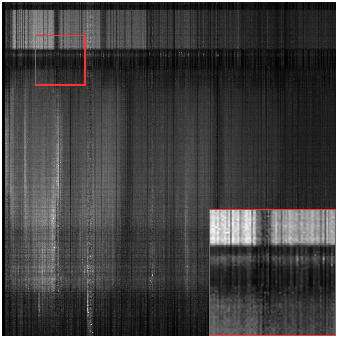}
			\includegraphics[width=1\linewidth]{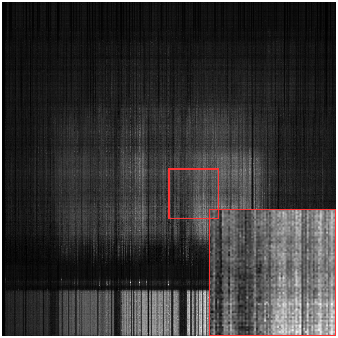}
			\includegraphics[width=1\linewidth]{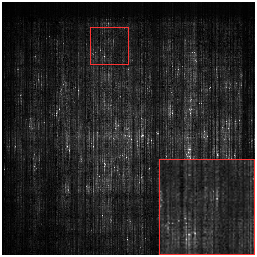}
			\includegraphics[width=1\linewidth]{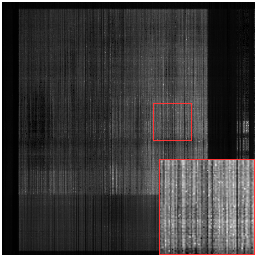}
	\end{minipage}}
	\subfigure[]{
		\begin{minipage}{0.056\textwidth}
			\includegraphics[width=1\linewidth]{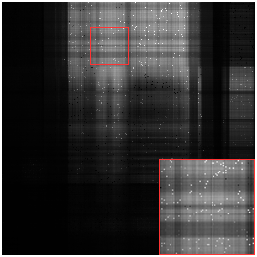}
			\includegraphics[width=1\linewidth]{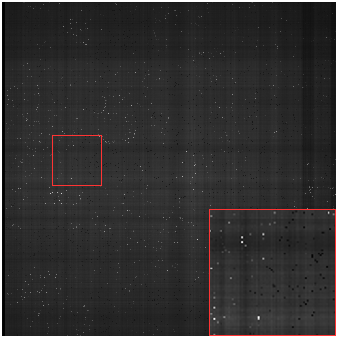}
			\includegraphics[width=1\linewidth]{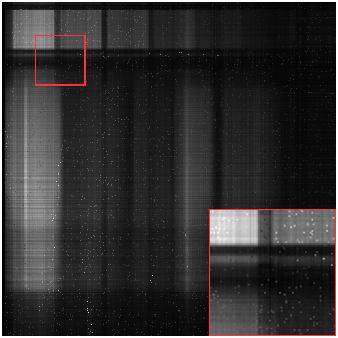}
			\includegraphics[width=1\linewidth]{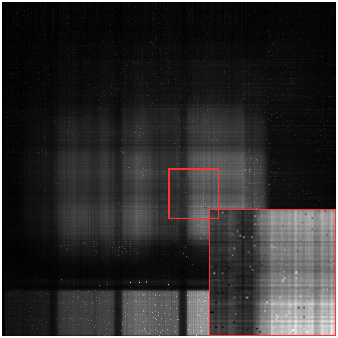}
			\includegraphics[width=1\linewidth]{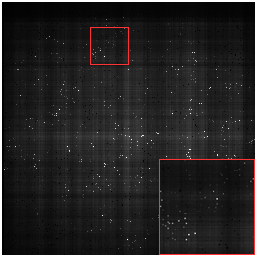}
			\includegraphics[width=1\linewidth]{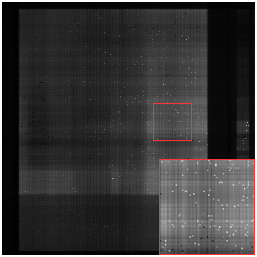}
	\end{minipage}}
	\subfigure[]{
		\begin{minipage}{0.056\textwidth}
			\includegraphics[width=1\linewidth]{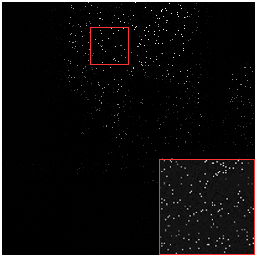}
			\includegraphics[width=1\linewidth]{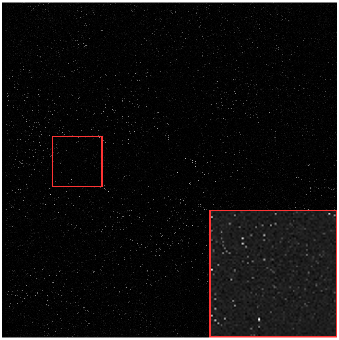}
			\includegraphics[width=1\linewidth]{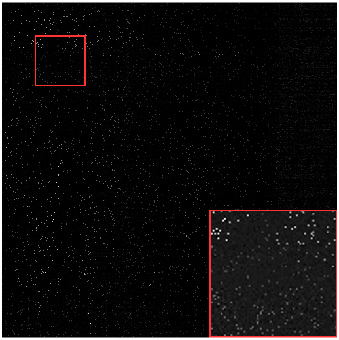}
			\includegraphics[width=1\linewidth]{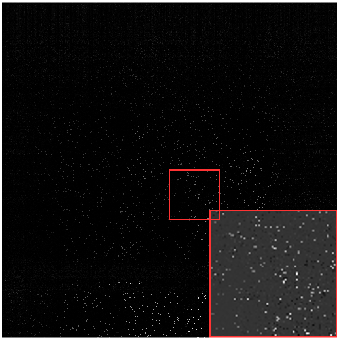}
			\includegraphics[width=1\linewidth]{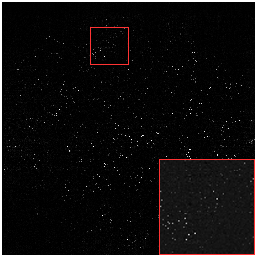}
			\includegraphics[width=1\linewidth]{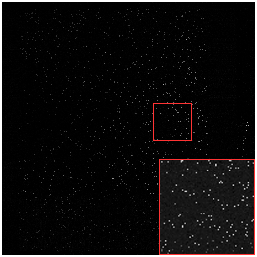}
	\end{minipage}}
	\subfigure[]{
		\begin{minipage}{0.056\textwidth}
			\includegraphics[width=1\linewidth]{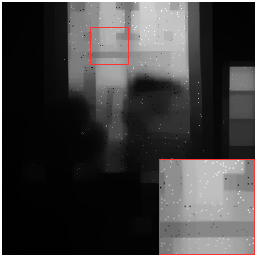}
			\includegraphics[width=1\linewidth]{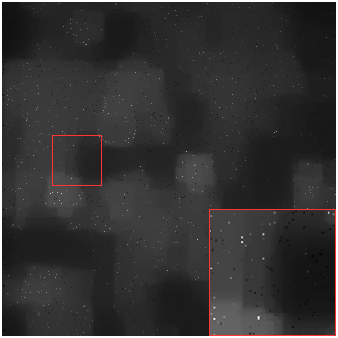}
			\includegraphics[width=1\linewidth]{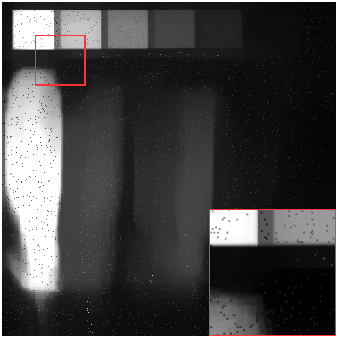}
			\includegraphics[width=1\linewidth]{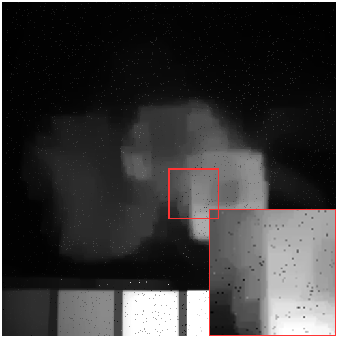}
			\includegraphics[width=1\linewidth]{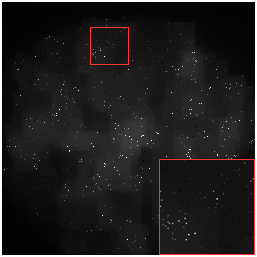}
			\includegraphics[width=1\linewidth]{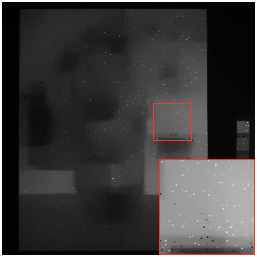}
	\end{minipage}}
	\subfigure[]{
		\begin{minipage}{0.056\textwidth}
			\includegraphics[width=1\linewidth]{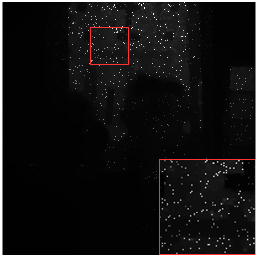}
			\includegraphics[width=1\linewidth]{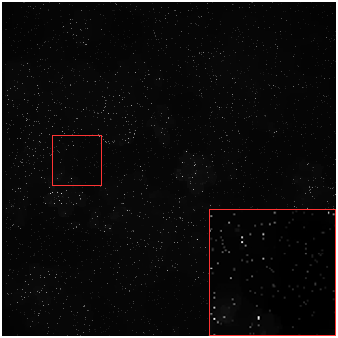}
			\includegraphics[width=1\linewidth]{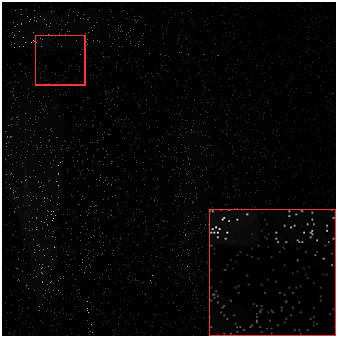}
			\includegraphics[width=1\linewidth]{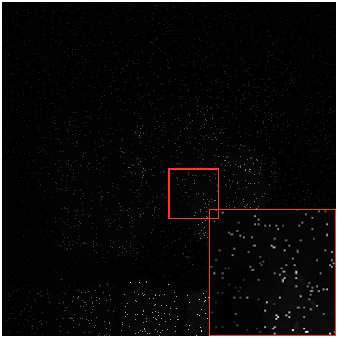}
			\includegraphics[width=1\linewidth]{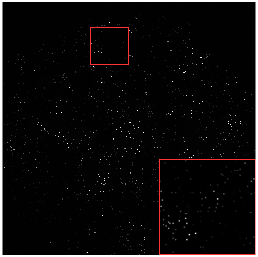}
			\includegraphics[width=1\linewidth]{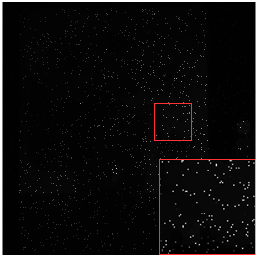}
	\end{minipage}}
	\subfigure[]{
		\begin{minipage}{0.056\textwidth}
			\includegraphics[width=1\linewidth]{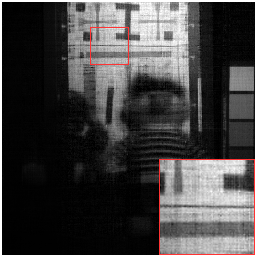}
			\includegraphics[width=1\linewidth]{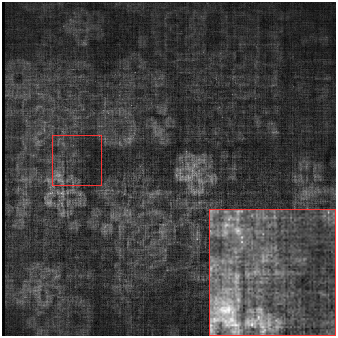}
			\includegraphics[width=1\linewidth]{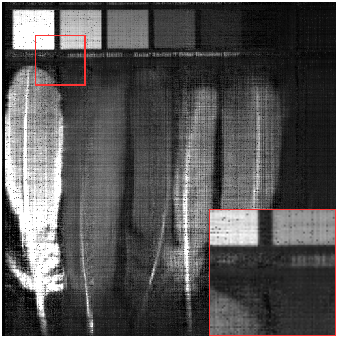}
			\includegraphics[width=1\linewidth]{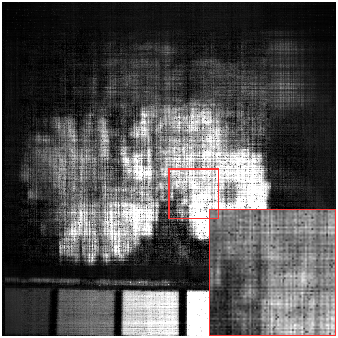}
			\includegraphics[width=1\linewidth]{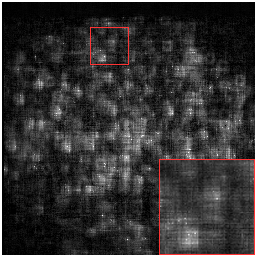}
			\includegraphics[width=1\linewidth]{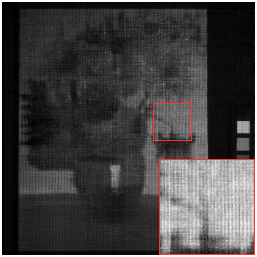}
	\end{minipage}}
	\subfigure[]{
		\begin{minipage}{0.056\textwidth}
			\includegraphics[width=1\linewidth]{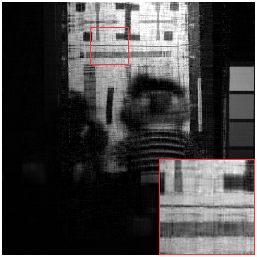}
			\includegraphics[width=1\linewidth]{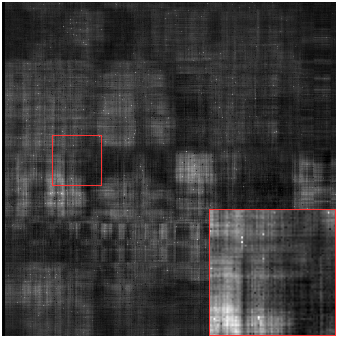}
			\includegraphics[width=1\linewidth]{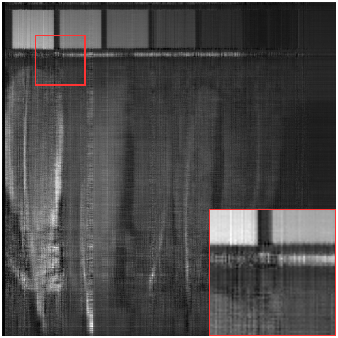}
			\includegraphics[width=1\linewidth]{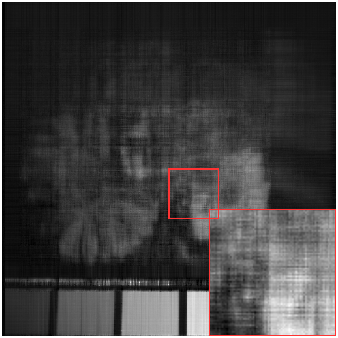}
			\includegraphics[width=1\linewidth]{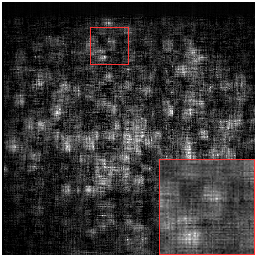}
			\includegraphics[width=1\linewidth]{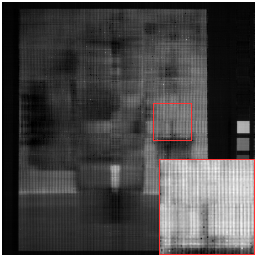}
	\end{minipage}}
	\subfigure[]{
		\begin{minipage}{0.056\textwidth}
			\includegraphics[width=1\linewidth]{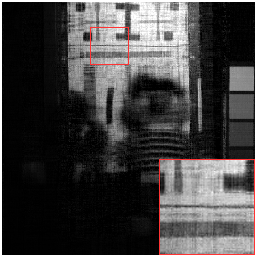}
			\includegraphics[width=1\linewidth]{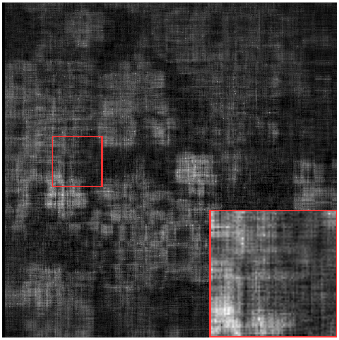}
			\includegraphics[width=1\linewidth]{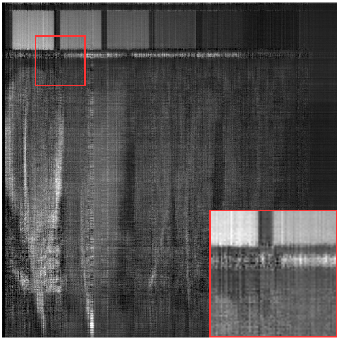}
			\includegraphics[width=1\linewidth]{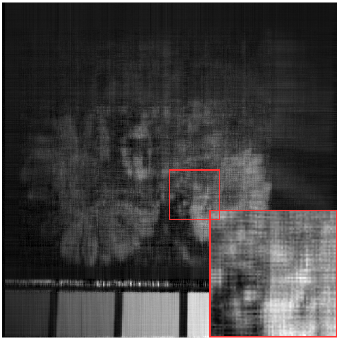}
			\includegraphics[width=1\linewidth]{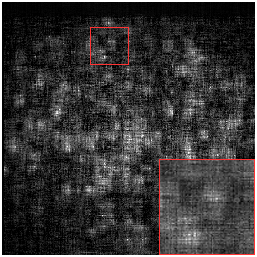}
			\includegraphics[width=1\linewidth]{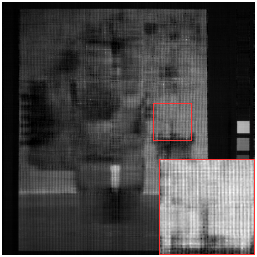}
	\end{minipage}}
	\subfigure[]{
		\begin{minipage}{0.056\textwidth}
			\includegraphics[width=1\linewidth]{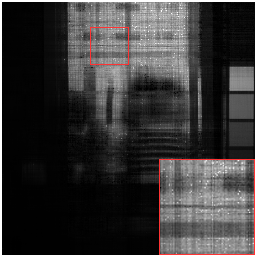}
			\includegraphics[width=1\linewidth]{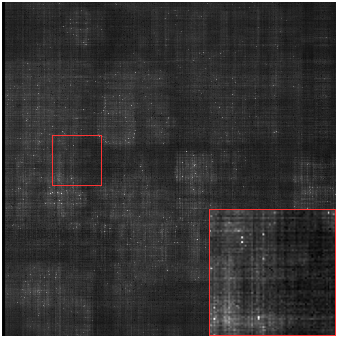}
			\includegraphics[width=1\linewidth]{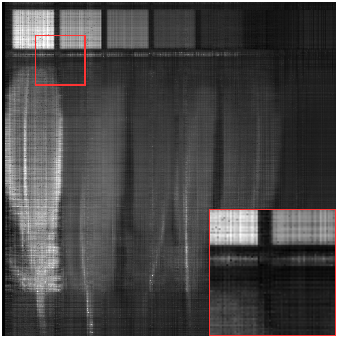}
			\includegraphics[width=1\linewidth]{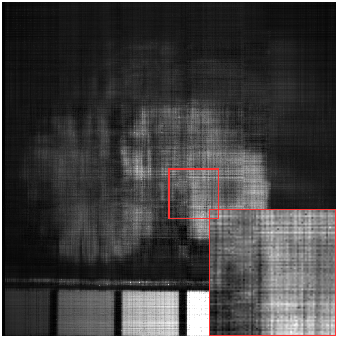}
			\includegraphics[width=1\linewidth]{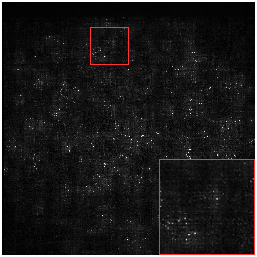}
			\includegraphics[width=1\linewidth]{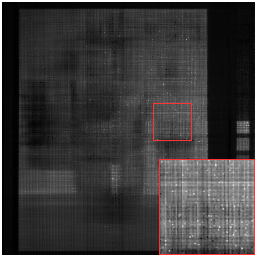}
	\end{minipage}}
	\subfigure[]{
		\begin{minipage}{0.056\textwidth}
			\includegraphics[width=1\linewidth]{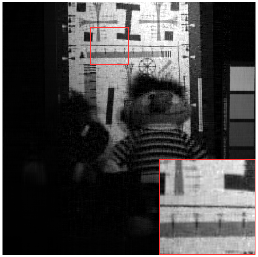}
			\includegraphics[width=1\linewidth]{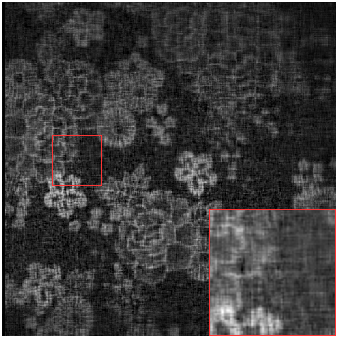}
			\includegraphics[width=1\linewidth]{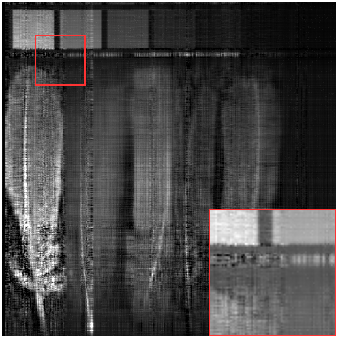}
			\includegraphics[width=1\linewidth]{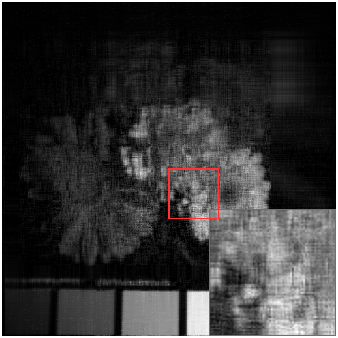}
			\includegraphics[width=1\linewidth]{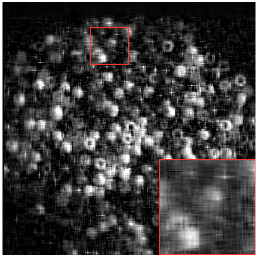}
			\includegraphics[width=1\linewidth]{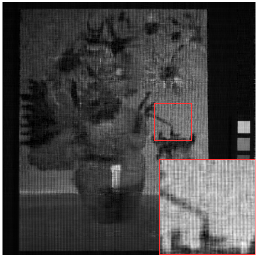}
	\end{minipage}}
	\subfigure[]{
		\begin{minipage}{0.056\textwidth}
			\includegraphics[width=1\linewidth]{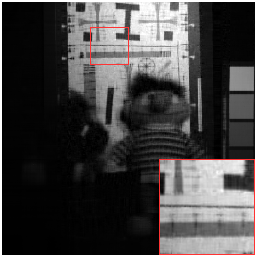}
			\includegraphics[width=1\linewidth]{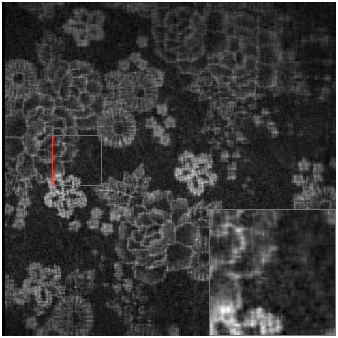}
			\includegraphics[width=1\linewidth]{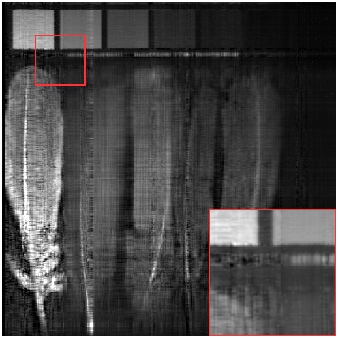}
			\includegraphics[width=1\linewidth]{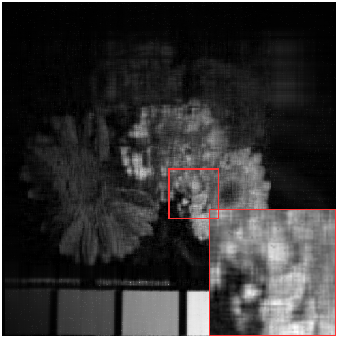}
			\includegraphics[width=1\linewidth]{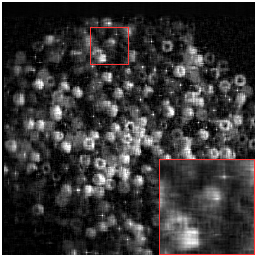}
			\includegraphics[width=1\linewidth]{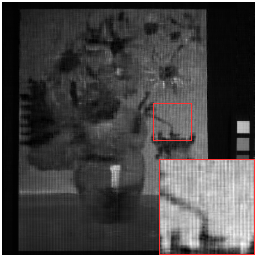}
	\end{minipage}}
	\caption{(a) Original images. (b) Corresponding sampled images with SR 2.5\%. (c)-(l) and (m) Completed images achieved by nine competing methods and proposed TLNM method and TLNMTV method, respectively. (a) Original image. (b) Obeserved image. (c) MC-ALM. (d) HaLRTC. (e) TMac. (f) LRTC-TV. (g) Trace/TV. (h) t-SVD. (i) McpTC. (j) ScadTC. (k) FTNN. (l) TLNM. (m) TLNMTV.}
	\label{HSITC025}
\end{figure*}

\begin{figure*}[!h]
	\centering
	\vspace{0cm}
	\subfigure[]{
		\includegraphics[width=1\linewidth, height=0.24\textheight]{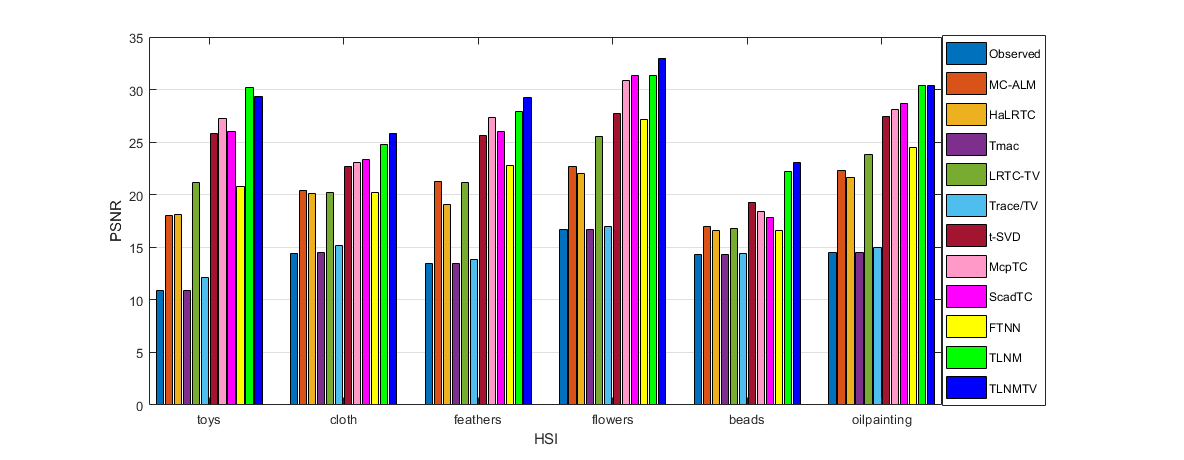}}
	\\
	\subfigure[]{
		\includegraphics[width=1\linewidth, height=0.24\textheight]{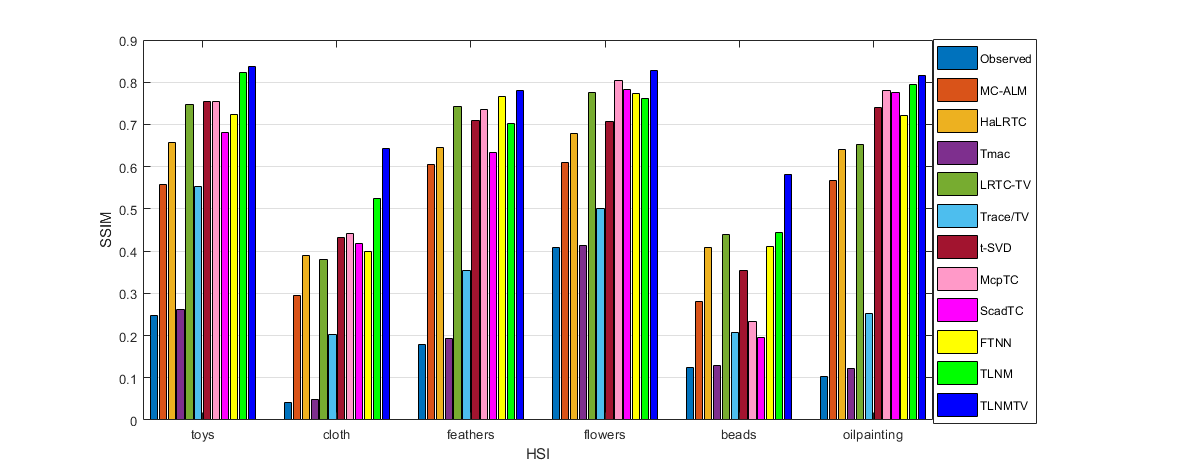}}
	\\
	\subfigure[]{
		\includegraphics[width=1\linewidth, height=0.24\textheight]{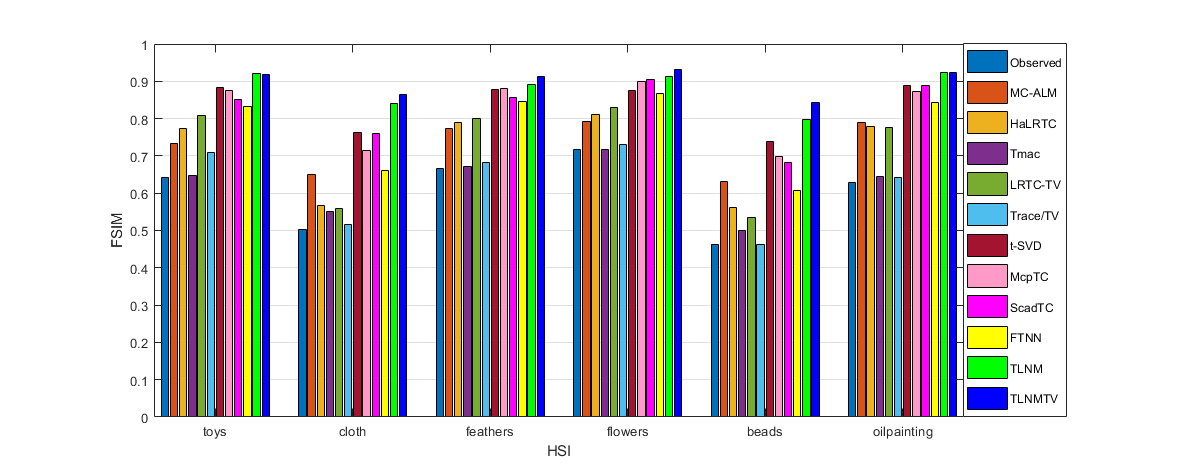}}
	\caption{The PSNR, SSIM, and FSIM of the results by different methods on all the HSI data with the sampling rate 2.5\%.}
	\label{HSI025}
\end{figure*}

\subsection{HSI Completion}
The HSI\footnote{http://www.cs.columbia.edu/CAVE/databases/multispectral/} test data is used in the experiment came from the open source CAVE data sets (respectively named "chart\_and\_stuffed\_toy", "cloth", "feathers", "flowers", "beads", "oil\_painting"). The size of HSI is $512 \times 512 \times 31$, indicating that the spatial resolution is $512 \times 512$ and the spectral resolution is 31, respectively. We show the visual effects of HSI recovery with sampling rates of 5\%, 10\% and 20\% respectively in Figs.\ref{HSITC5}-\ref{HSITC20}, in which the small figure is a local enlargement. It can be seen from the figure that the two proposed methods are superior to other methods, and TLNMTV method is more accurate on the details. Table \ref{THSI} lists the average quantitative numerical results of HSI at three sampling rates. It is not difficult to find that the proposed TLNM method is far superior to the classical HaLRTC method based on SNN and t-SVD method based on TNN. This shows that the proposed SNN method combined with QR decomposition and $L_{2,1}$ norm can improve the recovery effect. At 5\% sampling rate, the PSNR value of TLNMTV method is 3dB higher than that of suboptimal method. Even at 20\% sampling rate, the TLNMTV method is still 1dB higher than the suboptimal FTNN method.

Through the above experimental study, it is found that our method has a better effect at a lower sampling rate. In view of this, we have done experiments at a sampling rate of 2.5\%. Fig. \ref{HSITC025} shows the visual effect at a 2.5\% sampling rate, and Fig. \ref{HSI025} further illustrates the PSNR, SSIM, and FSIM of all methods on all HSIs. It can be seen from Fig. \ref{HSITC025} that our method can still obtain more significant recovery effect at a very low sampling rate of 2.5\%, which is quite effective in extracting both global visual and structural information and local detail information. And Fig. \ref{HSI025} shows that when the sampling rate is 2.5\%, the results obtained by our method are better than those obtained by other state-of-the-art methods.
\begin{figure*}[h] %??????????, ???????, ??????????????, ????
	\centering  %??????
	\vspace{0cm} %??????????
	\subfigure[]{
		\begin{minipage}{0.056\textwidth}
			\includegraphics[width=1\linewidth]{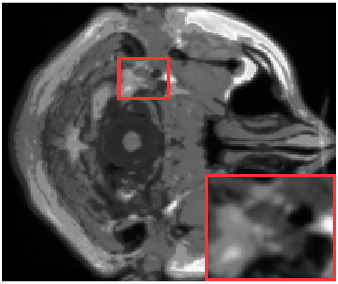}
			\includegraphics[width=1\linewidth]{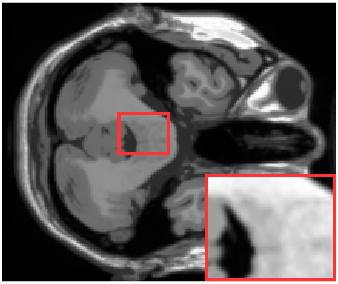}
			\includegraphics[width=1\linewidth]{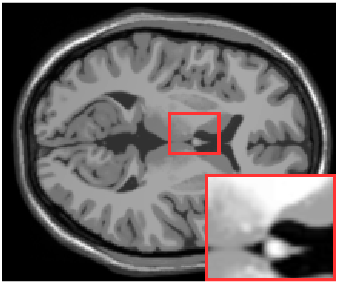}
			\includegraphics[width=1\linewidth]{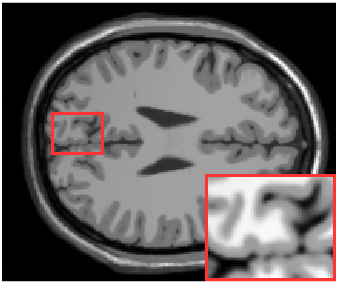}
	\end{minipage}}
	\subfigure[]{
		\begin{minipage}{0.056\textwidth}
			\includegraphics[width=1\linewidth]{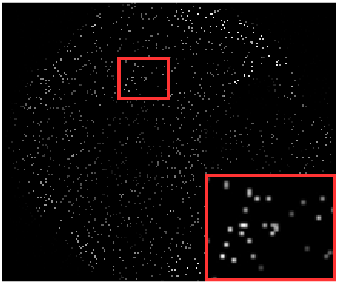}
			\includegraphics[width=1\linewidth]{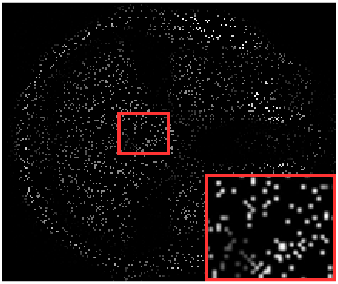}
			\includegraphics[width=1\linewidth]{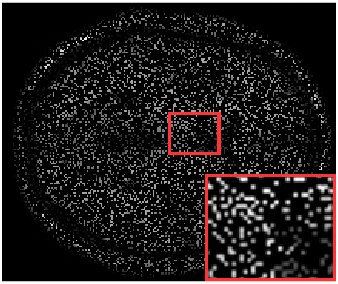}
			\includegraphics[width=1\linewidth]{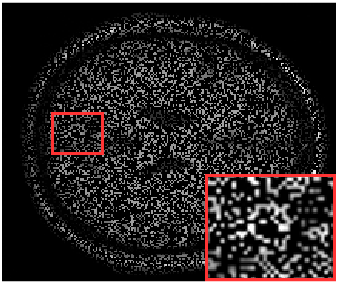}
	\end{minipage}}
	\subfigure[]{
		\begin{minipage}{0.056\textwidth}
			\includegraphics[width=1\linewidth]{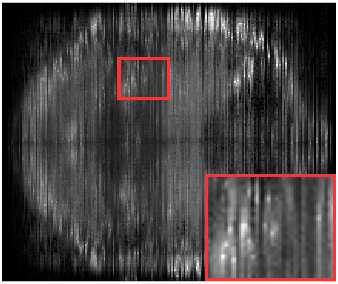}
			\includegraphics[width=1\linewidth]{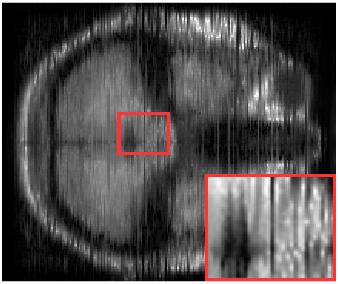}
			\includegraphics[width=1\linewidth]{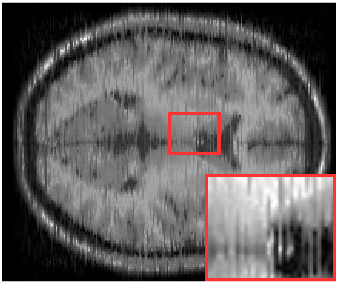}
			\includegraphics[width=1\linewidth]{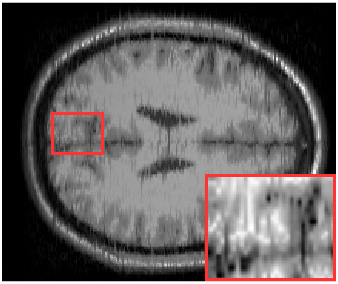}
	\end{minipage}}
	\subfigure[]{
		\begin{minipage}{0.056\textwidth}
			\includegraphics[width=1\linewidth]{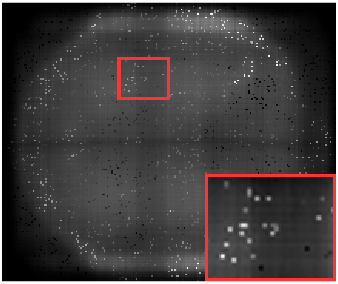}
			\includegraphics[width=1\linewidth]{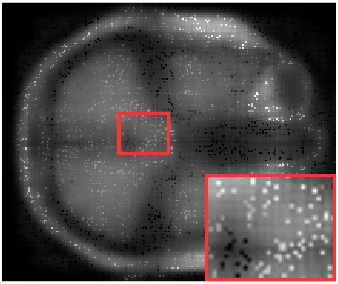}
			\includegraphics[width=1\linewidth]{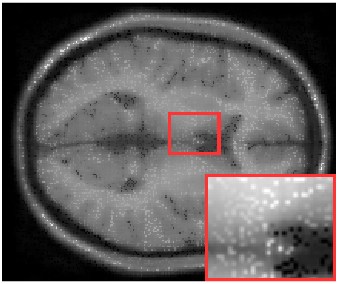}
			\includegraphics[width=1\linewidth]{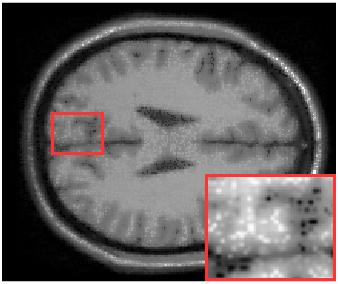}
	\end{minipage}}
	\subfigure[]{
		\begin{minipage}{0.056\textwidth}
			\includegraphics[width=1\linewidth]{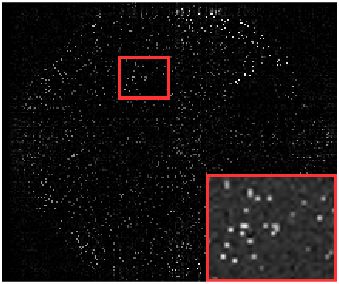}
			\includegraphics[width=1\linewidth]{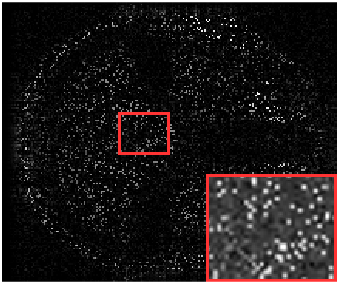}
			\includegraphics[width=1\linewidth]{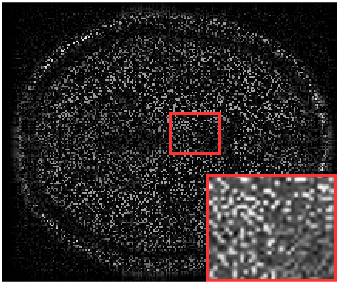}
			\includegraphics[width=1\linewidth]{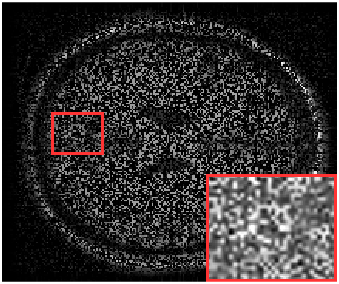}
	\end{minipage}}
	\subfigure[]{
		\begin{minipage}{0.056\textwidth}
			\includegraphics[width=1\linewidth]{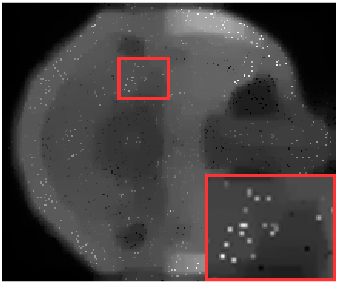}
			\includegraphics[width=1\linewidth]{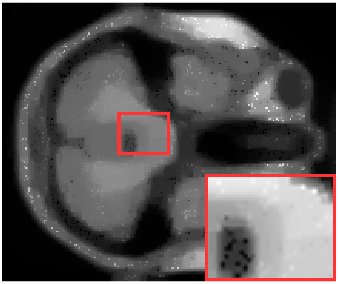}
			\includegraphics[width=1\linewidth]{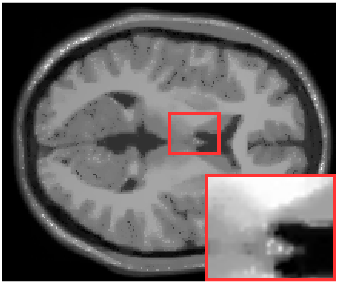}
			\includegraphics[width=1\linewidth]{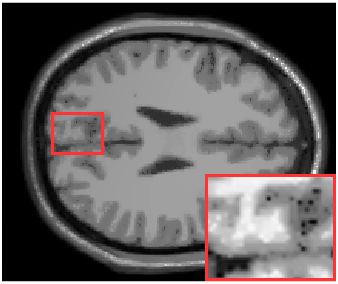}
	\end{minipage}}
	\subfigure[]{
		\begin{minipage}{0.056\textwidth}
			\includegraphics[width=1\linewidth]{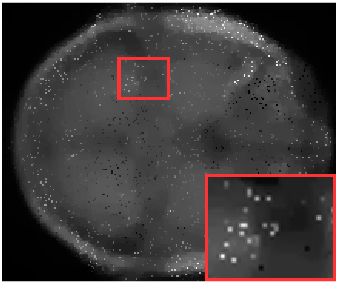}
			\includegraphics[width=1\linewidth]{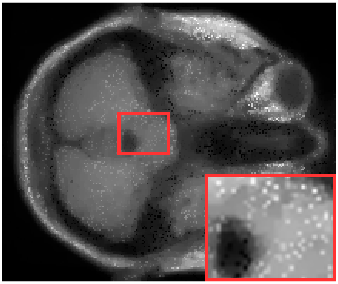}
			\includegraphics[width=1\linewidth]{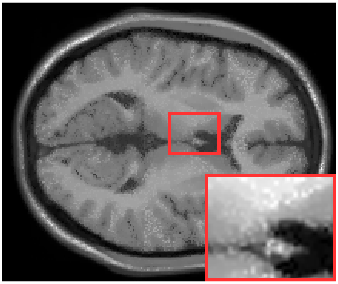}
			\includegraphics[width=1\linewidth]{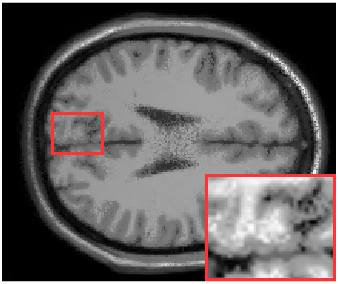}
	\end{minipage}}
	\subfigure[]{
		\begin{minipage}{0.056\textwidth}
			\includegraphics[width=1\linewidth]{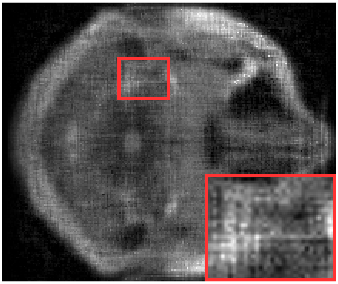}
			\includegraphics[width=1\linewidth]{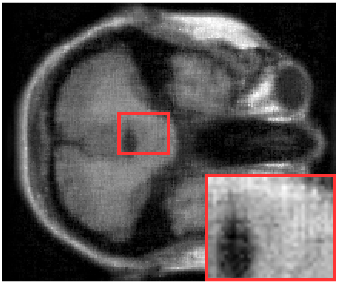}
			\includegraphics[width=1\linewidth]{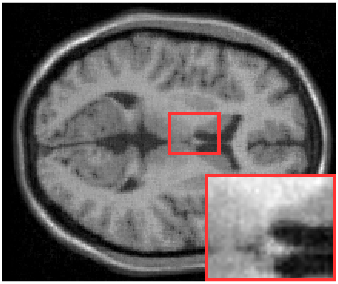}
			\includegraphics[width=1\linewidth]{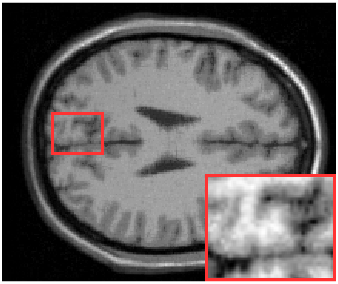}
	\end{minipage}}
	\subfigure[]{
		\begin{minipage}{0.056\textwidth}
			\includegraphics[width=1\linewidth]{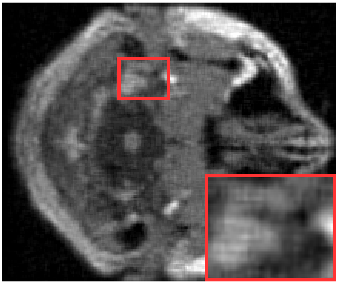}
			\includegraphics[width=1\linewidth]{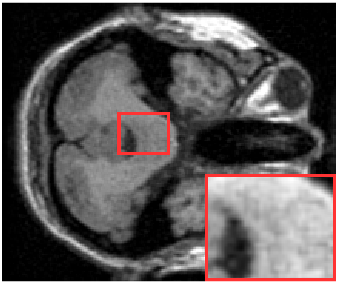}
			\includegraphics[width=1\linewidth]{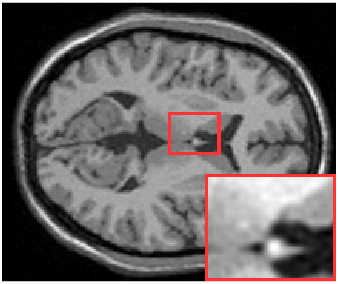}
			\includegraphics[width=1\linewidth]{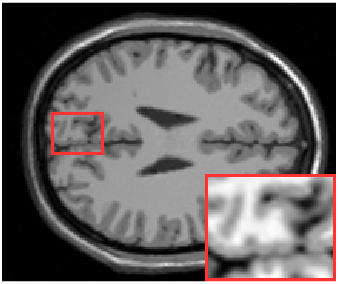}
	\end{minipage}}
	\subfigure[]{
		\begin{minipage}{0.056\textwidth}
			\includegraphics[width=1\linewidth]{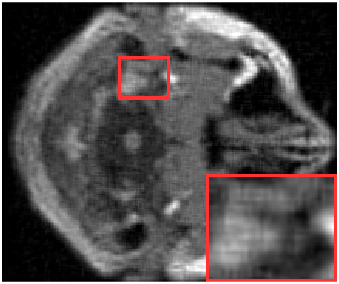}
			\includegraphics[width=1\linewidth]{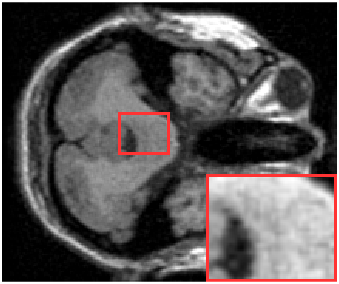}
			\includegraphics[width=1\linewidth]{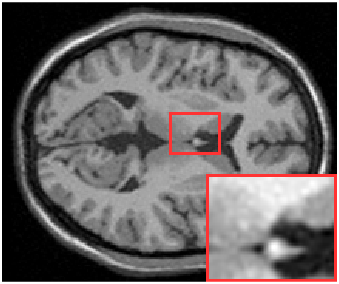}
			\includegraphics[width=1\linewidth]{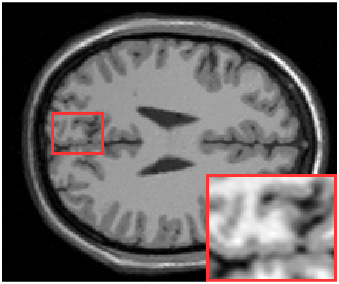}
	\end{minipage}}
	\subfigure[]{
		\begin{minipage}{0.056\textwidth}
			\includegraphics[width=1\linewidth]{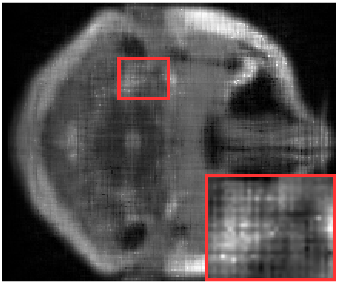}
			\includegraphics[width=1\linewidth]{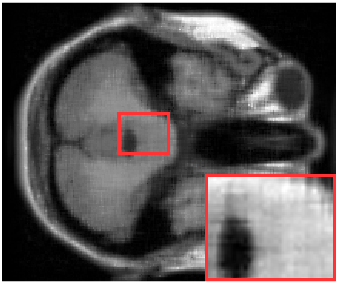}
			\includegraphics[width=1\linewidth]{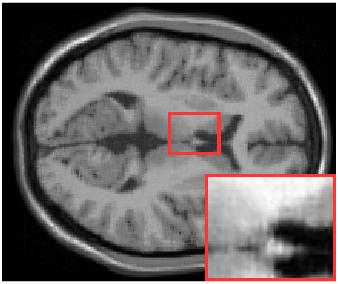}
			\includegraphics[width=1\linewidth]{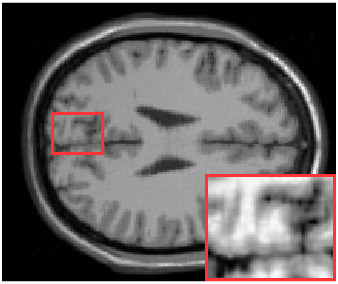}
	\end{minipage}}
	\subfigure[]{
		\begin{minipage}{0.056\textwidth}
			\includegraphics[width=1\linewidth]{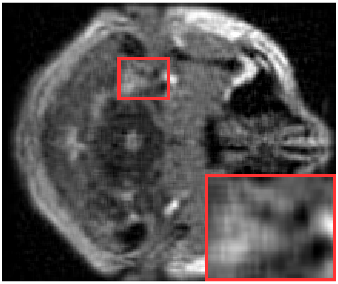}
			\includegraphics[width=1\linewidth]{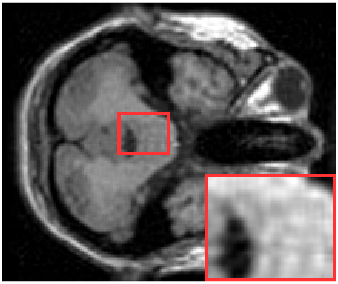}
			\includegraphics[width=1\linewidth]{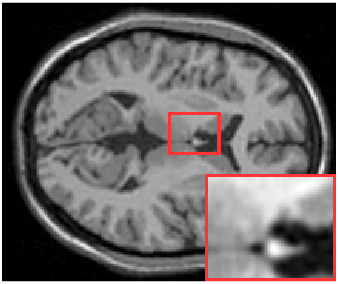}
			\includegraphics[width=1\linewidth]{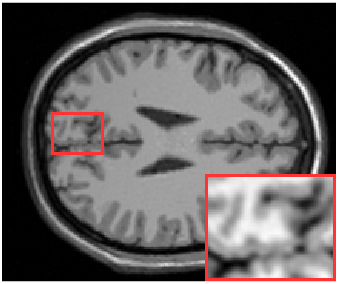}
	\end{minipage}}
	\subfigure[]{
		\begin{minipage}{0.056\textwidth}
			\includegraphics[width=1\linewidth]{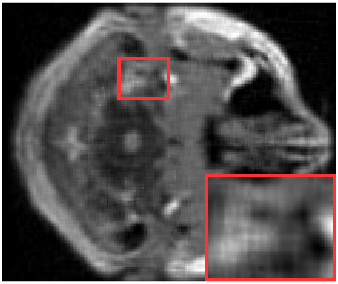}
			\includegraphics[width=1\linewidth]{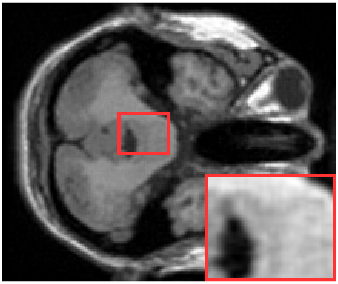}
			\includegraphics[width=1\linewidth]{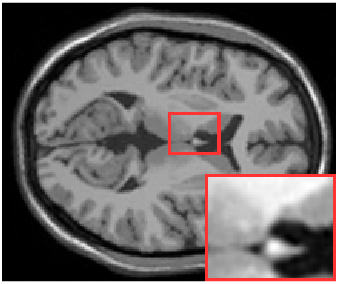}
			\includegraphics[width=1\linewidth]{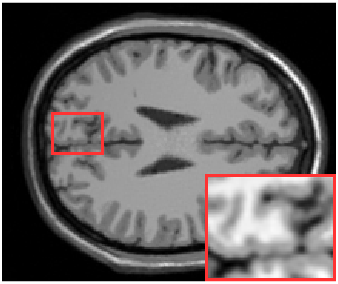}
	\end{minipage}}
	\caption{(a) Original images. (b) Corresponding sampled images with SR 5\% to 30\%. (c)-(l) and (m) Completed images achieved by nine competing methods and proposed TLNM method and TLNMTV method, respectively. (a) Original image. (b) Observed image. (c) MC-ALM. (d) HaLRTC. (e) TMac. (f) LRTC-TV. (g) Trace/TV. (h) t-SVD. (i) McpTC. (j) ScadTC. (k) FTNN. (l) TLNM. (m) TLNMTV.}
	\label{MRI}
\end{figure*}
\begin{table}[h]
	\caption{Quantitative evaluation of the results for MRI with different SRs.}
	\resizebox{\textwidth}{!}{
		\begin{tabular}{|c|c|c|c|c|c|c|c|c|c|c|c|c|c|}
			\hline
			SR                    & index & Observed & MC-ALM  & HaLRTC  & Tmac     & LRTC-TV & Trace-TV & t-SVD   & McpTC   & ScadTC  & FTNN    & TLNM    & TLNMTV  \\ \hline
			\multirow{4}{*}{5\%}  & PSNR  & 9.908    & 17.513  & 16.934  & 9.951    & 19.000  & 18.095   & 22.720  & 27.520  & 27.531  & 24.164  & 27.227  & 27.541  \\ \cline{2-14}
			& SSIM  & 0.174    & 0.289   & 0.299   & 0.103    & 0.530   & 0.464    & 0.514   & 0.748   & 0.748   & 0.653   & 0.724   & 0.782   \\ \cline{2-14}
			& FSIM  & 0.457    & 0.670   & 0.614   & 0.553    & 0.663   & 0.657    & 0.770   & 0.855   & 0.855   & 0.819   & 0.845   & 0.857   \\ \cline{2-14}
			& ERGAS & 1025.997 & 431.950 & 463.986 & 1021.063 & 363.138 & 414.979  & 242.057 & 135.859 & 135.675 & 203.091 & 141.034 & 135.572 \\ \hline
			\multirow{4}{*}{10\%} & PSNR  & 10.143   & 20.171  & 19.910  & 10.284   & 22.579  & 21.763   & 25.431  & 30.154  & 30.156  & 27.010  & 30.330  & 31.797  \\ \cline{2-14}
			& SSIM  & 0.188    & 0.452   & 0.448   & 0.093    & 0.712   & 0.637    & 0.656   & 0.830   & 0.830   & 0.780   & 0.846   & 0.914   \\ \cline{2-14}
			& FSIM  & 0.492    & 0.746   & 0.717   & 0.564    & 0.787   & 0.777    & 0.830   & 0.895   & 0.895   & 0.874   & 0.897   & 0.927   \\ \cline{2-14}
			& ERGAS & 998.616  & 317.697 & 328.310 & 982.744  & 240.306 & 270.797  & 178.582 & 100.559 & 100.528 & 146.333 & 99.056  & 83.339  \\ \hline
			\multirow{4}{*}{20\%} & PSNR  & 10.656   & 23.458  & 24.151  & 11.143   & 27.659  & 26.219   & 29.224  & 34.587  & 34.448  & 30.880  & 34.130  & 35.778  \\ \cline{2-14}
			& SSIM  & 0.219    & 0.640   & 0.671   & 0.112    & 0.874   & 0.814    & 0.809   & 0.937   & 0.932   & 0.892   & 0.936   & 0.965   \\ \cline{2-14}
			& FSIM  & 0.544    & 0.824   & 0.828   & 0.569    & 0.899   & 0.882    & 0.899   & 0.953   & 0.950   & 0.930   & 0.946   & 0.967   \\ \cline{2-14}
			& ERGAS & 941.393  & 216.791 & 200.759 & 890.839  & 134.009 & 162.074  & 116.642 & 60.244  & 61.197  & 93.371  & 63.738  & 52.713  \\ \hline
			\multirow{4}{*}{30\%} & PSNR  & 11.231   & 25.921  & 27.584  & 11.950   & 30.867  & 29.330   & 32.172  & 37.987  & 37.206  & 33.712  & 37.591  & 38.364  \\ \cline{2-14}
			& SSIM  & 0.253    & 0.760   & 0.814   & 0.144    & 0.932   & 0.895    & 0.888   & 0.977   & 0.966   & 0.940   & 0.972   & 0.980   \\ \cline{2-14}
			& FSIM  & 0.575    & 0.875   & 0.895   & 0.580    & 0.944   & 0.930    & 0.938   & 0.979   & 0.974   & 0.958   & 0.975   & 0.981   \\ \cline{2-14}
			& ERGAS & 881.074  & 163.086 & 134.957 & 812.918  & 92.554  & 113.377  & 83.475  & 40.760  & 44.527  & 67.186  & 42.717  & 39.060  \\ \hline
	\end{tabular}}
	\label{TMRI}
\end{table}
\begin{figure}[h]
	\vspace{0cm}
	\subfigure[]{
		\includegraphics[width=0.47\linewidth,height=0.165\textheight]{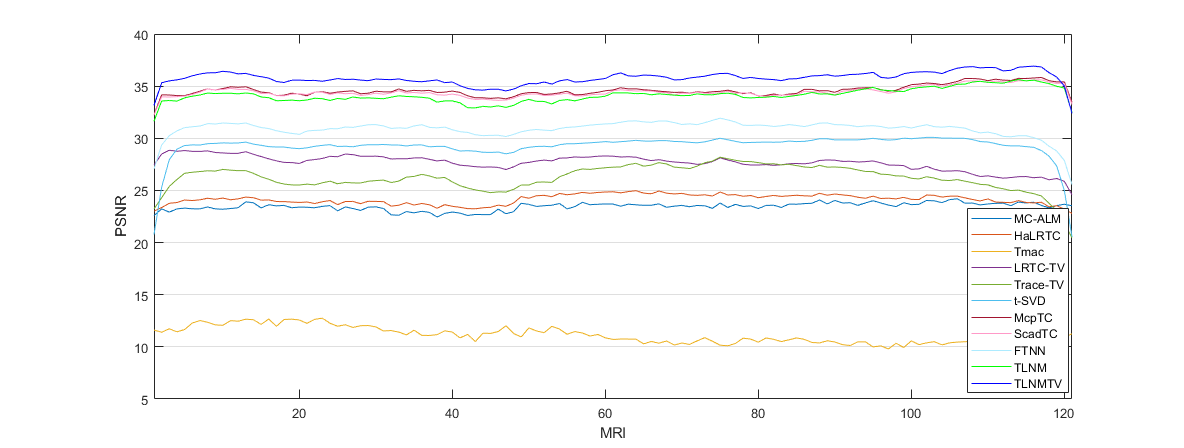}}
	\subfigure[]{
		\includegraphics[width=0.47\linewidth,height=0.165\textheight]{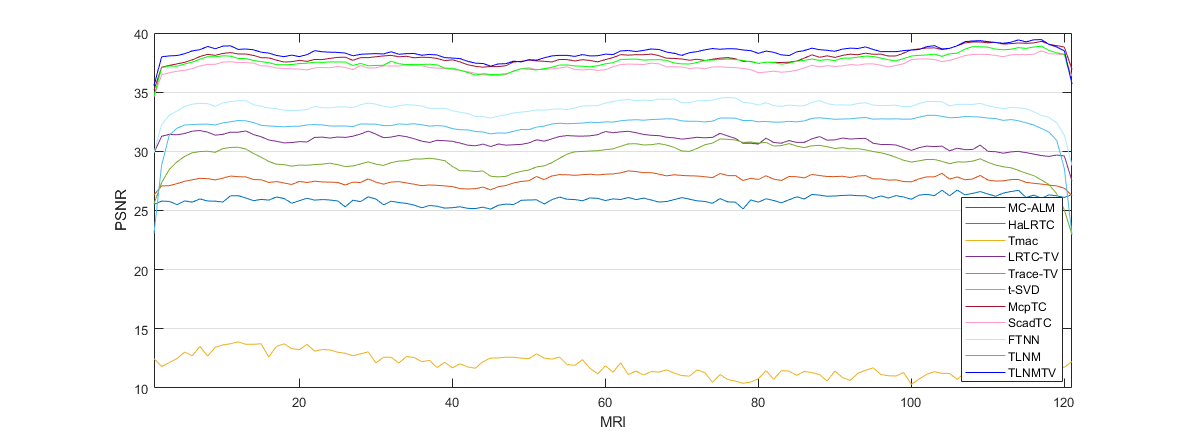}}
	\\
	\subfigure[]{
		\includegraphics[width=0.47\linewidth,height=0.165\textheight]{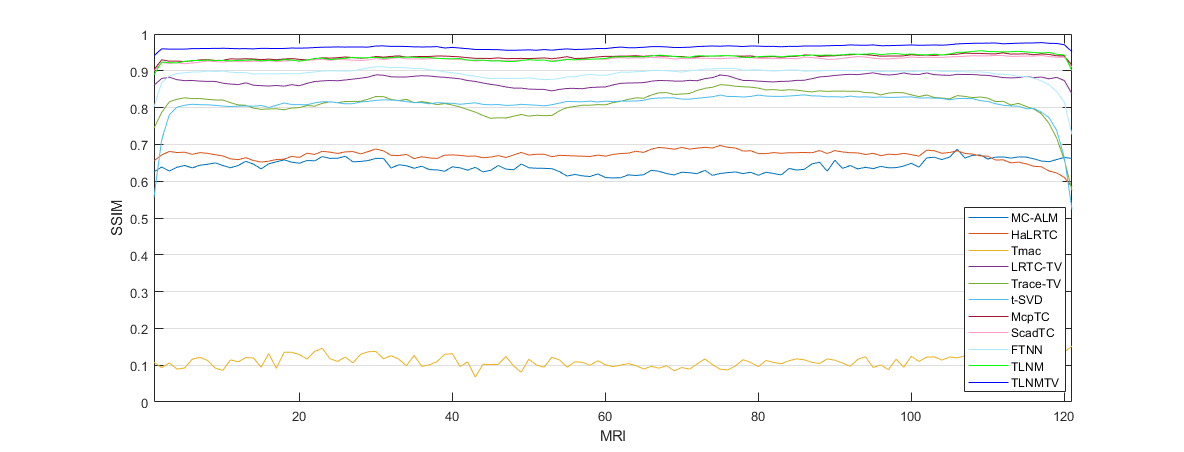}}
	\subfigure[]{
		\includegraphics[width=0.47\linewidth,height=0.165\textheight]{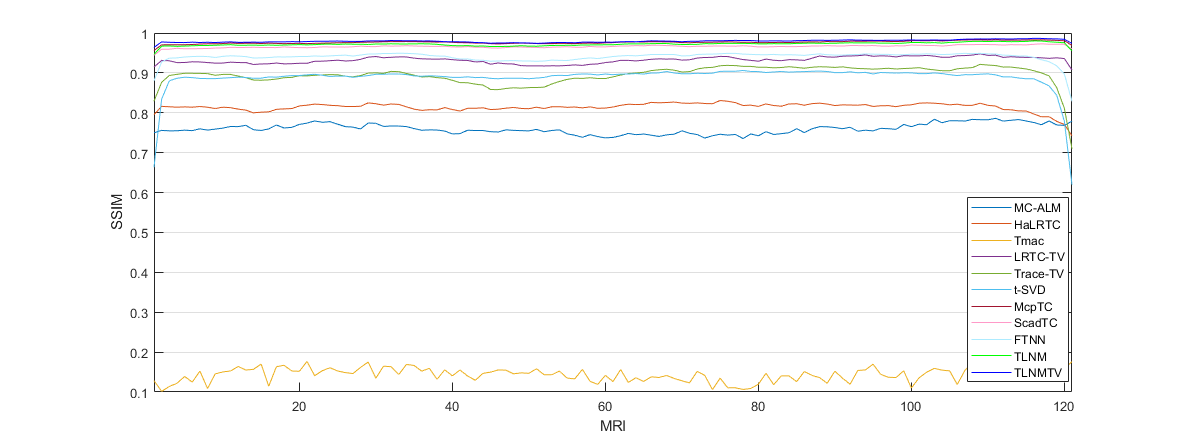}}
	\\
	\subfigure[]{
		\includegraphics[width=0.47\linewidth,height=0.165\textheight]{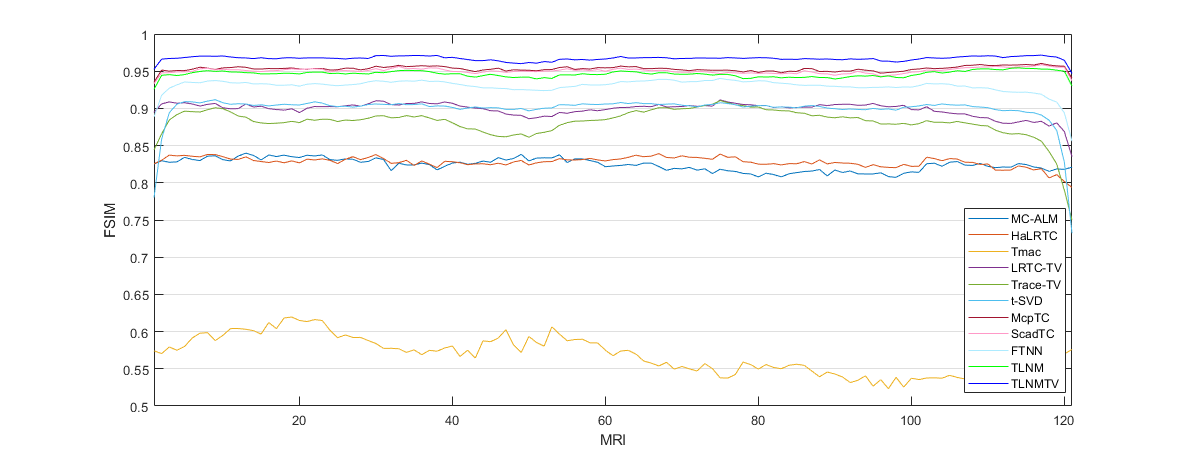}}
	\subfigure[]{
		\includegraphics[width=0.47\linewidth,height=0.165\textheight]{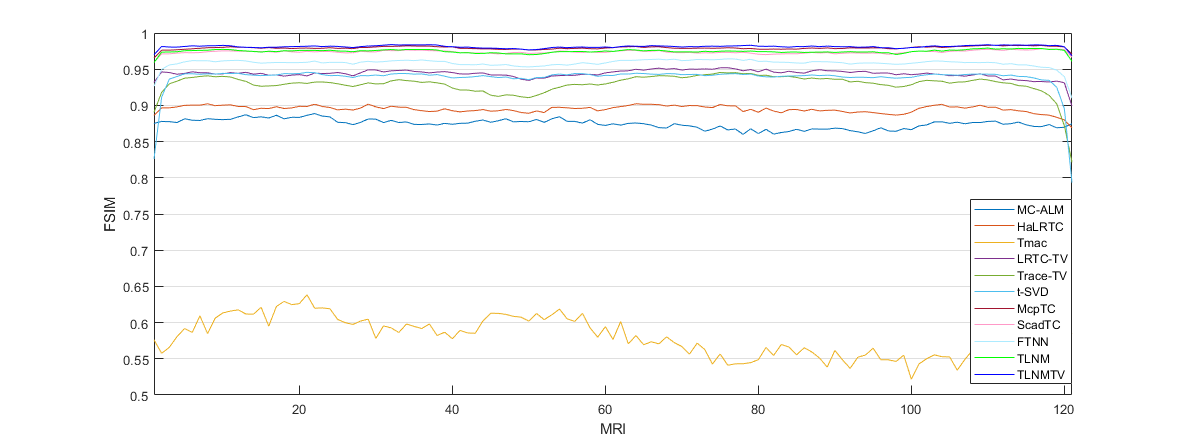}}
	\caption{ PSNR, SSIM, and FSIM values comparison of different methods for each slice on MRI data sets, (a), (c), (e) under SR 20\%, (b), (d), (f) under SR 30\%.}
	\label{MRI0302}
\end{figure}
\subsection{MRI Completion}
We test the performance of the proposed method and the comparison method on MRI\footnote{http://brainweb.bic.mni.mcgill.ca/brainweb/selection\_normal.html} data with the size of $181\times217\times121$. Fig. \ref{MRI} shows the visual effects of MRI images at the sampling rate of 5\%, 10\%, 20\% and 30\% respectively. It can be found from the figure that for the visual effect of MRI image recovery, our method is superior to the comparison method in both global and local information recovery. Further, Fig. \ref{MRI0302} describes the specific values of PSNR, SSIM and FSIM of each slice at the sampling rate of 20\% and 30\%. As one can see from Fig. \ref{MRI0302} that the TLNMTV method is the optimal method. Table \ref{TMRI} reports the PSNR, SSIM, FSIM, ERGAS values for each sample rate. It can be seen from the table that at different sampling rates, the TLNM method is obviously superior to the HaLRTC and t-SVD, and the results obtained are very close to those obtained by non-convex method McpTC. Both methods are based on SNN, which indicates that the convex TLNM method is very efficient. In addition, the TLNMTV method is obviously superior to the suboptimal method McpTC. Even when the sampling rate is 20\%, the PSNR value of TLNMTV method is 1.2 dB higher than that of McpTC.

\section{Parameters setting and Convergency Behaviours}
\begin{figure*}[h]
	\centering
	\subfigure[HSI]{
		%		\label{level.sub.2}
		\includegraphics[width=0.8\linewidth]{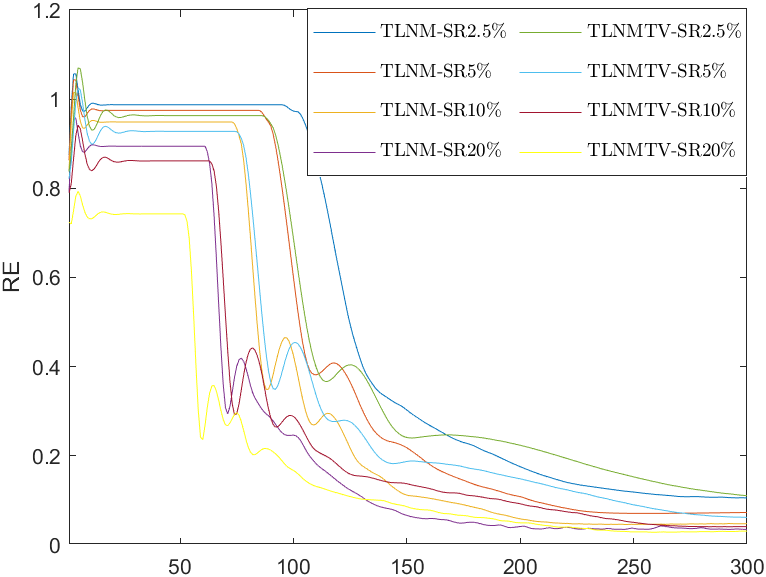}
	}
	\\
	\subfigure[MRI]{
		%		\label{level.sub.2}
		\includegraphics[width=0.8\linewidth]{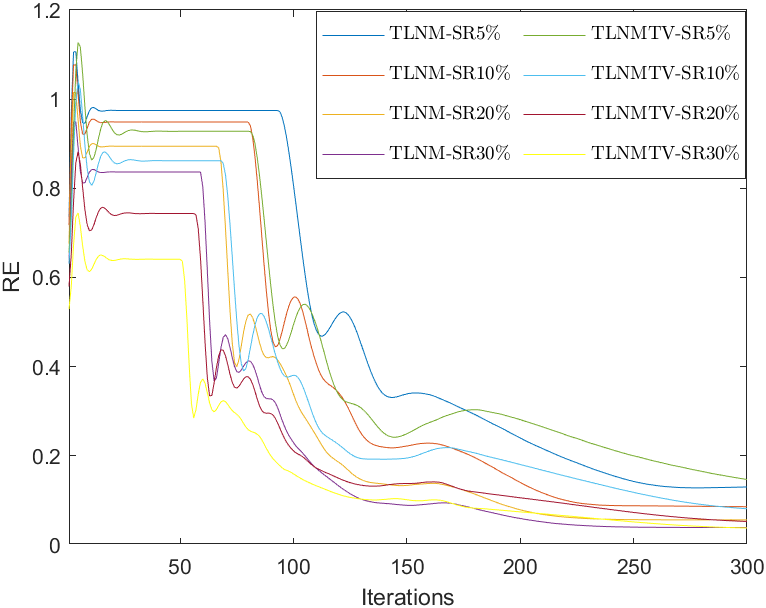}	
	}
	\caption{The convergence behaviours of LRTC Algorithm, with respect to different sampling rates.}
	\label{convergence}
\end{figure*}

We obtain the optimal solutions of the two algorithm parameters through experimental tests. Table \ref{rnset} lists the optimal parameters of $r_{n}$ for different images with different sampling rates. We let $\alpha$ and $\beta=[1/N, 1/N, \dots, 1/N]$. In order to observe the convergence behavior of the algorithm, the relative error (RE) of the recovery tensor is defined as follows:
\begin{eqnarray}
RE:=\frac{\|\mathcal{X}_{rec}-\mathcal{X}_{o}\|_{F}}{\|\mathcal{X}_{o}\|_{F}},
\end{eqnarray}
where $\mathcal{X}_{rec}$ is recovery tensor and $\mathcal{X}_{o}$ is origin tensor. Fig. \ref{convergence} shows the iterative RE of the proposed algorithm for HSI and MRI. As the iteration progresses, the RE can become smaller, which guarantees the convergence of the proposed algorithm.
\begin{table}[h]
	\caption{The value of $r_{n}$ is set for different images wiyh different SRs.}
	\resizebox{\textwidth}{!}{
		\begin{tabular}{|c|cccc|cccc|}
			\hline
			Method                & \multicolumn{4}{c|}{TLNM}                                                                                   & \multicolumn{4}{c|}{TLNMTV}                                                                                    \\ \hline
			SR                    & \multicolumn{1}{c|}{2.50\%}   & \multicolumn{1}{c|}{5\%}       & \multicolumn{1}{c|}{10\%}      & 20\%      & \multicolumn{1}{c|}{2.50\%}    & \multicolumn{1}{c|}{5\%}       & \multicolumn{1}{c|}{10\%}       & 20\%       \\ \hline
			chart and stuffed toy & \multicolumn{1}{c|}{50 50 4}  & \multicolumn{1}{c|}{85 85 4}   & \multicolumn{1}{c|}{125 125 4} & 205 205 4 & \multicolumn{1}{c|}{60 60 4}   & \multicolumn{1}{c|}{95 95 4}   & \multicolumn{1}{c|}{145 145 4}  & 205 205 4  \\ \hline
			cloth                 & \multicolumn{1}{c|}{50 50 4}  & \multicolumn{1}{c|}{85 85 4}   & \multicolumn{1}{c|}{125 125 4} & 205 205 4 & \multicolumn{1}{c|}{130 130 4} & \multicolumn{1}{c|}{160 160 4} & \multicolumn{1}{c|}{235 235 4}  & 275 275 4  \\ \hline
			feather               & \multicolumn{1}{c|}{40 40 4}  & \multicolumn{1}{c|}{65 65 4}   & \multicolumn{1}{c|}{105 105 4} & 205 205 4 & \multicolumn{1}{c|}{50 50 4}   & \multicolumn{1}{c|}{95 95 4}   & \multicolumn{1}{c|}{145 145 4}  & 205 205 4  \\ \hline
			flower                & \multicolumn{1}{c|}{50 50 4}  & \multicolumn{1}{c|}{75 75 4}   & \multicolumn{1}{c|}{125 125 4} & 205 205 4 & \multicolumn{1}{c|}{50 50 4}   & \multicolumn{1}{c|}{95 95 4}   & \multicolumn{1}{c|}{165 165 4}  & 205 205 4  \\ \hline
			beads                 & \multicolumn{1}{c|}{50 50 4}  & \multicolumn{1}{c|}{115 115 4} & \multicolumn{1}{c|}{175 175 4} & 255 255 4 & \multicolumn{1}{c|}{110 110 4} & \multicolumn{1}{c|}{155 155 4} & \multicolumn{1}{c|}{215 215 4}  & 235 235 4  \\ \hline
			oil painting          & \multicolumn{1}{c|}{60 60 4}  & \multicolumn{1}{c|}{95 95 4}   & \multicolumn{1}{c|}{185 185 4} & 235 235 4 & \multicolumn{1}{c|}{80 80 4}   & \multicolumn{1}{c|}{135 135 4} & \multicolumn{1}{c|}{205 205 4}  & 250 250 4  \\ \hline
			SR                    & \multicolumn{1}{c|}{5\%}      & \multicolumn{1}{c|}{10\%}      & \multicolumn{1}{c|}{20\%}      & 30\%      & \multicolumn{1}{c|}{5\%}       & \multicolumn{1}{c|}{10\%}      & \multicolumn{1}{c|}{20\%}       & 30\%       \\ \hline
			MRI                   & \multicolumn{1}{c|}{50 60 30} & \multicolumn{1}{c|}{60 70 40}  & \multicolumn{1}{c|}{70 90 55}  & 90 110 60 & \multicolumn{1}{c|}{50 60 30}  & \multicolumn{1}{c|}{80 90 40}  & \multicolumn{1}{c|}{100 120 55} & 110 130 60 \\ \hline
	\end{tabular}}
	\label{rnset}
\end{table}

\section{Conclusion}

This paper integrates SNN with QR decomposition and $L_{2,1}$ norm to solve the LRTC problem, and proposes the TLNM method. The TV regularization term is further introduced to improve the local prior information and make the recovery effect of the TLNM method more ideal. HSI and MRI are tested at different sampling rates. Experiments show that our TLNM method is more efficient than the TNN-based TLNM-TQR method and outperforms the classic HaLRTC and t-SVD methods. Besides, the TLMNTV method surpasses other compared methods in this paper, thus confirming the high efficiency of our proposed method. It is particularly worth emphasizing that our method still achieves fairly good HSI recovery when the SR is extremely low, i.e. 2.5\%. In future work, we will try to introduce other more efficient forms, such as non-convex transformations to replace the nuclear norm, to further improve the recovery performance of our method.

%In this paper, SNN is accommodated with QR decomposition and $L_{2,1}$ norm to solve LRTC problem, and the TLNM method is proposed. The TV regularization term is further introduced to improve the local prior information and make the recovery effect more accurate. HSI and MRI are tested at various sampling rates. Experiments show that our TLNM method outperforms the classical HaLRTC and t-SVD methods, which is consistent with our expectation that it is more effective than the TLNM-TQR method based on TNN method. The TLNMTV method is superior to the comparison method in this paper, which proves the high efficiency of our method. In particular, it is worth emphasizing that our method can still achieve good HSI recovery when SR is 2.5\%, which is extremely low. In future work, we will try to introduce other forms of substitution, such as non-convex transformation to replace the nuclear norm, to further improve the recovery effect of our method.

\bigbreak

\bibliography{bibtex/bare_jrnl_cs}

\end{document}